# GEOMETRIC CHARACTERIZATION OF INTERMITTENCY IN THE PARABOLIC ANDERSON MODEL

By Jürgen Gärtner,[1] Wolfgang König[1,2] and Stanislav Molchanov

*Technische Universität Berlin, Universität Leipzig and University of North Carolina*

We consider the parabolic Anderson problem $\partial_t u = \Delta u + \xi(x)u$ on $\mathbb{R}_+ \times \mathbb{Z}^d$ with localized initial condition $u(0,x) = \delta_0(x)$ and random i.i.d. potential $\xi$. Under the assumption that the distribution of $\xi(0)$ has a double-exponential, or slightly heavier, tail, we prove the following geometric characterization of intermittency: with probability one, as $t \to \infty$, the overwhelming contribution to the total mass $\sum_x u(t,x)$ comes from a slowly increasing number of "islands" which are located far from each other. These "islands" are local regions of those high exceedances of the field $\xi$ in a box of side length $2t\log^2 t$ for which the (local) principal Dirichlet eigenvalue of the random operator $\Delta + \xi$ is close to the top of the spectrum in the box. We also prove that the shape of $\xi$ in these regions is nonrandom and that $u(t,\cdot)$ is close to the corresponding positive eigenfunction. This is the geometric picture suggested by localization theory for the Anderson Hamiltonian.

## Contents

Received October 2005; revised May 2006.

[1]Supported in part by the German Science Foundation, Schwerpunkt Project SPP 1033.

[2]Supported in part by a Heisenberg grant awarded by the German Science Foundation, realized in 2003/04.

## 1. Introduction and main result.

1.1. *Main objective.* We consider the following Cauchy problem with localized initial datum:

$$\begin{aligned} \partial_t u(t,x) &= \Delta u(t,x) + \xi(x)u(t,x), && (t,x) \in (0,\infty)\times\mathbb{Z}^d, \\ u(0,x) &= \delta_0(x), && x \in \mathbb{Z}^d, \end{aligned} \tag{1.1}$$

where $\Delta$ is the lattice Laplacian on $\mathbb{Z}^d$, $\Delta f(x) = \sum_{y\sim x}[f(y)-f(x)]$, and $\xi = (\xi(x))_{x\in\mathbb{Z}^d}$ is a random i.i.d. potential modeling the environment. Let $\mathrm{Prob}(\cdot)$ and $\langle\cdot\rangle$ denote the underlying probability measure and expectation, respectively. Problem (1.1) is often referred to as the *parabolic Anderson problem* and appears in the description of population dynamics, catalytic reactions, and so forth. One interpretation of this model is in terms of branching random walk on $\mathbb{Z}^d$ with constant branching rate and spatially dependent branching mechanism governed by a typical realization of a homogeneous random field. If $\xi(z)$ denotes the mean offspring of a particle at site $z$, then the solution $u(t,x)$ to (1.1) describes the expected number of particles at $x$ at time $t$ for the given realization of the random field $\xi(\cdot)$. The references [2, 7, 12] provide more explanations and heuristics around the parabolic Anderson model. Recent related results are reviewed in [5].

The purpose of the present paper is to describe in detail the almost sure spatial structure of the solution $u(t,\cdot)$ for large $t$ for a large class of i.i.d. potentials which are unbounded from above. This is an attempt toward a mathematical foundation of the following descriptive manifestation of *intermittency*: *With probability one, as $t\to\infty$, the random field $u(t,\cdot)$ develops high peaks on islands which are located far from each other and which give the overwhelming contribution to the total mass*

$$U(t) = \sum_{x\in\mathbb{Z}^d} u(t,x). \tag{1.2}$$

Such a picture is suggested by localization theory for the Anderson Hamiltonian $\mathcal{H} = \Delta + \xi$. Indeed, since the (upper part of the) $\ell^2$-spectrum of $\mathcal{H}$ is a pure point spectrum, one may expand $u(t,\cdot)$ in a Fourier series with respect to the (random) eigenvalues $\lambda_k$ and the corresponding (random) eigenfunctions $\mathrm{e}_k$ of $\mathcal{H}$:

$$u(t,\cdot) = \sum_k e^{\lambda_k t}\mathrm{e}_k(0)\mathrm{e}_k(\cdot). \tag{1.3}$$

The exponentially localized eigenfunctions $\mathrm{e}_k(\cdot)$ corresponding to large eigenvalues $\lambda_k$ are expected to be sparsely distributed in space. Hence, for large $t$, the weighted superposition of only a few of them will contribute to the total mass of the solution. In particular, relevant for the total mass $U(t)$ are those high peaks of the solution $u(t,\cdot)$ that are caused by those local exceedances of the potential $\xi(\cdot)$ that nearly maximize the corresponding principal Dirichlet eigenvalue of the Hamiltonian $\mathcal{H}$. Moreover, locally the shape of $u(t,\cdot)$ should resemble the corresponding positive eigenfunction. It is this *geometric characterization of intermittency* which we make precise in the present paper. For directly analyzing (1.3), it would be necessary to control $\mathrm{e}_k(0)$. This task, however, appears difficult, since, for $\lambda_k$ large, the localization center of $\mathrm{e}_k$ is far away from the origin, and $\mathrm{e}_k(0)$ is relatively small and may even be negative. This difficulty is due to the localized initial datum. Hence, a direct approximation of the eigenvalues $\lambda_k$ and the eigenfunctions $\mathrm{e}_k$ via a resolvent cluster expansion of the same kind as that used in localization theory does not seem to be promising here. (Nevertheless, it may be successful for homogeneous initial data.) Instead, we control the relative total mass of the solution $u(t,\cdot)$ far away from the localization centers in terms of suitable *positive* eigenfunctions, compare Theorem 4.1 below. Each of them corresponds to the potential peak in one of these centers after removal of all the others. The exponential decay of these eigenfunctions is then controlled by a probabilistic cluster expansion, that is, by a decomposition of the paths in the Feynman–Kac representation of the solution into segments between successive visits to clusters of high potential peaks; compare Proposition 6.1 below.

Intermittency is often studied by comparing the large-$t$ asymptotics of the moments of $U(t)$ of successive orders. The heuristic relation to the above geometric picture of intermittency is explained in [7], for example. For various types of potentials, the second-order asymptotics, both of the moments and in an almost sure sense, of $U(t)$ have been investigated in several papers. In all these cases, these asymptotics are described in terms of a variational problem, which gives some insight in the geometry of the high peaks both of the potential and of the solution. For a large class of i.i.d. potentials bounded from above, this has been carried out in [1], and for certain unbounded from above i.i.d. fields with double-exponential tails and with heavier tails in [7, 8]. A fourth class of i.i.d. fields has been considered in [10]. There it is also shown that, under a natural regularity condition, these are *all* the four universality classes that can arise for i.i.d. fields. Correlated fields [9] and Gaussian and (high-peak) Poisson fields in the spatially continuous model have been looked at in [4] and [6].

One important model, which has been analyzed in greater detail, is Brownian motion among Poisson traps. In this model the random potential $\xi$ is equal to $-\infty$ in neighborhoods of the points of a homogeneous Poisson point process on $\mathbb{R}^d$ and 0 elsewhere. The solution $u(t,x)$ of the spatially continuous analog of problem (1.1) then equals the survival probability of a Brownian particle among those traps when moving from $x$ to 0 in time $t$. This model was studied in a series of papers by Sznitman (see his monograph [13] and the references therein). In particular, Sznitman proved the so-called *pinning effect*: the optimal survival strategy of the Brownian particle is to move in short time to one of the trap-free regions of optimal size and shape and to stay there for the remaining time. This statement may be considered as an alternative interpretation of intermittency for this type of potential.

In the present paper we continue the study of the parabolic Anderson problem (1.1) for potentials with nearly *double-exponential* and heavier tails [see (1.4) below]. Following [10], this comprises two of the four universality classes of asymptotic behavior for i.i.d. fields. For this class of potentials, the quenched and annealed long-time behavior of the total mass $U(t)$ has been studied by Gärtner and Molchanov [7, 8]. Furthermore, the asymptotic correlation structure of $u(t,\cdot)$ has been described by Gärtner and den Hollander [3]. We shall be able to use a number of results from those papers.

1.2. *Assumptions.* We are interested in the situation when the size of the islands of the relevant high peaks of the i.i.d. potential $\xi(\cdot)$ and of the solution $u(t,\cdot)$ stay bounded as $t\to\infty$. This will turn out to happen when the upper tail of the distribution of $\xi(0)$ lies in the vicinity of the *double-exponential distribution* with parameter $\varrho\in(0,\infty)$,

$$\operatorname{Prob}(\xi(0)>r)=\exp\{-e^{r/\varrho}\},\qquad r\in\mathbb{R}, \tag{1.4}$$

as well as for heavier tails ($\varrho=\infty$). In the latter case these islands are expected to shrink to single isolated lattice sites.

To be precise, let $F$ denote the distribution function of $\xi(0)$. Assume that $F$ is continuous and $F(r)<1$ for all $r\in\mathbb{R}$ (i.e., $\xi$ is unbounded from above). Introduce the nondecreasing function

$$\varphi(r)=\log\frac{1}{1-F(r)},\qquad r\in\mathbb{R}. \tag{1.5}$$

Then its left-continuous inverse $\psi$ is given by

$$\psi(s)=\min\{r\in\mathbb{R}:\varphi(r)\geq s\},\qquad s>0. \tag{1.6}$$

Note that $\psi$ is strictly increasing with $\varphi(\psi(s))=s$ for all $s>0$. The relevance of $\psi$ comes from the observation that $\xi$ has the same distribution as $\psi\circ\eta$, where $\eta=(\eta(x))_{x\in\mathbb{Z}^d}$ is an i.i.d. field of exponentially distributed random variables with mean one.

We now formulate our main assumption.

Assumption (F). *The distribution function $F$ is continuous, $F(r)<1$ for all $r\in\mathbb{R}$, and, in dimension $d=1$, $\int_{-\infty}^{-1}\log|r|\,F(dr)<\infty$. There exists $\varrho\in(0,\infty]$ such that*

$$\lim_{s\to\infty}[\psi(cs)-\psi(s)]=\varrho\log c,\qquad c\in(0,1). \tag{1.7}$$

*If $\varrho=\infty$, then $\psi$ satisfies, in addition,*

$$\lim_{s\to\infty}[\psi(s+\log s)-\psi(s)]=0. \tag{1.8}$$

This is Assumption (F) of [8]. The crucial supposition (1.7) specifies that the upper tail of the distribution of $\xi(0)$ is close to the double-exponential distribution (1.4) for $\varrho\in(0,\infty)$ and is heavier for $\varrho=\infty$. Assumption (1.8) excludes too heavy tails. Note that (1.8) is fulfilled for Gaussian but not for exponential tails. The extra assumption for $d=1$ rules out screening effects coming from extremely negative parts of the potential; these effects are not present in $d\ge 2$ by percolation.

The reader easily checks that (1.7) implies that

$$\lim_{t\to\infty}\frac{\psi(t)}{\log t}=\varrho. \tag{1.9}$$

1.3. *Spectral properties.* Before formulating our result, we introduce further notation and recall some results of [8].

Let $B_t=[-t,t]^d\cap\mathbb{Z}^d$ be the centered cube in $\mathbb{Z}^d$ with side length $2t$, and let

$$h_t=\max_{x\in B_t}\xi(x),\qquad t>0, \tag{1.10}$$

be the height of the potential $\xi$ in $B_t$. It can be easily seen [8], Corollary 2.7, that, under Assumption (F), almost surely,

$$h_t=\psi(d\log t)+o(1)\qquad \text{as } t\to\infty. \tag{1.11}$$

Let us remark that it is condition (1.8) which ensures that the almost sure asymptotics of $h_t$ in (1.11) is nonrandom up to order $o(1)$.

One of the main results in [8], Theorem 2.2, is the second-order asymptotics of the total mass $U(t)$ defined in (1.2). Under Assumption (F), with probability one,

$$\frac{1}{t}\log U(t)=h_t-\chi(\varrho)+o(1)\qquad \text{as } t\to\infty. \tag{1.12}$$

Here $\chi:(0,\infty]\to(0,2d]$ is a strictly increasing and surjective function. (Note that, by duality, $u(t,0)$ in [8] coincides with $U(t)$ in the present paper, and $\chi(\varrho)$ is identical to $2d\chi(\varrho)$ in the notation of [8].)

An analytic description of $\chi(\varrho)$ is as follows. Define $\mathcal{L}\colon[-\infty,0]^{\mathbb{Z}^d}\to[0,\infty]$ by

$$\mathcal{L}(V)=\begin{cases}\sum_{x\in\mathbb{Z}^d}e^{V(x)/\varrho}, & \text{if } \varrho\in(0,\infty),\\ |\{x\in\mathbb{Z}^d\colon V(x)>-\infty\}|, & \text{if } \varrho=\infty.\end{cases}\tag{1.13}$$

(We drop the dependence on $\varrho$ in this notation.) One should regard $\mathcal{L}$ as *large deviation rate function* for the fields $\xi-h_t$. Indeed, if the distribution of $\xi$ is exactly given by (1.4), then we have

$$\operatorname{Prob}(\xi(\cdot)-h>V(\cdot)\text{ in }\mathbb{Z}^d)=\exp\{-e^{h/\varrho}\mathcal{L}(V)\}$$

for any $V\colon\mathbb{Z}^d\to[-\infty,0]$ and any $h\in(0,\infty)$. For $V\in[-\infty,0]^{\mathbb{Z}^d}$, let $\lambda(V)\in[-\infty,0]$ be the top of the spectrum of the self-adjoint operator $\Delta+V$ in $\ell^2=\ell^2(\mathbb{Z}^d)$ in the domain $\{V>-\infty\}$ with zero boundary condition. In terms of the well-known Rayleigh–Ritz formula,

$$\lambda(V)=\sup_{f\in\ell^2(\mathbb{Z}^d)\colon\|f\|_2=1}\langle(\Delta+V)f,f\rangle,\tag{1.14}$$

where $\langle\cdot,\cdot\rangle$ and $\|\cdot\|_2$ denote the inner product and the norm in $\ell^2(\mathbb{Z}^d)$, respectively. Then

$$-\chi(\varrho)=\sup\{\lambda(V)\colon V\in[-\infty,0]^{\mathbb{Z}^d},\ \mathcal{L}(V)\le 1\};\tag{1.15}$$

see [8], Lemmas 2.17 and 1.10. This variational problem plays an important role in the study of our model.

For our deeper investigations, we introduce, in addition, an assumption about the optimal potential shape.

ASSUMPTION (M). *Up to spatial shifts, the variational problem in (1.15) possesses a unique maximizer, which has a unique maximum.*

By $V_\varrho$ we denote the unique maximizer of (1.15) which attains its unique maximum at the origin. We will call $V_\varrho$ *optimal potential shape*. Assumption (M) is satisfied at least for large $\varrho$. This fact, as well as further important properties of the variational problem (1.15), are stated in the next proposition.

PROPOSITION 1.1. (a) *For any $\varrho\in(0,\infty]$, the supremum in (1.15) is attained.*

(b) *If $\varrho$ is sufficiently large, then the maximizer in (1.15) is unique modulo shifts and has a unique maximum.*

(c) *If Assumption (M) is satisfied, then the optimal potential shape has the following properties:*

(i) *If* $\varrho \in (0,\infty)$, *then* $V_\varrho = f_\varrho \otimes \cdots \otimes f_\varrho$ *for some* $f_\varrho : \mathbb{Z} \to (-\infty, 0)$. *If* $\varrho = \infty$, *then* $V_\varrho$ *is degenerate in the sense that* $V_\varrho(0) = 0$ *and* $V_\varrho(x) = -\infty$ *for* $x \neq 0$.

(ii) *The operator* $\Delta + V_\varrho$ *has a unique nonnegative eigenfunction* $w_\varrho \in \ell^2(\mathbb{Z}^d)$ *with* $w_\varrho(0) = 1$ *corresponding to the eigenvalue* $\lambda(V_\varrho)$. *Moreover,* $w_\varrho \in \ell^1(\mathbb{Z}^d)$. *If* $\varrho \in (0,\infty)$, *then* $w_\varrho$ *is positive on* $\mathbb{Z}^d$, *while* $w_\varrho = \delta_0$ *for* $\varrho = \infty$.

The proof of Proposition 1.1 will be given in Section 3 below. There it will be shown that the variational problem in (1.15) is "dual" to a variational problem in terms of eigenfunctions, rather than potentials. The latter problem was studied in [3], Theorem 2, and reduces to finding a solution $v_\varrho : \mathbb{Z} \to (0,\infty)$ of the one-dimensional equation

$$\Delta v_\varrho + 2\varrho v_\varrho \log v_\varrho = 0$$

with minimal $\ell^2$-norm. For $\varrho$ large, $v_\varrho$ is unique up to shifts and has a unique maximum which we choose to be at the origin. Then the relation to our problem is given via $w_\varrho = \mathrm{const}\, v_\varrho \otimes \cdots \otimes v_\varrho$.

The question whether or not uniqueness in Assumption (M) is satisfied for *any* $\varrho$ is still open. Without uniqueness, the formulation and proofs of our results would be more cumbersome.

1.4. *Our result.* We shall show that the main contribution to the total mass $U(t)$ comes from a neighborhood of the set of best local coincidences of $\xi - h_t$ with shifts of $V_\varrho$ in the centered box of side length $2t\log^2 t$. These neighborhoods are widely separated from each other and hence not numerous. Furthermore, we show that we may further restrict ourselves to those neighborhoods in which, in addition, $u(t,\cdot)$, properly normalized, is close to $w_\varrho$.

Let us turn to a precise formulation. In order to include the degenerate case $\varrho = \infty$ in the statement, let

$$d_R(f,g) = \max_{x \in B_R} |e^{f(x)} - e^{g(x)}|, \qquad f, g \in [-\infty, \infty)^{B_R}.$$

Denote by $B_R(y) = y + B_R$ the closed cube of side length $2R$ centered at $y \in \mathbb{Z}^d$ and write

$$B_R(A) = \bigcup_{y \in A} B_R(y) \tag{1.16}$$

for the "$R$-cube neighborhood" of a set $A \subset \mathbb{Z}^d$. In particular, $B_0(A) = A$. Let $|A|$ denote the cardinality of $A$. We furnish $\mathbb{Z}^d$ with the lattice norm $|y| = \sum_{i=1}^d |y_i|$, where $y = (y_1, \ldots, y_d)$.

Let $\varrho \in (0,\infty]$ be so large that Assumption (M) is satisfied. For any $\varepsilon > 0$, let $r(\varrho,\varepsilon)$ denote the smallest $r \in \mathbb{N}_0$ such that

$$\|w_\varrho\|_2^2 \sum_{x\in\mathbb{Z}^d\setminus B_r} w_\varrho(x) < \varepsilon. \tag{1.17}$$

Note that $r(\infty,\varepsilon) = 0$, due to the degeneracy of $w_\infty$.

Our main result is the following.

Theorem 1.2. *Let the Assumptions (F) and (M) be satisfied. Then there exists a random $t$-dependent subset $\Gamma^* = \Gamma^*_{t\log^2 t}$ of $B_{t\log^2 t}$ such that almost surely:*

$$\text{(i)}\quad \liminf_{t\to\infty} \frac{1}{U(t)} \sum_{x\in B_{r(\varrho,\varepsilon)}(\Gamma^*)} u(t,x) \geq 1-\varepsilon, \qquad \varepsilon \in (0,1); \tag{1.18}$$

$$\text{(ii)}\quad |\Gamma^*| \leq t^{o(1)} \quad \textit{and} \quad \min_{y,\widetilde{y}\in\Gamma^* : y\neq\widetilde{y}} |y-\widetilde{y}| \geq t^{1-o(1)} \qquad \textit{as } t\to\infty; \tag{1.19}$$

$$\text{(iii)}\quad \lim_{t\to\infty} \max_{y\in\Gamma^*} d_R(\xi(y+\cdot) - h_t, V_\varrho(\cdot)) = 0, \qquad R > 0; \tag{1.20}$$

$$\text{(iv)}\quad \lim_{t\to\infty} \max_{y\in\Gamma^*} d_R\left(\frac{u(t,y+\cdot)}{u(t,y)}, w_\varrho(\cdot)\right) = 0, \qquad R > 0. \tag{1.21}$$

Theorem 1.2 states that, up to an arbitrarily small relative error $\varepsilon$, the islands with centers in $\Gamma^*$ and radius $r(\varrho,\varepsilon)$ carry the whole mass of the solution $u(t,\cdot)$. Locally, in an arbitrarily fixed $R$-neighborhood of each of these centers, the shapes of the potential and the normalized solution resemble $h_t + V_\varrho$ and $w_\varrho$, respectively. The number of these islands increases at most as an arbitrarily small power of $t$ and their distance increases almost like $t$. Note that, for $\varrho = \infty$, the set $B_{r(\varrho,\varepsilon)}(\Gamma^*)$ in (1.18) is equal to $\Gamma^*$, that is, the islands consist of single lattice sites.

In this paper we have made no attempt to choose $\Gamma^*$ as small as possible. We mention without proof that, for Weibull tails, when $\mathrm{Prob}(\xi(0) > r) = \exp\{-r^\alpha L(r)\}$ with $\alpha > 1$ and $L(r)$ slowly varying as $r\to\infty$, one may choose $|\Gamma^*| = o(\log^{1+\varepsilon} t)$ with $\varepsilon > 0$ arbitrarily small. In this case it might even turn out that $|\Gamma^*|$ stays bounded in probability as $t\to\infty$. But this will be the subject of another analysis.

**2. Strategy of the proof.** Our proof relies on the Feynman–Kac representation for the solution $u$ of (1.1),

$$u(t,x) = \mathbb{E}_x \exp\left\{\int_0^t \xi(X_s)\,ds\right\}\delta_0(X_t), \tag{2.1}$$

where $(X_t)_{t\in[0,\infty)}$ denotes the time-continuous nearest-neighbor random walk on $\mathbb{Z}^d$ with generator $\Delta$, starting at $x\in\mathbb{Z}^d$ under $\mathbb{P}_x$. We introduce the entrance time into a set $A\subset\mathbb{Z}^d$,

$$\tau_A=\inf\{t\geq 0\colon X_t\in A\} \tag{2.2}$$

and write $\tau_z$ instead of $\tau_{\{z\}}$ for $z\in\mathbb{Z}^d$.

In our proof we shall need a random auxiliary subset $\Gamma=\Gamma_{t\log^2 t}$ of the large box $B_{t\log^2 t}$, which will be introduced in Proposition 2.2 below. The set $\Gamma^*=\Gamma^*_{t\log^2 t}$ of Theorem 1.2 will later be constructed as a certain subset of $\Gamma_{t\log^2 t}$. Given $t>0$, we split $u$ into three terms, $u=u_1+u_2+u_3$, where, for $s\geq 0$,

$$u_1(s,x)=\mathbb{E}_x\exp\biggl\{\int_0^s\xi(X_u)\,du\biggr\}\delta_0(X_s)\mathbb{1}\{\tau_{B^{\mathrm{c}}_{t\log^2 t}}\leq s\}, \tag{2.3}$$

$$u_2(s,x)=\mathbb{E}_x\exp\biggl\{\int_0^s\xi(X_u)\,du\biggr\}\delta_0(X_s)\mathbb{1}\{\tau_{B^{\mathrm{c}}_{t\log^2 t}}> s\}\mathbb{1}\{\tau_{\Gamma_{t\log^2 t}}>s\}, \tag{2.4}$$

$$u_3(s,x)=\mathbb{E}_x\exp\biggl\{\int_0^s\xi(X_u)\,du\biggr\}\delta_0(X_s)\mathbb{1}\{\tau_{B^{\mathrm{c}}_{t\log^2 t}}> s\}\mathbb{1}\{\tau_{\Gamma_{t\log^2 t}}\leq s\}. \tag{2.5}$$

In words, in $u_1$ we have the contribution from the paths that reach the complement of the "macrobox" $B_{t\log^2 t}$ by time $s$. In $u_2$ we consider the paths that stay inside this box, but do not enter the set $\Gamma_{t\log^2 t}$ up to time $s$. In $u_3$ they stay inside $B_{t\log^2 t}$ and do enter $\Gamma_{t\log^2 t}$. Note that $u_2(t,\cdot)=0$ on $B^{\mathrm{c}}_{t\log^2 t}\cup\Gamma_{t\log^2 t}$, and $u_3(t,\cdot)=0$ on $B^{\mathrm{c}}_{t\log^2 t}$. The functions $u_1$, $u_2$ and $u_3$ depend on $t$ via the sets $B_{t\log^2 t}$ and $\Gamma_{t\log^2 t}$. We shall mainly be interested in the case $s=t$ as $t\to\infty$. Note that $u_1$, $u_2$ and $u_3$ are solutions to certain initial-boundary value problems for the parabolic differential equation in (1.1).

The contribution to the total mass $U(t)$ coming from $u_1$ turns out to be negligible:

PROPOSITION 2.1. *Let Assumption (F) be satisfied. Then, with probability one,*

$$\lim_{t\to\infty}\frac{1}{U(t)}\sum_{x\in\mathbb{Z}^d}u_1(t,x)=0. \tag{2.6}$$

PROOF. This follows from (2.41) and the text below (2.42) in [8]. Indeed, note that the function $u$ in [8] is the solution to (1.1) with initial datum $u(0,\cdot)=1$ instead of $u(0,\cdot)=\delta_0(\cdot)$. Hence, the second line of [8], (2.41), is equal to $\sum_{x\in\mathbb{Z}^d}u_1(t,x)$. Below (2.42) it is shown that this term vanishes as $t\to\infty$,

$$\lim_{t\to\infty}\sum_{x\in\mathbb{Z}^d}u_1(t,x)=0\qquad\text{almost surely.} \tag{2.7}$$

Since $U(t)\to\infty$ [see (1.12)], the assertion follows. □

One of the essential points in our proof of Theorem 1.2 is to construct the (random and $t$-dependent) set $\Gamma$ appearing in (2.4) and (2.5) in such a way that assertion (1.18) is satisfied for $u_3$ and $\Gamma$ in place of $u$ and $\Gamma^*$, (1.19) and (1.20) are satisfied for $\Gamma$ in place of $\Gamma^*$, and $u_2$ is asymptotically negligible.

PROPOSITION 2.2. *Let the Assumptions (F) and (M) be satisfied. Fix $\varepsilon,\eta,\gamma,R>0$ and $\alpha>0$ arbitrarily. Then there exists a random $t$-dependent subset $\Gamma=\Gamma_{t\log^2 t}$ of $B_{t\log^2 t}$ such that, almost surely, for $t$ sufficiently large:*

$$\text{(i}'\text{)}\quad \frac{1}{U(t)}\sum_{x\in B_{t\log^2 t}\setminus B_{r(\varrho,\varepsilon')}(\Gamma)} u_3(t,x)\le\varepsilon'+\alpha,\qquad \varepsilon'\in(\varepsilon,1);\tag{2.8}$$

$$\text{(ii}'\text{)}\quad |\Gamma|\le t^{\eta d}\quad\text{and}\quad \min_{y,\widetilde{y}\in\Gamma:\,y\neq\widetilde{y}}|y-\widetilde{y}|\ge t^{1-\eta};\tag{2.9}$$

$$\text{(iii}'\text{)}\quad d_R(\xi(y+\cdot)-h_t,V_\varrho(\cdot))<\gamma,\qquad y\in\Gamma;\tag{2.10}$$

*furthermore,*

$$\lim_{t\to\infty}\frac{1}{U(t)}\sum_{x\in B_{t\log^2 t}} u_2(t,x)=0.\tag{2.11}$$

The proof of this proposition will be carried out in Sections 4–7.

In the proof of our next and final proposition we shall construct the set $\Gamma^*$ of Theorem 1.2 as a certain subset of $\Gamma$ from Proposition 2.2 in such a way that the conditions (i) and (iv) of Theorem 1.2 are satisfied. [The conditions (ii) and (iii) are trivially satisfied for any subset $\Gamma^*$ of $\Gamma$.]

PROPOSITION 2.3. *Let the Assumptions (F) and (M) be satisfied. Fix $\varepsilon,\eta,\gamma,R>0$ and $\alpha>0$ arbitrarily, and let $\Gamma$ be the random $t$-dependent set constructed in the proof of Proposition 2.2. Then there exists a random $t$-dependent subset $\Gamma^*$ of $\Gamma$ such that, almost surely, for $t$ sufficiently large:*

$$\text{(i}''\text{)}\quad \sum_{x\in B_{r(\varrho,\varepsilon')}(\Gamma\setminus\Gamma^*)} u_3(t,x)\le t^{-\eta d}\sum_{x\in B_{r(\varrho,\varepsilon')}(\Gamma^*)} u(t,x),\qquad \varepsilon'\in(\varepsilon,1);\tag{2.12}$$

$$\text{(iv}'\text{)}\quad d_R\biggl(\frac{u(t,y+\cdot)}{u(t,y)},w_\varrho(\cdot)\biggr)<\gamma,\qquad y\in\Gamma^*.\tag{2.13}$$

The proof of this proposition will be given in Section 8.

Let us now finish the proof of our main result subject to Propositions 2.1–2.3.

PROOF OF THEOREM 1.2. Let $\varepsilon_n$, $\eta_n$, $\gamma_n$, $R_n$ and $\alpha_n$ be positive numbers such that $\varepsilon_n,\eta_n,\gamma_n,\alpha_n\to 0$ and $R_n\to\infty$ monotonically as $n\to\infty$. Let $\Gamma_n=\Gamma_{n,t\log^2 t}$ and $\Gamma^*_n=\Gamma^*_{n,t\log^2 t}$ be the random sets of Propositions 2.2 and 2.3 for $\varepsilon,\eta,\gamma,R$ and $\alpha$ replaced by $\varepsilon_n,\eta_n,\gamma_n,R_n$ and $\alpha_n$, respectively. Let $u_1^{(n)},u_2^{(n)},u_3^{(n)}$ be the functions in (2.3)–(2.5) with $\Gamma_{t\log^2 t}$ replaced by $\Gamma_{n,t\log^2 t}$. In accordance with Propositions 2.1 and 2.2, there exist random times $t_n\uparrow\infty$ such that, almost surely, for $t>t_n$,

$$\frac{1}{U(t)}\sum_{x\in\mathbb{Z}^d}u_1^{(n)}(t,x)\le\frac{1}{n};\tag{2.14}$$

$$\frac{1}{U(t)}\sum_{x\in B_{t\log^2 t}\setminus B_{r(\varrho,\varepsilon')}(\Gamma_n)}u_3^{(n)}(t,x)\le\varepsilon'+\alpha_n,\qquad \varepsilon'\in(\varepsilon_n,1);\tag{2.15}$$

$$|\Gamma_n|\le t^{\eta_n d}\quad\text{and}\quad\min_{y,\widetilde{y}\in\Gamma_n\,:\,y\ne\widetilde{y}}|y-\widetilde{y}|\ge t^{1-\eta_n};\tag{2.16}$$

$$d_{R_n}(\xi(y+\cdot)-h_t,V_\varrho(\cdot))<\gamma_n,\qquad y\in\Gamma_n;\tag{2.17}$$

$$\frac{1}{U(t)}\sum_{x\in B_{t\log^2 t}}u_2^{(n)}(t,x)\le\frac{1}{n}.\tag{2.18}$$

In accordance with Proposition 2.3, we also may assume that, for $t>t_n$,

$$\frac{1}{U(t)}\sum_{x\in B_{r(\varrho,\varepsilon')}(\Gamma_n\setminus\Gamma^*_n)}u_3^{(n)}(t,x)\le\frac{1}{n},\qquad \varepsilon'\in(\varepsilon_n,1);\tag{2.19}$$

$$d_{R_n}\left(\frac{u(t,y+\cdot)}{u(t,y)},w_\varrho(\cdot)\right)<\gamma_n,\qquad y\in\Gamma^*_n.\tag{2.20}$$

Now define $\Gamma^*_{t\log^2 t}=\Gamma^*_{n,t\log^2 t}$ if $t_n\le t<t_{n+1}$. Combining the bounds (2.14)–(2.20), one easily checks that $\Gamma^*_{t\log^2 t}$ satisfies the assertions (i)–(iv) of Theorem 1.2. □

The remainder of this paper is organized as follows. Section 3 provides important background material from functional analysis for the coming developments. Our strategy for estimating the left-hand side of (2.8) is developed in Section 4. In Section 5 we construct the set Γ of Proposition 2.2. In Section 6 we derive important properties of this set and of related objects. The proof of Proposition 2.2 will be completed in Section 7. Finally, in Section 8, we finish the proof of our main Theorem 1.2 by proving Proposition 2.3.

**3. Functional analytic background.** In this section we collect and prove a number of results evolving around the variational formula in (1.15). In particular, we prove Proposition 1.1 and provide important tools for the proof of Theorem 1.2. This section is purely functional analytic.

Recall that $\lambda(V)$ defined in (1.14) denotes the top of the spectrum of the operator $\Delta+V$ in $\ell^2(\mathbb{Z}^d)$. Note that, if $V(x)\to-\infty$ as $|x|\to\infty$ [which is the case, e.g., if $\mathcal{L}(V)\le 1$], then the operator $\Delta+V$ has compact resolvent and $\lambda(V)$ is an isolated simple eigenvalue corresponding to a positive eigenfunction.

LEMMA 3.1. *Assume that $\varrho\in(0,\infty)$. Then*

$$\lambda(V)\le-\chi(\varrho)+\varrho\log\mathcal{L}(V) \tag{3.1}$$

*for all $V\in[-\infty,\infty)^{\mathbb{Z}^d}$.*

PROOF. Assume, without loss of generality, that $\mathcal{L}(V)<\infty$. Then $V$ is bounded from above, and $\mathcal{L}(V+c)=\mathcal{L}(V)e^{c/\varrho}=1$ for $c=-\varrho\log\mathcal{L}(V)$. It follows from (1.15) that

$$\lambda(V+c)\le-\chi(\varrho). \tag{3.2}$$

Since $\lambda(V+c)=\lambda(V)+c$, this implies (3.1). □

For any nonempty set $A\subset\mathbb{Z}^d$ and any $V\in[-\infty,0]^A$, define

$$\mathcal{L}_A(V)=\begin{cases}\displaystyle\sum_{x\in A}e^{V(x)/\varrho}, & \text{if } \varrho\in(0,\infty),\\ |\{x\in A:V(x)>-\infty\}|, & \text{if } \varrho=\infty.\end{cases} \tag{3.3}$$

Let $\lambda_A(V)$ be the principal (i.e., largest) eigenvalue of the operator $\Delta+V$ in $\ell^2(A\cap\{V>-\infty\})$ with zero boundary condition. By the Rayleigh–Ritz formula,

$$\lambda_A(V)=\sup_{f\in\ell^2(\mathbb{Z}^d):\,\mathrm{supp}(f)\subset A,\|f\|_2=1}\langle(\Delta+V)f,f\rangle. \tag{3.4}$$

We define the finite-volume version of $\chi(\varrho)$ by

$$-\chi_R(\varrho)=\sup\{\lambda_{B_R}(V):V\in[-\infty,0]^{B_R},\mathcal{L}_{B_R}(V)\le 1\}. \tag{3.5}$$

Let $\mathcal{M}_\varrho$ denote the set of maximizers in (1.15). The topology of pointwise convergence on $[-\infty,0]^{\mathbb{Z}^d}$ is induced by some metric $d(\cdot,\cdot)$ and making $[-\infty,0]^{\mathbb{Z}^d}$ compact. Let $\arg\max(V)$ denote the set of sites at which $V$ attains its maximum.

LEMMA 3.2. (i) *The supremum in (1.15) is attained, that is, $\mathcal{M}_\varrho\neq\varnothing$. For any $\varepsilon>0$,*

$$\begin{aligned}&\sup\{\lambda(V):V\in[-\infty,0]^{\mathbb{Z}^d},\mathcal{L}(V)\le 1,\\ &\qquad 0\in\arg\max(V),d(V,\mathcal{M}_\varrho)\ge\varepsilon\}<-\chi(\varrho).\end{aligned} \tag{3.6}$$

(ii) *If $\varrho$ is large enough, then the maximizer $V_\varrho$ in (1.15) is unique up to shifts and possesses a single-point maximum, which we assume to be at the origin. For such $\varrho$, $V_\varrho = f_\varrho \oplus \cdots \oplus f_\varrho$ with $f_\varrho : \mathbb{Z} \to (-\infty, 0)$ if $\varrho < \infty$, and $f_\infty(0) = 0$ and $f_\infty(x) = -\infty$ for $x \neq 0$.*

PROOF. Since both assertions are trivial for $\varrho = \infty$, we assume that $\varrho \in (0, \infty)$. By $\partial A$ we denote the inner boundary of a set $A \subset \mathbb{Z}^d$.

*Proof of* (i). Let us begin with two preliminary steps.

STEP 1. *For any $V \in [-\infty, \infty)^{\mathbb{Z}^d}$ and any $R \in \mathbb{N}$,*

$$\lambda(V) \le \max\{\lambda_{B_R}(V + 2d\mathbb{1}_{\partial B_R}), \lambda_{B_R^c}(V + 2d\mathbb{1}_{\partial B_R^c})\}. \tag{3.7}$$

PROOF. Given $f \in \ell^2(\mathbb{Z}^d)$ with $\|f\|_2 = 1$, define $f_R = f\mathbb{1}_{B_R}$ and $f_R^c = f\mathbb{1}_{B_R^c}$. Then

$$\langle (\Delta + V)f, f\rangle = \langle (\Delta + V)f_R, f_R\rangle + \langle (\Delta + V)f_R^c, f_R^c\rangle + 2\langle \Delta f_R, f_R^c\rangle. \tag{3.8}$$

One easily checks that

$$\langle \Delta f_R, f_R^c\rangle \le d \sum_{x \in \partial B_R} f_R(x)^2 + d \sum_{x \in \partial B_R^c} f_R^c(x)^2. \tag{3.9}$$

Hence,

$$\begin{aligned} \langle (\Delta + V)f, f\rangle &\le \langle (\Delta + V + 2d\mathbb{1}_{\partial B_R})f_R, f_R\rangle \\ &\qquad + \langle (\Delta + V + 2d\mathbb{1}_{\partial B_R^c})f_R^c, f_R^c\rangle. \end{aligned} \tag{3.10}$$

Now pass to the supremum over all such $f$ and note that $\|f_R\|_2^2 + \|f_R^c\|_2^2 = \|f\|_2^2 = 1$ to arrive at

$$\lambda(V) \le \sup_{\alpha \in [0,1]} [\alpha \lambda_{B_R}(V + 2d\mathbb{1}_{\partial B_R}) + (1 - \alpha)\lambda_{B_R^c}(V + 2d\mathbb{1}_{\partial B_R^c})]. \tag{3.11}$$

This certainly implies assertion (3.7). □

STEP 2. *For any $W \in [-\infty, 0]^{\mathbb{Z}^d}$ and any $R \in \mathbb{N}$,*

$$\lambda(W) - \lambda_{B_R}(W) \le \frac{(2d - \lambda_{B_R}(W))^2}{\min_{B_{R-1}^c} |W|}. \tag{3.12}$$

PROOF. Fix $\varepsilon > 0$ arbitrarily and pick $f \in \ell^2(\mathbb{Z}^d)$ with $\|f\|_2 = 1$ such that

$$\lambda(W) - \varepsilon \le \langle (\Delta + W)f, f\rangle. \tag{3.13}$$

Combining this with (3.8) and (3.9) with $V$ and $R$ replaced by $W$ and $R-1$, respectively, and taking into account that $\Delta$ is negative definite and $W \le 0$, we see that

$$\lambda(W) - \varepsilon \le -\min_{B_{R-1}^{\rm c}} |W| \|f_{R-1}^{\rm c}\|_2^2 + 2d.$$

Hence, since $\lambda(W) \ge \lambda_{B_R}(W)$,

$$\|f_{R-1}^{\rm c}\|_2^2 \le \frac{2d + \varepsilon - \lambda_{B_R}(W)}{\min_{B_{R-1}^{\rm c}} |W|}. \tag{3.14}$$

Combining again (3.13) with (3.8) and (3.9) with $V$ replaced by $W$ and using the Rayleigh–Ritz formula for the first term on the right-hand side of (3.8), we find that

$$\begin{aligned} \lambda(W) - \varepsilon &\le \lambda_{B_R}(W)\|f_R\|_2^2 + 2d\|f_{R-1}^{\rm c}\|_2^2 \\ &\le \lambda_{B_R}(W) + (2d - \lambda_{B_R}(W))\|f_{R-1}^{\rm c}\|_2^2. \end{aligned} \tag{3.15}$$

Substituting (3.14) into (3.15) and letting $\varepsilon \downarrow 0$, we arrive at the desired assertion. □

Now let $(V_n)_{n\in\mathbb{N}}$ be a maximizing sequence for the variational problem in (1.15), that is, $V_n \in [-\infty, 0]^{\mathbb{Z}^d}$, $\mathcal{L}(V_n) \le 1$ for all $n$, and $\lim_{n\to\infty} \lambda(V_n) = -\chi(\varrho)$. We may assume that every $V_n$ attains its maximum at the origin. By compactness of $[-\infty, 0]^{\mathbb{Z}^d}$, we may also assume that $V_n$ converges toward some $V \in [-\infty, 0]^{\mathbb{Z}^d}$ pointwise. Clearly, $\mathcal{L}(V) \le 1$. By Lemma 3.1, $\lim_{n\to\infty} \mathcal{L}(V_n) = 1$. Since $\lambda(V_n) \le V_n(0)$, $V$ is not identically equal to $-\infty$.

For proving both assertions in (i), it only remains to show that $\lambda(V) \ge -\chi(\varrho)$. We will do this by checking that $\limsup_{n\to\infty} \lambda(V_n) \le \lambda(V)$. A preliminary step is the following.

STEP 3. $\lim_{R\to\infty} \limsup_{n\to\infty} \sup_{B_R^{\rm c}} V_n = -\infty$.

PROOF. It is clearly sufficient to prove that

$$\lim_{R\to\infty} \limsup_{n\to\infty} \mathcal{L}_{B_R^{\rm c}}(V_n) = 0. \tag{3.16}$$

Fix $R \in \mathbb{N}$. Combine Step 1 (for $V_n$ instead of $V$) and Lemma 3.1 to obtain

$$\lambda(V_n) \le -\chi(\varrho) + \varrho \log[\max\{\mathcal{L}_{B_R}(V_n + 2d\mathbb{1}_{\partial B_R}), \mathcal{L}_{B_R^{\rm c}}(V_n + 2d\mathbb{1}_{\partial B_R^{\rm c}})\}].$$

Since $\lambda(V_n) \to -\chi(\varrho)$, we deduce that

$$\liminf_{n\to\infty} \max\{\mathcal{L}_{B_R}(V_n + 2d\mathbb{1}_{\partial B_R}), \mathcal{L}_{B_R^{\rm c}}(V_n + 2d\mathbb{1}_{\partial B_R^{\rm c}})\} \ge 1.$$

We next show that, if $R$ is large enough, then, for large $n$, the maximum is equal to the first of the two terms. Indeed,

$$\begin{aligned}\mathcal{L}_{B_R^c}(V_n + 2d\mathbb{1}_{\partial B_R^c}) &= \mathcal{L}(V_n) - \mathcal{L}_{B_R}(V_n) + (e^{2d/\varrho} - 1) \sum_{x \in \partial B_R^c} e^{V_n(x)/\varrho} \\ &\to 1 - \mathcal{L}_{B_R}(V) + (e^{2d/\varrho} - 1) \sum_{x \in \partial B_R^c} e^{V(x)/\varrho},\end{aligned}$$

which is bounded away from 1 as $R$ gets large.

Consequently, if $R$ is large enough,

$$\begin{aligned}1 &\le \liminf_{n\to\infty} \mathcal{L}_{B_R}(V_n + 2d\mathbb{1}_{\partial B_R}) \\ &= \liminf_{n\to\infty} [\mathcal{L}(V_n) - \mathcal{L}_{B_R^c}(V_n)] + (e^{2d/\varrho} - 1) \sum_{x \in \partial B_R} e^{V(x)/\varrho}.\end{aligned}$$

Since $\mathcal{L}(V_n) \to 1$, we conclude that, for $R$ large,

$$\limsup_{n\to\infty} \mathcal{L}_{B_R^c}(V_n) \le (e^{2d/\varrho} - 1) \sum_{x \in \partial B_R^c} e^{V(x)/\varrho}.$$

Now (3.16) follows by letting $R \to \infty$. $\square$

Now use Step 2 for $W = V_n$ and Step 3 to see that

$$\lim_{R\to\infty} \limsup_{n\to\infty} [\lambda(V_n) - \lambda_{B_R}(V_n)] = 0.$$

Since $\lambda_{B_R}(V_n) \to \lambda_{B_R}(V)$ as $n \to \infty$ and $\lambda_{B_R}(V) \to \lambda(V)$ as $R \to \infty$, it follows that $\lambda(V) = \limsup_{n\to\infty} \lambda(V_n) = -\chi(\varrho)$. This finishes the proof of (i).

*Proof of* (ii). Let us rewrite our variational problem in (1.15) in the form

$$\chi(\varrho) = \inf\{-\lambda(V) : \mathcal{L}(V) = 1\}. \tag{3.17}$$

Since the functional $\lambda - \varrho \log \mathcal{L}$ remains invariant under shifts $V \mapsto V + \text{const}$,

$$\chi(\varrho) = \inf\{-\lambda(V) + \varrho \log \mathcal{L}(V) : \mathcal{L}(V) < \infty\}. \tag{3.18}$$

Let $\mathcal{P}$ denote the set of probability measures on $\mathbb{Z}^d$, and introduce the functionals $I, J : \mathcal{P} \to [0, \infty]$ by

$$I(p) = -\langle \Delta\sqrt{p}, \sqrt{p}\rangle \quad \text{and} \quad J(p) = -\langle p, \log p\rangle. \tag{3.19}$$

In [3] the variational problem

$$\widetilde{\chi}(\varrho) = \inf_{\mathcal{P}} [I + \varrho J] \tag{3.20}$$

has been considered (we write $\widetilde{\chi}$ for the $\chi$ of [3] and omit the factor $1/2d$ which appeared there) which will turn out to be "dual" to the problem in (3.18). It has been shown in [3], Theorem 2.II, that problem (3.20) possesses a minimizer, which is unique up to shifts and has a single-point maximum provided that $\varrho$ is large. Therefore, to prove the first part of assertion (ii), it will be enough to verify the following equivalence for $p^* \in \mathcal{P}$:

$$\begin{aligned} &p^* \text{ minimizes } I + \varrho J \\ &\qquad \Longleftrightarrow \quad V^* = \varrho \log p^* \text{ minimizes } -\lambda(\cdot) \text{ under } \mathcal{L}(\cdot) = 1. \end{aligned} \tag{3.21}$$

[Note that $\mathcal{L}(\varrho \log p^*) = 1$ if and only if $p^* \in \mathcal{P}$.] As a first step, let us show that

$$\widetilde{\chi}(\varrho) = \chi(\varrho). \tag{3.22}$$

To this end, we remark that

$$\varrho J(p) = \inf_{V : \mathcal{L}(V) < \infty} [-\langle V, p\rangle + \varrho \log \mathcal{L}(V)], \tag{3.23}$$

since the functional under the infimum is convex and attains its minimum at $V = \varrho \log p$. Hence, using (3.23), a slight modification of the Rayleigh–Ritz formula (1.14) and (3.18), we obtain

$$\begin{aligned} \widetilde{\chi}(\varrho) &= \inf_{\mathcal{P}} [I + \varrho J] \\ &= \inf_{p \in \mathcal{P}} \inf_{V : \mathcal{L}(V) < \infty} [\langle -\Delta \sqrt{p}, \sqrt{p}\rangle - \langle V, p\rangle + \varrho \log \mathcal{L}(V)] \\ &= \inf_{V : \mathcal{L}(V) < \infty} [-\lambda(V) + \varrho \log \mathcal{L}(V)] \\ &= \chi(\varrho). \end{aligned} \tag{3.24}$$

To prove "$\Longrightarrow$" in (3.21), assume that $p^* \in \mathcal{P}$ minimizes $I + \varrho J$, and put $V^* = \varrho \log p^*$. Then, again using the Rayleigh–Ritz formula, we find that

$$-\lambda(V^*) \le -\langle (\Delta + V^*) \sqrt{p^*}, \sqrt{p^*}\rangle = [I + \varrho J](p^*) = \widetilde{\chi}(\varrho) = \chi(\varrho).$$

Hence, $V^*$ minimizes $-\lambda(\cdot)$ under $\mathcal{L}(\cdot) = 1$.

Before proving the reversed implication in (3.21), let us remark that for $\mathcal{L}(V) < \infty$, by perturbation theory of linear operators, the Gateaux derivatives of $V \mapsto \lambda(V)$ and of the corresponding positive eigenfunction $V \mapsto v_V$ with $\|v_V\|_2 = 1$ in any direction $W \in \mathbb{R}^{\mathbb{Z}^d}$ with finite support, $\partial_W \lambda(V)$ and $\partial_W v_V$, exist; see [11], Theorems VIII.2.6 and VIII.2.9. Since $\langle v_V, v_V\rangle = 1$, we have $\langle v_V, \partial_W v_V\rangle = 0$. Differentiating the eigenvalue equation

$$(\Delta + V - \lambda(V)) v_V = 0$$

in direction $W$, we obtain

$$(\Delta + V - \lambda(V))\, \partial_W v_V = (\partial_W \lambda(V) - W) v_V.$$

Taking the inner product with $v_V$, we find that

$$\langle \partial_W \lambda(V) - W, v_V^2 \rangle = 0,$$

hence, $\nabla\lambda(V) = v_V^2$.

Now suppose that $V^*$ is a minimizer of $-\lambda(\cdot)$ under $\mathcal{L}(\cdot) = 1$, and put $p^* = e^{V^*/\varrho}$. Since $V^*$ also minimizes the functional $-\lambda + \varrho \log \mathcal{L}$, we get

$$0 = \nabla(-\lambda + \varrho \log \mathcal{L})(V^*) = -v_{V^*}^2 + e^{V^*/\varrho} = p^* - v_{V^*}^2,$$

that is, $p^* = v_{V^*}^2$. Hence,

$$\begin{aligned}[I + \varrho J](p^*) &= -\langle \Delta\sqrt{p^*}, \sqrt{p^*}\rangle - \varrho\langle \log p^*, p^*\rangle \\ &= \langle(-\Delta - V^*)v_{V^*}, v_{V^*}\rangle = -\lambda(V^*) = \chi(\varrho) = \widetilde{\chi}(\varrho).\end{aligned}$$

This shows that $p^*$ minimizes $I + \varrho J$.

If $\varrho$ is large, then we also know from [3] (see Proposition 3 and Theorem 2) that the unique centered minimizer $p^*$ of (3.20) has the form $p^* = p_0^* \otimes \cdots \otimes p_0^*$ for some probability measure $p_0^*$ on $\mathbb{Z}$, which has a unique maximum at the origin. Therefore, the second part of assertion (ii) follows from this by using again the one-to-one correspondence in (3.21).

PROOF OF PROPOSITION 1.1. For finite $\varrho$, all assertions of Proposition 1.1 are either standard or have been proved in Lemma 3.2. For $\varrho = \infty$, the problem in (1.15) is trivial since there is (up to shifts) only one admissible element $V \in [-\infty, 0]^{\mathbb{Z}^d}$, that is, only one $V$ satisfies $\mathcal{L}(V) \le 1$. □

Now we provide several finite-space approximation results for the maximizer $V_\varrho$ in (1.15) and corresponding eigenfunction $w_\varrho$ of $\Delta + V_\varrho$. For $p \in \{1, 2, \infty\}$, let us write $\|\cdot\|_{p,R}$ for the restriction of the $\ell^p$-norm $\|\cdot\|_p$ to $B_R$.

Let $w_\varrho^{(R)}: \mathbb{Z}^d \to [0, \infty)$ be the positive eigenfunction of the operator $\Delta + V_\varrho$ with zero boundary condition in $B_R$ corresponding to the eigenvalue $\lambda_{B_R}(V_\varrho)$. We norm $w_\varrho^{(R)}$ such that $w_\varrho^{(R)}(0) = 1$. In particular, $w_\varrho^{(\infty)} = w_\varrho$ with obvious notation. Note that $w_\varrho^{(R)}$ vanishes outside of $B_R$.

We also need eigenvalues in *dotted* sets, more precisely, in boxes whose center point has been removed. For any centered box $B$ and any $V \in [-\infty, \infty)^B$, we introduce the notation $\dot{V}$ for the function in $[-\infty, \infty)^B$ that is identical to $V$ in $B \setminus \{0\}$ and satisfies $\dot{V}(0) = -\infty$. Then $\lambda_B(\dot{V})$ is the principal Dirichlet eigenvalue of $\Delta + V$ in $B \setminus \{0\}$. Analogous notation is used for $B = \mathbb{Z}^d$. Clearly, $\lambda(\dot{V}_\varrho) < \lambda(V_\varrho)$.

LEMMA 3.3. *Assume that Assumptions* (F) *and* (M) *are satisfied. Then:*

(i) $\lim_{R\to\infty}\lambda_{B_R}(V_\varrho)=\lambda(V_\varrho)$ *and* $\lim_{R\to\infty}\lambda_{B_R}(\dot V_\varrho)=\lambda(\dot V_\varrho)$*;*
(ii) $\lim_{R\to\infty}\chi_R(\varrho)=\chi(\varrho)$*;*
(iii) $\lim_{R\to\infty}\|w_\varrho^{(R)}-w_\varrho\|_2=0=\lim_{R\to\infty}\|w_\varrho^{(R)}-w_\varrho\|_1$.

Proof. Since $\lim_{x\to\infty}V_\varrho(x)=-\infty$, assertion (i) follows by a standard compactness argument. Assertion (ii) easily follows from (i) since $-\chi(\varrho)=\lambda(V_\varrho)$ and since $\chi_R(\varrho)\geq\chi(\varrho)$ for any $R$. From [3], Theorem 2I.(3)(iii), it follows that $w_\varrho$ lies in $\ell^1(\mathbb{Z}^d)$ [and, hence, also in $\ell^2(\mathbb{Z}^d)$]. The proof of the two assertions in (iii) is hence standard. □

In our proof of Theorem 1.2 below we shall need the following assertion. It says that optimality of the eigenvalue $\lambda_{B_\mathcal{R}}(\cdot)$ in a sufficiently large box $B_\mathcal{R}$ implies closeness to $V_\varrho$ in a given box $B_R$. Recall that $d_R$ is the uniform metric on $[-\infty,\infty)^{B_R}$.

Corollary 3.4. *Fix* $\gamma>0$ *and* $R\in(0,\infty)$ *arbitrarily. Then there exists* $\delta_0>0$ *and* $\mathcal{R}_0>0$ *such that, for any* $\mathcal{R}>\mathcal{R}_0$*, any* $\delta\in(0,\delta_0)$ *and any* $V\in[-\infty,0]^{B_\mathcal{R}}$ *satisfying* $0\in\arg\max(V)$*, the following implication holds:*

$$[\mathcal{L}_{B_\mathcal{R}}(V)\leq 1 \text{ and } \lambda_{B_\mathcal{R}}(V)>-\chi_\mathcal{R}(\varrho)-3\delta]\quad\Longrightarrow\quad d_R(V,V_\varrho)<\frac{\gamma}{2}. \tag{3.25}$$

Proof. According to Lemma 3.2(i), we can choose $\delta>0$ so small that $d_R(V,V_\varrho)<\gamma/2$ for any $V\in[-\infty,0]^{\mathbb{Z}^d}$ satisfying $0\in\arg\max(V)$, $\mathcal{L}(V)\leq 1$ and $\lambda(V)>-\chi(\varrho)-4\delta$. According to Lemma 3.3(i) we may choose $\mathcal{R}_0>0$ so large that $-\chi_{\mathcal{R}_0}(\varrho)\geq-\chi(\varrho)-\delta$.

Now let $\mathcal{R}>\mathcal{R}_0\vee R$ and let $V$ be in $[-\infty,0]^{B_\mathcal{R}}$ with $0\in\arg\max(V)$, $\mathcal{L}_{B_\mathcal{R}}(V)\leq 1$ and $\lambda_{B_\mathcal{R}}(V)>-\chi_\mathcal{R}(\varrho)-3\delta$. Consider $\widetilde V\in[-\infty,0]^{\mathbb{Z}^d}$ given by $\widetilde V=V$ on $B_\mathcal{R}$ and $\widetilde V=-\infty$ on $B_\mathcal{R}^c$. Then we have $\mathcal{L}(\widetilde V)\leq 1$, $0\in\arg\max(\widetilde V)$, and $\lambda(\widetilde V)=\lambda_{B_\mathcal{R}}(V)>-\chi_\mathcal{R}(\varrho)-3\delta\geq-\chi(\varrho)-4\delta$ by our choice of $\mathcal{R}_0$. By our choice of $\delta$, we may conclude that $d_R(V,V_\varrho)=d_R(\widetilde V,V_\varrho)<\gamma/2$. Hence, $\delta$ and $\mathcal{R}_0$ possess the claimed property. □

**4. Spectral bounds.** In this section we derive crucial estimates for $u_3$ which will enable us to reduce the bound in (2.8) to the investigation of the spectral properties of the Hamiltonian $\Delta+\xi$. This result provides the key idea for the construction of the random time-dependent set $\Gamma_{t\log^2 t}$ which will be carried out in Section 5.

Randomness is of no relevance in this section. Fix a box $B\subset\mathbb{Z}^d$ containing the origin, a potential $V:B\to\mathbb{R}$, and a nonempty set $\Gamma\subset B$ arbitrarily. Given $y\in\Gamma$, we denote by $\lambda_y$ and $v_y$ the principal eigenvalue and corresponding positive eigenfunction of $\Delta+V$ in $(B\setminus\Gamma)\cup\{y\}$ with zero boundary condition. We assume that $v_y$ is normalized to $v_y(y) = 1$ rather than in $\ell^2$-sense.

We consider the function $w$ given by the Feynman–Kac formula

$$w(t,x) = \mathbb{E}_x \exp\left\{\int_0^t V(X_s)\,ds\right\}\delta_0(X_t)\mathbb{1}\{\tau_{B^c} > t\}\mathbb{1}\{\tau_\Gamma \le t\}, \qquad t \ge 0, x \in B, \tag{4.1}$$

where we recall that $\tau_A$ is the entrance time into a set $A \subset \mathbb{Z}^d$.

THEOREM 4.1. *For any $t > 0$,*

$$w(t,\cdot) \le \sum_{y\in\Gamma} w(t,y)\|v_y\|_2^2 v_y(\cdot). \tag{4.2}$$

*In particular, for any $r \ge 0$ and $t > 0$,*

$$\frac{\sum_{x\in B\setminus B_r(\Gamma)} w(t,x)}{\sum_{x\in B} w(t,x)} \le \max_{y\in\Gamma}\left[\|v_y\|_2^2 \sum_{x\in B\setminus B_r(\Gamma)} v_y(x)\right]. \tag{4.3}$$

In order to further bound the expression on the right, we may use the following probabilistic representation of the eigenfunction $v_y$:

$$v_y(x) = \mathbb{E}_x \exp\left\{\int_0^{\tau_y} [V(X_s) - \lambda_y]\,ds\right\}\mathbb{1}\{\tau_y = \tau_\Gamma < \tau_{B^c}\}, \qquad y \in \Gamma, x \in B. \tag{4.4}$$

Expectations of this kind can be estimated with the help of the following lemma:

LEMMA 4.2. *For any finite set $A \subset \mathbb{Z}^d$, any potential $V : A \to \mathbb{R}$ and any number $\gamma > \lambda(A)$,*

$$\mathbb{E}_x \exp\left\{\int_0^{\tau_{A^c}} [V(X(s)) - \gamma]\,ds\right\} \le 1 + 2d\frac{|A|}{\gamma - \lambda(A)}, \qquad x \in A, \tag{4.5}$$

*where $\lambda(A)$ denotes the principal Dirichlet eigenvalue of $\Delta + V$ in $A$ with zero boundary condition.*

Later on we choose $B = B_{t\log^2 t}$, $r$ as defined in (1.17), $\Gamma$ the time-dependent random set $\Gamma_{t\log^2 t}$ constructed in Section 5, and $V = \xi$. Then $w(t,\cdot)$ coincides with $u_3(t,\cdot)$ defined in (2.5). Hence, our main task in proving the crucial estimate in (2.8) will consist in controlling the tails of the eigenfunctions $v_y$ uniformly in $y \in \Gamma_{t\log^2 t}$. Indeed, we will show a uniform exponential decay of these eigenfunctions away from their centers; see Proposition 6.1 below. In order to achieve this, we shall use (4.4) and Lemma 4.2 for $\gamma = \lambda_y$ and appropriate sets $A$. To this end, we will need lower bounds for the "spectral gap" $\lambda_y - \lambda(A)$, which will be derived in Lemma 5.4.

Now we turn to the proofs of Theorem 4.1 and Lemma 4.2. In order to prove (4.2), we first need a technical assertion.

LEMMA 4.3. *For $0 < s < t$ and any $y \in \Gamma$,*

$$\mathbb{E}_y \exp\left\{\int_0^{t-s} V(X_u)\,du\right\}\delta_0(X_{t-s})\mathbb{1}\{\tau_{B^c} > t-s\} \le e^{-\lambda_y s}\|v_y\|_2^2 w(t,y). \tag{4.6}$$

PROOF. We obtain a lower bound for $w(t,y)$ by requiring that the random walker is in $y$ at time $s$ and has not entered $\Gamma \setminus \{y\}$ before. Using the Markov property at time $s$, we obtain

$$\begin{aligned} w(t,y) \ge {} & \mathbb{E}_y \exp\left\{\int_0^{s} V(X_u)\,du\right\}\delta_y(X_s)\mathbb{1}\{\tau_{B^c} > s\}\mathbb{1}\{\tau_{\Gamma\setminus\{y\}} > s\} \\ & \times \mathbb{E}_y \exp\left\{\int_0^{t-s} V(X_u)\,du\right\}\delta_0(X_{t-s})\mathbb{1}\{\tau_{B^c} > t-s\}. \end{aligned} \tag{4.7}$$

Using an eigenvalue expansion for the parabolic problem in $(B \setminus \Gamma) \cup \{y\}$ representing the first factor on the right-hand side of (4.7), one obtains the bound

$$\begin{aligned} \mathbb{E}_y \exp\left\{\int_0^{s} V(X_u)\,du\right\}\delta_y(X_s)\mathbb{1}\{\tau_{B^c} > s\}\mathbb{1}\{\tau_{\Gamma\setminus\{y\}} > s\} & \ge e^{\lambda_y s}\frac{v_y(y)^2}{\|v_y\|_2^2} \\ & = e^{\lambda_y s}\|v_y\|_2^{-2}. \end{aligned} \tag{4.8}$$

[Recall that we normed $v_y$ by $v_y(y) = 1$ rather than in $\ell_2$-sense.] Now combine the two estimates to arrive at the assertion. □

PROOF OF THEOREM 4.1. Clearly, (4.3) follows from (4.2) by summing over $x \in B \setminus B_r(\Gamma)$ and estimating elementarily.

Let us turn to the proof of (4.2). Fix $x \in B \setminus \Gamma$. In the Feynman–Kac formula for $w$ in (4.1), we sum over the entrance points in $\Gamma$ and use the strong Markov property at time $\tau_\Gamma$ to obtain

$$\begin{aligned} w(t,x) = {} & \sum_{y\in\Gamma} \mathbb{E}_x \exp\left\{\int_0^{\tau_\Gamma} V(X_u)\,du\right\}\mathbb{1}\{\tau_{B^c} > \tau_\Gamma\}\mathbb{1}\{\tau_\Gamma \le t\}\delta_y(X_{\tau_\Gamma}) \\ & \times \left[\mathbb{E}_y \exp\left\{\int_0^{t-s} V(X_u)\,du\right\}\delta_0(X_{t-s})\mathbb{1}\{\tau_{B^c} > t-s\}\right]_{s=\tau_\Gamma}. \end{aligned} \tag{4.9}$$

Now use Lemma 4.3 and observe that we may replace $\tau_\Gamma$ by $\tau_y$, the first entrance time into $\{y\}$, and $\delta_y(X_{\tau_\Gamma})$ by $v_y(X_{\tau_y})\mathbb{1}\{\tau_y < \tau_{\Gamma\setminus\{y\}}\}$ to obtain

$$w(t,x) \le \sum_{y\in\Gamma} w(t,y)\|v_y\|_2^2 \mathbb{E}_x \exp\left\{\int_0^{\tau_y} [V(X_u) - \lambda_y]\,du\right\}$$

(4.10)
$$\times\, v_y(X_{\tau_y})\mathbb{1}\{\tau_y < \tau_{\Gamma\setminus\{y\}}\}.$$

By the eigenvalue relation at the stopping time $\tau_y$, the latter expectation equals $v_y(x)$. This yields (4.2). □

Proof of Lemma 4.2. This is essentially taken from the proof of Lemma 2.18 in [8]. Denote the left-hand side of (4.5) by $1+v(x)$, which is finite since $\gamma > \lambda(A)$. Then $v$ is the solution of the boundary value problem

$$[\Delta + V - \gamma]v = \gamma - V$$

in $A$ with zero boundary condition on $A^{\rm c}$. Hence,

$$v = \mathcal{R}_\gamma(V-\gamma),$$

where $\mathcal{R}_\gamma$ is the resolvent of $\Delta + V(\cdot)$ in $A$ with zero boundary condition. Since $V-\gamma \le \lambda(A) + 2d - \gamma \le 2d$ on $A$, one derives

$$v(x) \le 2d\mathcal{R}_\gamma\mathbb{1}(x) \le 2d(\mathcal{R}_\gamma\mathbb{1},\mathbb{1})_A \le 2d\frac{|A|}{\gamma-\lambda(A)}, \qquad x\in A,$$

where $(\cdot,\cdot)_A$ denotes the inner product in $\ell^2(A)$, and the last estimate follows from a Fourier expansion of the resolvent. □

## 5. Construction of $\Gamma$.

In this section we introduce the random $t$-dependent set $\Gamma = \Gamma_{t\log^2 t}$ to which we want to apply Theorem 4.1 and with which we shall later prove Proposition 2.2. We are going to define $\Gamma_t$ and switch from $t$ to $t\log^2 t$ in Section 7 only. The actual definition of $\Gamma_t$ appears in Section 5.1 in (5.11). Properties of the set $\Gamma$ are derived in Sections 5.2 and 5.3. Comments on the construction can be found in Section 5.4. Finally, Section 5.5 contains a technical proof.

### 5.1. *Definition of $\Gamma_t$.*

Let $\varrho \in (0,\infty]$ be fixed, and let Assumptions (F) and (M) be satisfied with this $\varrho$. We abbreviate $\chi = \chi(\varrho)$ for the quantity defined in (1.15). Let, furthermore, $a$ be a fixed positive number.

Recall that $B_t = [-t,t]^d \cap \mathbb{Z}^d$ is the centered box of side length $2t$, and that $h_t = \max_{B_t}\xi$ is the maximal value of the random field $\xi$ in $B_t$. Introduce the *set of high exceedances* of the field $\xi$ in $B_t$,

(5.1)
$$Z^{(t)} = Z^{(t)}(\xi) = \{x \in B_t : \xi(x) > h_t - \chi - a\}.$$

Fix $\mathcal{R} > 0$ arbitrarily. Decompose $Z^{(t)}$ into its $2\mathcal{R}$-connected components [two points $x,y$ in $\mathbb{Z}^d$ are called *$r$-neighbors* if $y \in B_r(x)$, and a subset of $\mathbb{Z}^d$ is called *$r$-connected* if any two points in that set may be connected by a path of $r$-neighbors in that set] called *$\mathcal{R}$-islands*. Denote by $Z^{(t)}_{\mathcal{R}}[z]$ the $2\mathcal{R}$-connected component of $Z^{(t)}$ that contains $z$. We put $Z^{(t)}_{\mathcal{R}}[z] = \varnothing$ if $z$ is not in $Z^{(t)}$. The neighborhoods $B_{\mathcal{R}}(Z^{(t)}_{\mathcal{R}}[z])$ are connected in the nearest-neighbor sense and pairwise disjoint.

According to Corollary 2.10 in [8], with probability one, for all sufficiently large $t$, each $\mathcal{R}$-island has no more than

$$K = \lfloor e^{(\chi+a)/\varrho} \rfloor \tag{5.2}$$

elements. Note that $K$ does not depend on $\mathcal{R}$. In each $\mathcal{R}$-island we pick one site with maximal value of the potential $\xi$, and we call this site the *capital* of the $\mathcal{R}$-island. Denote by

$$\mathcal{C}^{(t)}_{\mathcal{R}} = \{z \in Z^{(t)} : z \text{ is the capital of an } \mathcal{R}\text{-island}\} \tag{5.3}$$

the set of capitals.

Let us introduce the terminology of (spectral) *optimality* of a set. We use the abbreviation $\lambda^{(t)}(A) = \lambda_{A\cap B_t}(\xi)$ for the principal eigenvalue of $\Delta + \xi$ in $A \cap B_t$ with zero boundary condition for a finite set $A \subset \mathbb{Z}^d$. Given $t > 0$, a threshold $\delta > 0$ and a radius $\mathcal{R} > 0$, we say that a set

$$A \subset B_t \text{ is } (\delta, \mathcal{R})\text{-optimal} \quad \Longleftrightarrow \quad \lambda^{(t)}(B_{\mathcal{R}}(A)) > h_t - \chi - \delta. \tag{5.4}$$

We denote by

$$\mathcal{C}^{(t)}_{\delta,\mathcal{R}} = \{z \in \mathcal{C}^{(t)}_{\mathcal{R}} : Z^{(t)}_{\mathcal{R}}[z] \text{ is } (\delta, \mathcal{R})\text{-optimal}\} \tag{5.5}$$

the set of capitals whose $\mathcal{R}$-island is $(\delta, \mathcal{R})$-optimal. By

$$Z^{(t)}_{\delta,\mathcal{R}} = \bigcup_{z \in \mathcal{C}^{(t)}_{\delta,\mathcal{R}}} Z^{(t)}_{\mathcal{R}}[z], \tag{5.6}$$

we denote the union of all $(\delta, \mathcal{R})$-optimal $\mathcal{R}$-islands. We next introduce the minimal distance between these islands:

$$\mathcal{D}^{(t)}_{\delta,\mathcal{R}} = \min\{\operatorname{dist}(Z^{(t)}_{\mathcal{R}}[z], Z^{(t)}_{\mathcal{R}}[\tilde{z}]) : z, \tilde{z} \in \mathcal{C}^{(t)}_{\delta,\mathcal{R}} \text{ and } Z^{(t)}_{\mathcal{R}}[z] \neq Z^{(t)}_{\mathcal{R}}[\tilde{z}]\}, \tag{5.7}$$

where dist refers to the lattice distance. It turns out that this distance grows rather fast:

LEMMA 5.1. *For any $\delta \in (0, \varrho \log 2)$ and any $\mathcal{R} > 0$, almost surely as $t \to \infty$,*

$$\mathcal{D}^{(t)}_{\delta,\mathcal{R}} \geq \mathtt{d}^{\delta}_t t^{-o(1)} \qquad \textit{where } \mathtt{d}^{\delta}_t = t^{2e^{-\delta/\varrho}-1}. \tag{5.8}$$

The proof will be given in Section 5.5 below.

Hence, the distance between $(\delta, \mathcal{R})$-optimal $\mathcal{R}$-islands grows like a power of $t$ which can be made arbitrarily close to one by choosing $\delta$ sufficiently small. Note that this growth rate does not depend on the radius $\mathcal{R}$.

We are going to define certain $t$-dependent large neighborhoods of the capitals of all the $(\delta,\mathcal{R})$-optimal $\mathcal{R}$-islands. We abbreviate

$$q = \frac{2d}{2d+a/2} \in (0,1) \quad \text{and} \quad \mathfrak{R}_t = \log^2 t \quad \text{and} \quad \mathfrak{g}_t^0 = t^{-2d}. \tag{5.9}$$

Then $\mathcal{R} \ll \mathfrak{R}_t \ll \mathsf{d}_t^\delta$, hence, the $\mathfrak{R}_t$-neighborhoods around the sites of $\mathcal{C}_{\delta,\mathcal{R}}^{(t)}$ do not intersect each other and have even a large distance to each other. [Our choice of $\mathfrak{R}_t$ is for definiteness only; e.g., any $\mathfrak{R}_t$ satisfying $\log t = o(\mathfrak{R}_t)$ and $\mathfrak{R}_t = t^{o(t)}$ would also work.] We consider the eigenvalues $\lambda^{(t)}(B_{\mathfrak{R}_t}(z))$ with $z \in \mathcal{C}_{\delta,\mathcal{R}}^{(t)}$. Let $I_t^{\mathrm{gap}}$ be the largest subinterval of $[h_t - \chi - \delta/2, h_t - \chi - \delta/4]$ that contains no eigenvalue $\lambda^{(t)}(B_{\mathfrak{R}_t}(z))$, $z \in \mathcal{C}_{\delta,\mathcal{R}}^{(t)}$. We shall refer to $I_t^{\mathrm{gap}}$ as to the *spectral gap*. Let

$$\mathfrak{g} = \mathfrak{g}_t(\xi;\delta,\mathcal{R}) = |I_t^{\mathrm{gap}}| \tag{5.10}$$

denote its length. Now we finally define the set $\Gamma$ by

$$\Gamma = \Gamma_t(\xi;\delta,\mathcal{R}) = \{z \in \mathcal{C}_{\delta,\mathcal{R}}^{(t)} : \lambda^{(t)}(B_{\mathfrak{R}_t}(z)) \geq \sup I_t^{\mathrm{gap}}\}. \tag{5.11}$$

In words, $\Gamma_t$ is the set of capitals $z$ of those $(\delta,\mathcal{R})$-optimal $\mathcal{R}$-islands $Z_{\mathcal{R}}^{(t)}[z]$ such that the eigenvalue $\lambda^{(t)}(B_{\mathfrak{R}_t}(z))$ lies above the spectral gap. Two obvious properties of $\Gamma_t$ are the following. By Lemma 5.1,

$$\min_{y,\widetilde{y}\in\Gamma_t : y\neq\widetilde{y}} |y-\widetilde{y}| \geq \mathcal{D}_{\delta,\mathcal{R}}^{(t)} \geq t^{2e^{-\delta/\varrho}-1-o(1)} \qquad \text{as } t\to\infty. \tag{5.12}$$

By construction of the spectral gap,

$$\min_{z\in\Gamma_t} \lambda^{(t)}(B_{\mathfrak{R}_t}(z)) - \max_{z\in\mathcal{C}_{\delta,\mathcal{R}}^{(t)}\setminus\Gamma_t} \lambda^{(t)}(B_{\mathfrak{R}_t}(z)) \geq \mathfrak{g}_t. \tag{5.13}$$

It turns out that $\mathfrak{g}_t^0 = t^{-2d}$ is a lower bound for the size of the spectral gap:

LEMMA 5.2. *Let $\delta$, $\mathcal{R}$ and $\mathfrak{R}_t$ be as above. Then, with probability one, $\mathfrak{g}_t \geq \mathfrak{g}_t^0$ for $t$ large.*

PROOF. The maximal distance between adjacent eigenvalues $\lambda^{(t)}(B_{\mathfrak{R}_t}(z))$, $z \in \mathcal{C}_{\delta,\mathcal{R}}^{(t)}$, in the interval $[h_t - \chi - \delta/2, h_t - \chi - \delta/4]$ is at least $\frac{\delta}{4}/(1+|\mathcal{C}_{\delta,\mathcal{R}}^{(t)}|)$. As a consequence of Lemma 5.1, the set $\mathcal{C}_{\delta,\mathcal{R}}^{(t)}$ has no more than $(t^{1+o(1)}/\mathsf{d}_t^\delta)^d$ elements. Hence, the assertion follows from the definition of $\mathfrak{g}_t^0$. $\square$

5.2. *Some properties of* $\Gamma_t$. Our construction of $\Gamma_t$ relies on the quantities

$$\varrho \in (0,\infty], \qquad \chi \in (0,2d], \qquad a>0, \qquad K \in \mathbb{N}, \tag{5.14}$$

which we regard as fixed, and on the parameters

$$\delta>0, \qquad \mathcal{R}>0, \qquad t>0, \tag{5.15}$$

which will be chosen appropriately in the sequel, that is, small enough respectively large enough, depending on the requirements of Proposition 2.2. The next lemma shows, in particular, that conditions (ii′) and (iii′) of that proposition can be met.

LEMMA 5.3. (i) *Let* $R>0$ *and* $\gamma>0$ *be given. Choose* $\delta_0>0$ *and* $\mathcal{R}_0>0$ *in accordance with Corollary 3.4 with* $KR$ *instead of* $R$. *Then, for any* $\delta \in (0,\delta_0)$ *and* $\mathcal{R}>\mathcal{R}_0$, *with probability one, for* $t$ *sufficiently large,* $d_{KR}(\xi(y+\cdot)-h_t, V_\varrho(\cdot))<\gamma$ *for any* $y \in \mathcal{C}^{(t)}_{\delta,\mathcal{R}}$. *In particular, the set* $\Gamma=\Gamma_t(\xi;\delta,\mathcal{R})$ *satisfies condition* (iii′) *of Proposition 2.2 with probability one for all* $t$ *sufficiently large.*

(ii) *Let* $\eta>0$ *be given. Then one may choose* $\delta>0$ *so small that, for any* $\mathcal{R}>0$, $\Gamma=\Gamma_t(\xi;\delta,\mathcal{R})$ *satisfies condition* (ii′) *of Proposition 2.2 with probability one for all* $t$ *sufficiently large.*

(iii) *For any* $\delta>0$, *one may choose* $\mathcal{R}>0$ *so large that, with probability one,* $\Gamma=\Gamma_t(\xi;\delta,\mathcal{R})$ *is not empty for* $t$ *large.*

PROOF. Recall from (1.11) that $h_t=\psi(d\log t)+o(1)$ almost surely as $t\to\infty$.

(i) Put $\eta=\frac{\gamma}{2}\wedge 2\delta$. We know from Corollary 2.12 in [8] that, with probability one, for $t$ large,

$$\max_{z\in Z^{(t)}_{\delta,\mathcal{R}}} \mathcal{L}_{B_{2K\mathcal{R}}(z)}(\xi-h_t-\eta)\le 1.$$

It follows from the definition of $Z^{(t)}_{\delta,\mathcal{R}}$ and $B_{\mathcal{R}}(Z^{(t)}_{\mathcal{R}}[z])$ that

$$\min_{z\in Z^{(t)}_{\delta,\mathcal{R}}} \lambda^{(t)}(B_{2K\mathcal{R}}(z))-h_t-\eta>-\chi-3\delta\ge -\chi_{2K\mathcal{R}}-3\delta.$$

Hence, an application of Corollary 3.4 with $KR$ instead of $R$ implies that $d_{KR}(\xi(y+\cdot)-h_t-\eta, V_\varrho(\cdot))<\gamma/2$. Since $\eta<\gamma/2$, the assertion follows.

(ii) The first part of (2.9) follows from the second, and the second is immediate from (5.12).

(iii) According to Corollary 2.19 in [8] one may choose $\mathcal{R}$ so large that, with probability one, for large $t$, there exists a (random and $t$-dependent) $x \in B_t$ such that

$$\lambda^{(t)}(B_{\mathcal{R}/2}(x)) > h_t - \chi - \left(a \wedge \frac{\delta}{4}\right).$$

This implies that $B_{\mathcal{R}/2}(x)$ intersects $Z^{(t)}$ since $\max_{B_{\mathcal{R}/2}(x)} \xi \geq \lambda^{(t)}(B_{\mathcal{R}/2}(x))$. Hence, there exists a $z \in \mathcal{C}_{\mathcal{R}}^{(t)}$ such that the $\mathcal{R}$-island $Z_{\mathcal{R}}^{(t)}[z]$ intersects $B_{\mathcal{R}/2}(x)$. Therefore, $B_{\mathcal{R}}(Z_{\mathcal{R}}^{(t)}[z])$ contains $B_{\mathcal{R}/2}(x)$ and, consequently,

$$\lambda^{(t)}(B_{\mathcal{R}}(Z_{\mathcal{R}}^{(t)}[z])) > h_t - \chi - \frac{\delta}{4}.$$

This clearly implies that $Z_{\mathcal{R}}^{(t)}[z]$ is $(\delta, \mathcal{R})$-optimal, and $\lambda^{(t)}(B_{\mathfrak{R}_t}(z))$ lies above the spectral gap. Hence, $z \in \Gamma_t$, and we are done. □

5.3. *Spectral properties of* $\Gamma$. The fact that also condition (i′) of Proposition 2.2 is satisfied under appropriate choice of the parameters will be proved in Section 7.1 below. An important preparation is presented in the following lemma. Analogously to Section 4, $\lambda_y$ denotes the principal Dirichlet eigenvalue of $\Delta + \xi$ in $(B_t \setminus \Gamma) \cup \{y\}$, for $y \in \Gamma$. As we have indicated in Section 4, it will be crucial in Section 6 to have lower bounds for the gaps between the eigenvalues $\lambda_y$ with $y \in \Gamma$ and the principal eigenvalues in certain neighborhoods of certain islands. These bounds are provided in the following lemma. Recall that $\mathfrak{g}_t^0 = t^{-2d}$, that the optimal potential shape $V_\varrho$ was introduced in Assumption (M), and that $\dot{V}_\varrho$ is $V_\varrho$ dotted (i.e., put equal to $-\infty$) at the origin.

LEMMA 5.4. *Put $b = \varrho \log \frac{8}{7}$ if $\varrho < \infty$ and $b = 1$ if $\varrho = \infty$. One may choose first $R$ sufficiently large, then $\delta > 0$ sufficiently small, and afterward $\mathcal{R} > R$ large enough so that the following is true with probability one. For $t$ sufficiently large and any $y \in \Gamma_t$,*

$$\lambda_y - \max_{z \in Z_{\delta,\mathcal{R}}^{(t)}} \lambda(B_R(Z_R^{(t)}[z]) \setminus \mathcal{C}_{\delta,\mathcal{R}}^{(t)}) \geq \frac{1}{4}(\lambda(V_\varrho) - \lambda(\dot{V}_\varrho)) \vee \frac{b}{2}, \tag{5.16}$$

$$\lambda_y - \max_{z \in \mathcal{C}_{\mathcal{R}}^{(t)} \setminus \mathcal{C}_{\delta,\mathcal{R}}^{(t)}} \lambda(B_{\mathcal{R}}(Z_{\mathcal{R}}^{(t)}[z])) \geq \delta/2, \tag{5.17}$$

$$\lambda_y - \max_{z \in \mathcal{C}_{\delta,\mathcal{R}}^{(t)} \setminus \Gamma_t} \lambda(B_{\mathfrak{R}_t}(z)) \geq \mathfrak{g}_t^0. \tag{5.18}$$

PROOF. Assertion (5.18) follows from (5.13) and the observation that $\lambda_y \geq \lambda(B_{\mathfrak{R}_t}(y))$ by construction and $\mathfrak{g}_t \geq \mathfrak{g}_t^0$ by Lemma 5.2.

To prove (5.17), observe that

$$\lambda_y \ge \lambda(B_{\mathfrak{R}_t}(y)) \ge \sup I_t^{\rm gap} \ge h_t - \chi - \frac{\delta}{2}, \qquad y \in \Gamma_t. \tag{5.19}$$

On the other hand, by definition of $\mathcal{C}^{(t)}_{\delta,\mathcal{R}}$ and $(\delta,\mathcal{R})$-optimality, we have

$$\lambda(B_{\mathcal{R}}(Z^{(t)}_{\mathcal{R}}[z])) \le h_t - \chi - \delta, \qquad z \in \mathcal{C}^{(t)}_{\mathcal{R}} \setminus \mathcal{C}^{(t)}_{\delta,\mathcal{R}}. \tag{5.20}$$

Combining the two estimates, we arrive at (5.17).

In order to derive (5.16), we choose $R$ so large that $\lambda_{B_R}(V_\varrho) \ge \lambda(V_\varrho) - b/4$. Further, we choose an auxiliary parameter $\gamma > 0$ so small that $\gamma < b/4$ and such that the following implication holds for any $V : B_{KR} \to \mathbb{R}$:

$$\begin{aligned} &\max_{B_{KR}} |V - V_\varrho| < \gamma \\ &\qquad \Longrightarrow \quad [|\lambda_A(V) - \lambda_A(V_\varrho)| \le \tfrac{1}{2}(\lambda(V_\varrho) - \lambda(\dot V_\varrho)), \text{ for all } A \subset B_{KR}]. \end{aligned} \tag{5.21}$$

This may be achieved by using the continuity of the eigenvalue $\lambda_A(\cdot)$.

Next, we require that $\delta \in (0, \delta_0)$ and $\mathcal{R} > \mathcal{R}_0$, where $\delta_0, \mathcal{R}_0$ are chosen in accordance with Lemma 5.3, and $\delta < \frac{1}{2}[(\lambda(V_\varrho) - \lambda(\dot V_\varrho)) \wedge b]$. According to Lemma 5.3, we have that $d_R(\xi(z + \cdot) - h_t, V_\varrho)) < \gamma$ for any $z \in Z^{(t)}_{\delta,\mathcal{R}}$, if $t$ is sufficiently large. Now let $t$ additionally be so large that $\mathcal{D}^{(t)}_{b,R} > 2\mathcal{R}$ (recall Lemma 5.1). Pick $z \in Z^{(t)}_{\delta,\mathcal{R}}$. We shall show that

$$\begin{aligned} &\lambda(B_R(Z^{(t)}_R[z]) \setminus \mathcal{C}^{(t)}_{\delta,\mathcal{R}}) \\ &\qquad \le \begin{cases} h_t - \chi - b, & \text{if } Z^{(t)}_R[z] \cap \mathcal{C}^{(t)}_{\delta,\mathcal{R}} = \varnothing, \\ h_t - \chi - \frac{1}{2}(\lambda(V_\varrho) - \lambda(\dot V_\varrho)), & \text{otherwise.} \end{cases} \end{aligned} \tag{5.22}$$

[Certainly, (5.22) implies (5.16) because of (5.19).] First assume that the $R$-island $Z^{(t)}_R[z]$ does not contain the capital $\widetilde z$ of the $\mathcal{R}$-island $Z^{(t)}_{\mathcal{R}}[z]$, that is, $B_R(Z^{(t)}_R[z])$ and $B_R(Z^{(t)}_R[\widetilde z])$ are disjoint, and $B_R(Z^{(t)}_R[z]) \setminus \mathcal{C}^{(t)}_{\delta,\mathcal{R}} = B_R(Z^{(t)}_R[z])$. We now show that the $R$-island $Z^{(t)}_R[\widetilde z]$ is $(b,R)$-optimal in the sense of definition (5.5). Then $Z^{(t)}_R[z]$ cannot be $(b,R)$-optimal as well, since the distance $\mathcal{D}^{(t)}_{b,R}$ between $(b,R)$-optimal $R$-islands is larger than $2\mathcal{R}$, and this implies the first line of (5.22).

To show the $(b,R)$-optimality of $Z^{(t)}_R[\widetilde z]$, recall that $d_R(\xi(\widetilde z + \cdot) - h_t, V_\varrho) < \gamma$ (because $\widetilde z \in \mathcal{C}^{(t)}_{\delta,\mathcal{R}}$) to estimate

$$\begin{aligned} \lambda(B_R(Z^{(t)}_R[\widetilde z])) &\ge \lambda(B_R(\widetilde z)) \ge h_t + \lambda_{B_R}(V_\varrho) - \gamma \\ &\ge h_t + \lambda(V_\varrho) - \frac{b}{4} - \gamma \ge h_t - \chi - \frac{b}{2}. \end{aligned} \tag{5.23}$$

Here we have also used that $\lambda(V_\varrho) = -\chi$ and that $\gamma < b/4$. This shows the $(b, R)$-optimality of $Z_R^{(t)}[\widetilde{z}]$ and ends the proof of the first line of (5.22).

Turning to the second case, we may assume that $z \in \mathcal{C}_{\delta,\mathcal{R}}^{(t)}$. Hence, $B_R(Z_R^{(t)}[z]) \setminus \mathcal{C}_{\delta,\mathcal{R}}^{(t)} = B_R(Z_R^{(t)}[z]) \setminus \{z\}$. We apply the implication in (5.21) for $V = \xi(z + \cdot) - h_t$ and $A = [B_R(Z_R^{(t)}[z]) \setminus \{z\}] - z$. [The assumption in (5.21) is satisfied by Lemma 5.3.] This implies

$$\begin{aligned}\lambda(B_R(Z_R^{(t)}[z]) \setminus \{z\}) &\le h_t + \lambda_{B_{KR}}(\dot{V}_\varrho) + \tfrac{1}{2}(\lambda(V_\varrho) - \lambda(\dot{V}_\varrho)) \\ &\le h_t - \chi - \tfrac{1}{2}(\lambda(V_\varrho) - \lambda(\dot{V}_\varrho)),\end{aligned}$$

where we recall that $-\chi = \lambda(V_\varrho)$. This implies the second line in (5.22) and ends the proof.

5.4. *Informal description.* Let us repeat in words what properties the field $\xi$ and the set $\Gamma = \Gamma_t(\xi; \delta, \mathcal{R})$ satisfy almost surely for large $t$, provided that the parameters $\delta$, $\mathcal{R}$ and $t$ are chosen appropriately.

We recall that we regard the quantities in (5.14) as fixed. As in Proposition 2.2, let parameters $\gamma, \eta, R > 0$ be given. We may assume that $R$ is sufficiently large, at least as large as is required in Lemma 5.4. Suppose that the parameters $\delta$ and $\mathcal{R}$ are chosen sufficiently small respectively large, in accordance with Lemmas 5.3 and 5.4. Furthermore, pick $\mathfrak{R}_t = \log^2 t$ as in (5.9). We shall assume, in addition, $R > 0$ sufficiently large, $\delta$ sufficiently small and $\mathcal{R} > 0$ sufficiently large such that certain additional conditions be satisfied which depend on the quantities in (5.14) only.

Then the following assertions hold almost surely if $t$ is sufficiently large. (For the sake of simplicity, we suppress the dependence on $t$ from the notation.)

The set of high exceedances, $Z = \{x \in B : \xi(x) > h - \chi - a\}$, consists of $\mathcal{R}$-islands which split into $R$-islands. Henceforth, we call the $\mathcal{R}$-islands *archipelagos* and the $R$-islands just *islands*. Every archipelago has no more than $K$ elements (and hence, no more than $K$ islands) and contains a capital in which the potential $\xi$ is maximal by definition.

We call

$$\begin{array}{rll} B_R(Z_R[z]) & \text{a } \textit{big cluster} & \text{if } z \in Z, \\ B_{\mathcal{R}}(Z_{\mathcal{R}}[z]) & \text{a } \textit{large cluster} & \text{if } z \in Z, \\ B_{\mathfrak{R}}(z) & \text{a } \textit{huge cluster} & \text{if } z \in \mathcal{C}_{\mathcal{R}}. \end{array} \tag{5.24}$$

The big clusters are disjoint, and the large clusters as well. Any large cluster contains no more than $K$ big clusters. We call the big clusters $(b, R)$-optimal and the large clusters $(\delta, \mathcal{R})$-optimal if the island, respectively, the archipelago, that forms the cluster has this property. The $(b, R)$-optimal big clusters have distance $\ge t^{3/4 - o(1)}$, and the $(\delta, \mathcal{R})$-optimal large clusters have an even much larger distance. Since $\mathfrak{R}$ is much smaller than the distance between $(\delta,\mathcal{R})$-optimal large clusters, but much larger than $\mathcal{R}$, any huge cluster $B_{\mathfrak{R}}(z)$ with $z\in\mathcal{C}_{\delta,\mathcal{R}}$ contains precisely one $(\delta,\mathcal{R})$-optimal large cluster, which is the one that contains $z$. Analogously, any $(\delta,\mathcal{R})$-optimal large cluster contains at most one $(b,R)$-optimal big cluster, which then is the one that contains the capital of the archipelago. However, the big, large and huge clusters may lie everywhere in the box $B$, and non-$(\delta,\mathcal{R})$-optimal archipelagos or non-$(b,R)$-optimal islands may even be neighboring; in particular, a huge cluster may contain many (nonoptimal) big and large clusters.

A site $y\in Z$ belongs to $\Gamma$ if and only if the potential $\xi$ is maximal at $y$ within the archipelago (i.e., $y$ is its capital), the large-cluster eigenvalue satisfies the lower bound $\lambda(B_{\mathcal{R}}(Z_{\mathcal{R}}[y]))>h-\chi-\delta$ [i.e., the archipelago of $y$ is $(\delta,\mathcal{R})$-optimal], and the huge-cluster eigenvalue $\lambda(B_{\mathfrak{R}}(y))$ lies above the spectral gap $I^{\mathrm{gap}}$. Since any huge cluster contains at most one point of the set $\Gamma$, we have the bound $\lambda_y\geq\lambda(B_{\mathfrak{R}}(y))$ for any $y\in\Gamma$. It is worth remarking that

$$\min_{y\in\Gamma}\lambda(B_{\mathfrak{R}}(y))-\max_{z\in\mathcal{C}_{\delta,\mathcal{R}}\setminus\Gamma}\lambda(B_{\mathfrak{R}}(z))\geq\mathfrak{g}\geq\mathfrak{g}^0. \tag{5.25}$$

Furthermore, by construction of the set $\Gamma$, if $\delta<a$, we also have that

$$\lambda_y\geq\lambda(B_{\mathfrak{R}}(y))\geq h-\chi-\frac{\delta}{2}\geq h-\chi-\frac{a}{2}. \tag{5.26}$$

An important issue is the eigenvalues associated with the big, large and huge clusters and the gaps in (5.16)–(5.18) between these eigenvalues and the eigenvalues $\lambda_y$ with $y\in\Gamma$. The necessity of *three* types of neighborhoods of points of $Z$ is partially due to our approach in Section 6 below. It may intuitively be explained as follows. In the big clusters around points of $\Gamma$, the potential $\xi$ is required to approximate the optimal shape $h_t+V_\varrho$. The eigenvalue of the surrounding large cluster guarantees this property via $(\delta,\mathcal{R})$-optimality (see Corollary 3.4). The minimal gap between the eigenvalues of any two $(\delta,\mathcal{R})$-optimal large clusters depends on $t$ and shrinks to $0$ as $t\to\infty$. In order to compensate for that, we have to introduce huge clusters whose size depends on $t$.

In order to successfully apply Theorem 4.1, we need to guarantee the following. For each $y\in\Gamma$, the positive eigenfunction $v_y$ corresponding to $\lambda_y$ is concentrated in a neighborhood of $y$. Therefore one needs to avoid "resonances" between the local eigenvalues corresponding to the huge clusters around $y$ and all the ones corresponding to other huge clusters around $(\delta,\mathcal{R})$-optimal capitals not belonging to $\Gamma$, as well as the ones corresponding to *dotted* huge clusters around the other capitals of $\Gamma$. This is expected to be satisfied, provided that the distance $\mathcal{D}_{\delta,\mathcal{R}}$ between such clusters is large enough, depending on the smallness of the spectral gap $\mathfrak{g}$. This will turn out to be guaranteed by the choice $\mathfrak{R}_t=\log^2 t$, thanks to Lemmas 5.2 and 5.3.

5.5. *Proof of Lemma* 5.1. Fix $\eta > 0$ small and recall from (5.2) [resp. (1.11)] that, with probability one, if $t$ is large enough, $|Z_{\mathcal{R}}^{(t)}[z]| \le K$ for any $z \in \mathcal{C}^{(t)}$, and $h_t \ge \psi(d\log t) - \eta$. Hence, our assertion clearly follows from the following statement: with probability one, for any $t$ sufficiently large, any two $\mathcal{R}$-connected subsets $A$, $\widetilde{A}$ of $B_t$ having not more than $K$ elements each and having distance larger than $\mathcal{R}$ to each other, such that the eigenvalues $\lambda(B_{\mathcal{R}}(A))$ and $\lambda(B_{\mathcal{R}}(\widetilde{A}))$ are both larger than $\psi(d\log t) - \eta - \chi - \delta$, have distance even larger than $\mathtt{d}_t^{\delta} t^{-\zeta}$ to each other. Here $\zeta > 0$ is a small number.

In order to prove this, according to the Borel–Cantelli lemma, it suffices to show the summability over $n \in \mathbb{N}$ of

$$\begin{aligned} p_n = \operatorname{Prob}(\exists\, 2\mathcal{R}\text{-connected sets } A, \widetilde{A} \subset B_{2^{n+1}} \text{ such that} \\ \mathcal{R} < \operatorname{dist}(A, \widetilde{A}) \le \mathtt{d}_{2^{n+1}}^{\delta} 2^{-\zeta(n+1)} \text{ and } |A| \vee |\widetilde{A}| \le K \text{ and} \\ \lambda(B_{\mathcal{R}}(A)) \wedge \lambda(B_{\mathcal{R}}(\widetilde{A})) > \psi(d\log(2^n)) - \eta - \chi - \delta). \end{aligned}$$

This is done as follows. Estimate

$$\begin{aligned} p_n &\le (2^{n+2}+1)^d (2\mathcal{R}+1)^{2Kd} (2\mathtt{d}_{2^{n+1}}^{\delta} 2^{-\zeta(n+1)} + 1)^d \\ &\qquad \times \operatorname{Prob}(\lambda(B_{K\mathcal{R}}) > \psi(d\log(2^n)) - \chi - \delta - \eta)^2 \\ &\le \operatorname{const}(2^n \mathtt{d}_{2^n}^{\delta} 2^{-\zeta n})^d \operatorname{Prob}(\lambda(B_{K\mathcal{R}}) > \psi(d\log(2^n)) - \chi - \delta - \eta)^2. \end{aligned} \tag{5.27}$$

Now let us assume that $\varrho$ is finite. We estimate the tails of $\lambda(B_{K\mathcal{R}})$, with the help of the exponential Chebyshev inequality as follows. For any $\gamma > 0$ and $t > 0$,

$$\begin{aligned} &\operatorname{Prob}(\lambda(B_{K\mathcal{R}}) > \psi(d\log t) - \chi - \delta - \eta) \\ &\qquad \le \exp\{-\gamma[\psi(d\log t) - \chi - \delta - \eta]\} \langle e^{\gamma\lambda(B_{K\mathcal{R}})} \rangle. \end{aligned} \tag{5.28}$$

Now we use a Fourier expansion, in terms of the eigenvalues $\lambda^{(k)}(B_{K\mathcal{R}})$ and the corresponding $\ell^2$-normalized eigenfunctions $v_k$ of $\Delta + \xi$ in $B_{K\mathcal{R}}$, and Parseval's identity to obtain

$$e^{\gamma\lambda(B_{K\mathcal{R}})} \le \sum_{k=0}^{\infty} e^{\gamma\lambda^{(k)}(B_{K\mathcal{R}})} \sum_{x \in B_{K\mathcal{R}}} v_k(x)^2 = \sum_{x \in B_{K\mathcal{R}}} q(\gamma, x, x), \tag{5.29}$$

where $q(t,x,y)$ denotes the fundamental solution of $\partial_t = \Delta + \xi(\cdot)$ in $\mathbb{Z}^d$. Taking expectations, we obtain, as $\gamma \to \infty$,

$$\begin{aligned} \langle e^{\gamma\lambda(B_{K\mathcal{R}})} \rangle &\le (2K\mathcal{R}+1)^d \langle q(\gamma, 0, 0) \rangle \\ &\le \langle e^{\gamma\xi(0)} \rangle e^{-\gamma\chi + o(\gamma)} \\ &\le \exp\{\gamma[\psi(\gamma) - \chi + \varrho\log\varrho - \varrho + o(1)]\}, \end{aligned} \tag{5.30}$$

where the second estimate is taken from Theorem 1.2 (see Remark 1.3.a) in [8] with $p=1$ [recall that our $\chi$ is identical to $2d\chi(\varrho)$ there], and the last estimate is taken from [8], Lemma 2.3.

We now use (5.30) in (5.28) and substitute $\gamma = \frac{c}{\varrho} d\log t$ with some $c>0$. For $\gamma$ sufficiently large, this gives

$$\begin{aligned}
&\operatorname{Prob}(\lambda(B_{K\mathcal{R}}) > \psi(d\log t) - \chi - \delta - \eta) \\
&\quad \le \exp\{-\gamma[\psi(d\log t) - \delta - \eta - \psi(\gamma) - \varrho\log\varrho + \varrho - \eta]\} \\
&\quad \le \exp\biggl\{-dc\log t\biggl[-\log c - \frac{\delta}{\varrho} + 1 - \frac{3\eta}{\varrho}\biggr]\biggr\},
\end{aligned} \tag{5.31}$$

where we have used also that $\psi(d\log t) - \psi(\frac{c}{\varrho} d\log t) = -\varrho\log\frac{c}{\varrho} + o(1)$, in accordance with (1.7) in our Assumption (F). Now we see that the choice $c = e^{-\delta/\varrho}$ is asymptotically close to optimal and yields the upper bound

$$\begin{aligned}
\operatorname{Prob}(\lambda(B_{K\mathcal{R}}) > \psi(d\log t) - \chi - \delta - \eta) &\le \exp\biggl\{-e^{-\delta/\varrho} d\log t + \frac{3\eta}{\varrho} d\log t\biggr\} \\
&= t^{-d(e^{-\delta/\varrho} - 3\eta/\varrho)}.
\end{aligned}$$

Using this bound with $t = 2^n$ in (5.27), we get that

$$p_n \le \operatorname{const}(\mathsf{d}_{2^n}^{\delta} 2^{-\zeta n})^d (2^n)^{d[1-2e^{-\delta/\varrho} + 6\eta/\varrho]} = \operatorname{const} 2^{-nd(6\eta/\varrho - \zeta)}.$$

For $\eta$ sufficiently small, this is summable over $n \in \mathbb{N}$. This ends the proof of the lemma in the case that $\varrho \in (0,\infty)$.

Let us turn to the case that $\varrho = \infty$. We go back to (5.27) and estimate $\lambda(B_{K\mathcal{R}}) \le \max_{B_{K\mathcal{R}}} \xi$ and recall that $\chi(\infty) = 2d$ to obtain, if $\delta + \eta < 2d$,

$$\begin{aligned}
&\operatorname{Prob}(\lambda(B_{K\mathcal{R}}) > \psi(d\log(2^n)) - \chi - \delta - \eta)^2 \\
&\quad \le |B_{K\mathcal{R}}|^2 \operatorname{Prob}(\xi(0) > \psi(d\log(2^n)) - 4d)^2 \\
&\quad \le \operatorname{const} e^{-2\varphi(\psi(d\log(2^n)) - 4d)}.
\end{aligned} \tag{5.32}$$

Here we recall that $\varphi$ and $\psi$ are defined in (1.5), respectively (1.6). Now we show that, because of Assumption (F), we have

$$\lim_{h\to\infty} [\varphi(h-c) - \vartheta\varphi(h)] = \infty, \qquad c>0, \vartheta \in (0,1). \tag{5.33}$$

In order to prove (5.33), pick $\beta \in (\vartheta, 1)$ arbitrarily and note that $\alpha_s = \psi(s) - \psi(\beta s)$ tends to $\infty$ as $s \to \infty$, according to Assumption (F). Since $\varphi \circ \psi$ is the identity, we have, for $s$ large,

$$\begin{aligned}
\varphi(\psi(s) - c) - \vartheta\varphi(\psi(s)) &\ge \varphi(\psi(s) - \alpha_s) - \vartheta\varphi(\psi(s)) \\
&= \varphi(\psi(\beta s)) - \vartheta s = (\beta - \vartheta)s \to \infty.
\end{aligned} \tag{5.34}$$

Substituting $h = \psi(s)$, this implies that (5.33) holds.

Now we use (5.33) in (5.32) with $\vartheta = 1 - \zeta/4$ and $h = \psi(d\log(2^n))$ and $c = 4d$ and obtain

$$\begin{aligned}
&\mathrm{Prob}(\lambda(B_{K\mathcal{R}}) > \psi(d\log(2^n)) - \chi - \delta - \eta)^2 \\
&\quad \le \mathrm{const}\, e^{-(1-\zeta/2)\varphi(\psi(d\log(2^n)))} \\
&\quad \le \mathrm{const}\, 2^{-nd(1-\zeta/2)}.
\end{aligned} \tag{5.35}$$

Using this in (5.27), we obtain

$$p_n \le \mathrm{const}\, 2^{nd(2-\zeta)}(2^{-nd})^{1-\zeta/2} = \mathrm{const}\, 2^{-nd\zeta/2}, \tag{5.36}$$

and this is summable over $n$. □

**6. Exponential decay of the eigenfunctions.** In this section we prove that, for the set $\Gamma = \Gamma_t(\xi;\delta,\mathcal{R})$ defined in Section 5, the eigenfunctions $v_y$ for $\Delta + \xi$ in $(B_t \setminus \Gamma_t) \cup \{y\}$ introduced in Section 4 decay exponentially away from their centers, uniformly in $y \in \Gamma$ and $t \gg 1$, provided that the parameters $\delta$, $R$ and $\mathcal{R}$ are chosen appropriately. This property is fundamental for estimating the right-hand side of (4.3) for our random Hamiltonian $\Delta + \xi$.

The main result of this section is the following. Recall the quantities in (5.14) and $q = 2d/(2d + a/2) \in (0,1)$.

PROPOSITION 6.1. *One may choose first $R$ sufficiently large, then $\delta > 0$ sufficiently small, and afterward $\mathcal{R}$ large enough so that, with probability one, for all sufficiently large $t$,*

$$v_y(x) \le q^{c|x-y|}, \qquad y \in \Gamma_t(\xi;\delta,\mathcal{R}), x \in B_t \setminus B_R(y), \tag{6.1}$$

*where $c > 0$ is a constant that depends on the quantities in (5.14) only.*

The rest of this section is devoted to the proof of Proposition 6.1. The main tools are probabilistic cluster expansions based on a decomposition of the trajectory of the random walk $X$ in the Feynman–Kac representation of $v_y$ into the subpaths between different archipelagos.

In Section 6.1 we explain the main idea of the proof of Proposition 6.1 in words. The proof of Proposition 6.1 is finished in Section 6.2, subject to three fundamental lemmas, the three types of expansions: for big, large and huge clusters, respectively. In Sections 6.3, 6.4 and 6.5 we prove the three cluster expansions, respectively. Finally, in Section 6.6 we provide a corollary of Proposition 6.1.

6.1. *Heuristic explanation of the proof of Proposition* 6.1. Let us describe the idea of the proof of Proposition 6.1 in words. For simplicity, we shall suppress the dependence on $t$ from the notation. Recall that $X$ is the simple random walk on $\mathbb{Z}^d$, and $\tau_A$ is its entrance time into a set $A$, and $\tau_y = \tau_{\{y\}}$.

Our proof relies on the following probabilistic representation for the eigenfunction $v_y$ [cf. (4.4)]:

$$v_y(x) = \mathbb{E}_x \exp\left\{\int_0^{\tau_y} [\xi(X_s) - \lambda_y]\, ds\right\} \mathbb{1}\{\tau_y = \tau_\Gamma < \tau_{B^c}\}, \qquad y \in \Gamma, x \in B. \tag{6.2}$$

In words, the walker starts at $x \in B$, stays inside the large box $B$, and ends up in $y$ without having visited any other site of $\Gamma$ before.

The main idea that leads to the exponential decay in (6.1) is that the particle has to make at least const $|x-y|$ steps outside the set of high exceedances, $Z$. But in this area, we may estimate $\xi(\cdot) - \lambda_y \leq h - \chi - a - \lambda_y \leq -a/2$, as is seen from the bound in (5.26). Since the first jump time of the random walk is exponentially distributed with mean $1/(2d)$, we have, for any $z \in \mathbb{Z}^d$,

$$\mathbb{E}_z \exp\left\{-\frac{a}{2} \inf\{t > 0 : X_t \neq X_0\}\right\} = \frac{2d}{2d + a/2} = q. \tag{6.3}$$

Hence, for any step outside $Z$, we may estimate the contribution to the expectation in (6.2) by a factor of $q$. Therefore, we obtain

$$\mathbb{E}_{z_0} \exp\left\{\int_0^{\tau_Z} [\xi(X_s) - \lambda_y]\, ds\right\} \mathbb{1}\{\tau_Z = \tau_z < \tau_{B^c}\} \leq q^{|z-z_0|}, \qquad z_0 \in B \setminus Z, z \in Z. \tag{6.4}$$

The big technical difficulty is to control the contribution from trajectories that do not go straight to $y$ through $B \setminus Z$, but take a detour through one or several other archipelagos, where they might gain a larger contribution. However, because of the far distances between the huge clusters $B_{\mathfrak{R}}(\widetilde{y})$ with $\widetilde{y} \in \Gamma$, it will finally turn out that the price for traveling from one to another is too high and spoils the gain from staying inside these favorable regions.

We have to control the contribution from trajectories that visit $Z \setminus \Gamma$ before time $\tau_\Gamma$. We shall distinguish three different subsets of $Z \setminus \Gamma$: (1) $Z_{\delta,\mathcal{R}} \setminus \mathcal{C}_{\delta,\mathcal{R}}$, the union of all islands in $(\delta, \mathcal{R})$-optimal archipelagos, after removing the capitals, (2) $Z \setminus Z_{\delta,\mathcal{R}}$, the union of all non-$(\delta, \mathcal{R})$-optimal archipelagos, and (3) $\mathcal{C}_{\delta,\mathcal{R}} \setminus \Gamma$, the set of optimal capitals whose huge-cluster eigenvalue lies below the spectral gap. Note that $Z \setminus \Gamma$ is the disjoint union of these three sets.

Let us consider the walker on his way within $B$ from $x \in B$ to $y \in \Gamma$. We follow his path until he first visits the set $Z$. The contribution of the path until $\tau_Z$ is controlled using (6.4). If the walker is already at $y$, the journey is finished, and there is nothing to do anymore. The remaining part of the path, that is, the piece from $Z$ to $y$, will be split into several subpaths of three different kinds, which will be controlled by appropriate cluster expansions:

*Big-cluster expansion*: This expansion handles paths that start from $Z_{\delta,\mathcal{R}} \setminus \mathcal{C}_{\delta,\mathcal{R}}$ and end in $Z \setminus Z_{\delta,\mathcal{R}}$ or $\mathcal{C}_{\delta,\mathcal{R}}$ without visiting the union of these two sets before.

*Large-cluster expansion*: This expansion handles paths that start from $Z \setminus Z_{\delta,\mathcal{R}}$ and end at their first visit to $\mathcal{C}_{\delta,\mathcal{R}}$.

*Huge-cluster expansion*: This expansion handles paths that start from $\mathcal{C}_{\delta,\mathcal{R}}$ and end at their first visit to the final destination $y$.

Clearly, the paths considered in the large-cluster expansion may contain subpaths handled in the big-cluster expansion, and the paths in the huge-cluster expansion may contain subpaths in the big-cluster and large-cluster expansions. In all three expansions the main difficulty is to control the contribution coming from the times the walker spends in $Z$. This will be done as follows. If the walker is in $z_0 \in Z_{\delta,\mathcal{R}} \setminus \mathcal{C}_{\delta,\mathcal{R}}$, then we control the contribution until the walker leaves the big cluster $B_R(Z_R[z_0]) \setminus \mathcal{C}_{\delta,\mathcal{R}}$ using Lemma 4.2 in combination with the big-cluster spectral gap bound (5.16). Similarly, if the walker is in $z_0 \in Z \setminus Z_{\delta,\mathcal{R}}$, then we control the contribution until the walker leaves the large cluster $B_{\mathcal{R}}(Z_{\mathcal{R}}[z_0])$ again using Lemma 4.2, together with the large-cluster spectral gap bound (5.17). If the walker is in $z_0 \in \mathcal{C}_{\delta,\mathcal{R}}$, then we use the huge-cluster $B_{\mathfrak{R}}(Z_{\mathcal{R}}[z_0])$ and the corresponding spectral gap bound in (5.18). These contributions will be more than compensated by applications of (6.4) to the path segments outside $Z$, whose total length dominates the cluster contributions.

6.2. *Proof of Proposition* 6.1. The precise formulations of the above arguments are given in Lemmas 6.4, 6.3 and 6.2. Afterward, we finish the proof of Proposition 6.1, subject to the three assertions.

Let us introduce some terminology. For any $z \in Z$, we call $\dot{Z}_{\mathcal{R}}[z] = Z_{\mathcal{R}}[z] \setminus \mathcal{C}$ a *dotted archipelago*, that is, we remove from the archipelago its capital. Analogously, we call $\dot{Z}_R[z] = Z_R[z] \setminus \mathcal{C}$ a *dotted island*. Denote the union of dotted $(\delta,\mathcal{R})$-optimal archipelagos by

$$\dot{Z}_{\delta,\mathcal{R}} = Z_{\delta,\mathcal{R}} \setminus \mathcal{C}_{\delta,\mathcal{R}} = \bigcup_{z\in\mathcal{C}_{\delta,\mathcal{R}}} \dot{Z}_{\mathcal{R}}[z]. \tag{6.5}$$

For notational convenience, we shall abbreviate $\exp\{\int[\xi(X_s) - \lambda_y]\,ds\}$ by $\exp\{\int\}$.

LEMMA 6.2 (Big-cluster expansion). (a) *One can choose $\delta > 0$ small enough and afterward $\mathcal{R} > 0$ large enough, such that, with probability one,* *for any sufficiently large $t$, for any $z_0 \in \dot{Z}_{\delta,\mathcal{R}}$ and $z \in Z \setminus Z_{\delta,\mathcal{R}}$,*

$$\mathbb{E}_{z_0} \exp\left\{\int_0^{\tau_z}\right\} \mathbb{1}\{\tau_z = \tau_{Z\setminus Z_{\delta,\mathcal{R}}} < \tau_{\mathcal{C}_{\delta,\mathcal{R}}} \wedge \tau_{B^c}\} \le q^{|z_0 - z|/(8K)}. \tag{6.6}$$

(b) *There exists $M > 1$, depending on the quantities in (5.14) only, such that one can choose $\delta > 0$ small enough and afterward $\mathcal{R} > 0$ large enough such that, with probability one, for any sufficiently large $t$, for any $z_0 \in \dot{Z}_{\delta,\mathcal{R}}$ and $z \in \mathcal{C}_{\delta,\mathcal{R}}$,*

$$\mathbb{E}_{z_0} \exp\left\{\int_0^{\tau_z}\right\} \mathbb{1}\{\tau_z = \tau_{\mathcal{C}_{\delta,\mathcal{R}}} < \tau_{Z\setminus Z_{\delta,\mathcal{R}}} \wedge \tau_{B^c}\} \le M q^{|z_0 - z|/(8K)}. \tag{6.7}$$

LEMMA 6.3 (Large-cluster expansion). *One can choose $\delta > 0$ small enough and afterward $\mathcal{R} > 0$ large enough, such that, with probability one, for any sufficiently large $t$, for any $z_0 \in Z \setminus Z_{\delta,\mathcal{R}}$ and any $z \in \mathcal{C}_{\delta,\mathcal{R}}$,*

$$\mathbb{E}_{z_0} \exp\left\{\int_0^{\tau_z}\right\} \mathbb{1}\{\tau_z = \tau_{\mathcal{C}_{\delta,\mathcal{R}}} < \tau_{B^c}\} \le q^{|z - z_0|/(256K^2)}. \tag{6.8}$$

LEMMA 6.4 (Huge-cluster expansion). *One can choose $\delta > 0$ small enough and afterward $\mathcal{R} > 0$ large enough, such that, with probability one, for any sufficiently large $t$, for any $z_0 \in \mathcal{C}_{\delta,\mathcal{R}} \setminus \Gamma$ and any $y \in \Gamma$,*

$$\mathbb{E}_{z_0} \exp\left\{\int_0^{\tau_y}\right\} \mathbb{1}\{\tau_y = \tau_\Gamma < \tau_{B^c}\} \le q^{|z_0 - y|/(4096K^2)}. \tag{6.9}$$

The proofs of Lemmas 6.2–6.4 are given in Sections 6.3–6.5, respectively. The proof of Lemma 6.2 is independent of the other two lemmas, while the proof of Lemma 6.3 uses Lemma 6.2, and the proof of Lemma 6.4 uses Lemmas 6.2 and 6.3.

Let us now finish the proof of Proposition 6.1 subject to Lemmas 6.2–6.4. Fix $y \in \Gamma$ and $x \in B \setminus B_R(y)$. In (6.2) we use the Markov property at time $\tau_Z$ and apply (6.4). Distinguishing the two cases whether or not the walker is already at the site $y$ at time $\tau_Z$, we obtain

$$\begin{aligned} v_y(x) &= \mathbb{E}_x \exp\left\{\int_0^{\tau_y}\right\} \mathbb{1}\{\tau_y = \tau_\Gamma < \tau_{B^c}\} \\ &\le q^{|x-y|} + \sum_{z_0 \in Z\setminus\Gamma} q^{|x-z_0|} \mathbb{E}_{z_0} \exp\left\{\int_0^{\tau_y}\right\} \mathbb{1}\{\tau_y = \tau_\Gamma < \tau_{B^c}\}. \end{aligned} \tag{6.10}$$

Now we distinguish the three cases: (1) $z_0 \in \mathcal{C}_{\delta,\mathcal{R}} \setminus \Gamma$, (2) $z_0 \in Z \setminus Z_{\delta,\mathcal{R}}$ and (3) $z_0 \in \dot{Z}_{\delta,\mathcal{R}}$. In case (2) we use the strong Markov property at time $\tau_{\mathcal{C}_{\delta,\mathcal{R}}}$ and sum on all sites $z_1 \in [\mathcal{C}_{\delta,\mathcal{R}} \setminus \Gamma] \cup \{y\}$ the walker can occupy at that time. In case (3) we use the strong Markov property at the entrance time to the set $[Z \setminus Z_{\delta,\mathcal{R}}] \cup [\mathcal{C}_{\delta,\mathcal{R}} \setminus \Gamma] \cup \{y\}$ and distinguish if the walker is in $Z \setminus Z_{\delta,\mathcal{R}}$ or in $[\mathcal{C}_{\delta,\mathcal{R}} \setminus \Gamma] \cup \{y\}$ at that time. In the first case, we use afterward the strong Markov property at time $\tau_{\mathcal{C}_{\delta,\mathcal{R}}}$. This yields

$$
\begin{aligned}
v_y(x) \le q^{|x-y|} &+ \sum_{z_0 \in \mathcal{C}_{\delta,\mathcal{R}} \setminus \Gamma} q^{|x-z_0|} \mathbb{E}_{z_0} \exp\Big\{ \int_0^{\tau_y} \Big\} \mathbb{1}\{\tau_y = \tau_\Gamma < \tau_{B^c}\} \\
&+ \sum_{z_0 \in Z \setminus Z_{\delta,\mathcal{R}}} q^{|x-z_0|} \sum_{z_1 \in [\mathcal{C}_{\delta,\mathcal{R}} \setminus \Gamma] \cup \{y\}} \mathbb{E}_{z_0} \exp\Big\{ \int_0^{\tau_{z_1}} \Big\} \\
&\qquad \times \mathbb{1}\{\tau_{z_1} = \tau_{\mathcal{C}_{\delta,\mathcal{R}}} < \tau_{\Gamma \setminus \{y\}} \wedge \tau_{B^c}\} \\
&\qquad \times \mathbb{E}_{z_1} \exp\Big\{ \int_0^{\tau_y} \Big\} \mathbb{1}\{\tau_y = \tau_\Gamma < \tau_{B^c}\} \qquad (6.11) \\
&+ \sum_{z_0 \in \dot{Z}_{\delta,\mathcal{R}}} q^{|x-z_0|} \Bigg[ \sum_{z_1 \in Z \setminus Z_{\delta,\mathcal{R}}} \mathbb{E}_{z_0} \exp\Big\{ \int_0^{\tau_{z_1}} \Big\} \\
&\qquad \times \mathbb{1}\{\tau_{z_1} = \tau_{Z \setminus Z_{\delta,\mathcal{R}}} < \tau_{\mathcal{C}_{\delta,\mathcal{R}}} \wedge \tau_{B^c}\} \\
&\qquad \times \sum_{z_2 \in [\mathcal{C}_{\delta,\mathcal{R}} \setminus \Gamma] \cup \{y\}} \mathbb{E}_{z_1} \exp\Big\{ \int_0^{\tau_{z_2}} \Big\} \\
&\qquad\qquad \times \mathbb{1}\{\tau_{z_2} = \tau_{\mathcal{C}_{\delta,\mathcal{R}}} < \tau_{B^c}\} \\
&\qquad \times \mathbb{E}_{z_2} \exp\Big\{ \int_0^{\tau_y} \Big\} \mathbb{1}\{\tau_y = \tau_\Gamma < \tau_{B^c}\} \\
&\qquad + \sum_{z_1 \in [\mathcal{C}_{\delta,\mathcal{R}} \setminus \Gamma] \cup \{y\}} \mathbb{E}_{z_0} \exp\Big\{ \int_0^{\tau_{z_1}} \Big\} \\
&\qquad\qquad \times \mathbb{1}\{\tau_{\mathcal{C}_{\delta,\mathcal{R}}} < \tau_{Z \setminus Z_{\delta,\mathcal{R}}} \wedge \tau_{B^c}\} \\
&\qquad\qquad \times \mathbb{E}_{z_1} \exp\Big\{ \int_0^{\tau_y} \Big\} \mathbb{1}\{\tau_y = \tau_\Gamma < \tau_{B^c}\} \Bigg].
\end{aligned}
$$

We note that the term in the fourth, respectively ninth, respectively last line is equal to one if $z_1$, respectively $z_2$, is equal to $y$. Now we use Lemmas 6.2–6.4 for the respective expectations on the right and obtain, for some $c > 0$, depending on the quantities in (5.14) only,

$$
\begin{aligned}
v_y(x) \le q^{|x-y|} &+ \sum_{z_0 \in \mathcal{C}_{\delta,\mathcal{R}} \setminus \Gamma} q^{c|x-z_0|} q^{c|z_0-y|} \\
&+ \sum_{z_0 \in Z \setminus Z_{\delta,\mathcal{R}}} \sum_{z_1 \in [\mathcal{C}_{\delta,\mathcal{R}} \setminus \Gamma] \cup y} q^{c|x-z_0|} q^{c|z_0-z_1|} q^{c|z_1-y|}
\end{aligned}
$$

$$
\begin{aligned}
&+\sum_{z_0\in\dot Z_{\delta,\mathcal R}}\sum_{z_1\in Z\setminus Z_{\delta,\mathcal R}}\sum_{z_2\in[\mathcal C_{\delta,\mathcal R}\setminus\Gamma]\cup\{y\}} q^{c|x-z_0|}q^{c|z_0-z_1|}\\
&\qquad\qquad\times q^{c|z_1-z_2|}q^{c|z_2-y|}\\
&+M\sum_{z_0\in\dot Z_{\delta,\mathcal R}}\sum_{z_1\in[\mathcal C_{\delta,\mathcal R}\setminus\Gamma]\cup\{y\}} q^{c|x-z_0|}q^{c|z_0-z_1|}q^{c|z_1-y|}.
\end{aligned}\tag{6.12}
$$

We use the triangle inequality to obtain the estimate

$$
\begin{aligned}
v_y(x)\le q^{c/4|x-y|}\Bigg[&1+\sum_{z\in\mathbb Z^d}q^{c/4|z|}+\Bigg(\sum_{z\in\mathbb Z^d}q^{c/4|z|}\Bigg)^2\\
&+\Bigg(\sum_{z\in\mathbb Z^d}q^{c/4|z|}\Bigg)^3+M\Bigg(\sum_{z\in\mathbb Z^d}q^{c/4|z|}\Bigg)^2\Bigg].
\end{aligned}\tag{6.13}
$$

The quantity in the square brackets depends on the quantities in (5.14) and $M$ only. Hence, by choosing $R$ large enough (recall that $|x-y|\ge R$) and altering the value of $c$, we arrive at (6.1). This ends the proof of Proposition 6.1.

6.3. *Big-cluster expansion*: *Proof of Lemma* 6.2. In this section we carry out the details of the big-cluster expansion, that is, we prove Lemma 6.2.

We introduce a large auxiliary parameter $R\in\mathbb N$ which we choose in accordance with Lemma 5.4 [see (6.22) below] and additionally so large that (6.28) is satisfied and such that $q^{R/(24K)}\le 1-\sum_{|y|\ge R}q^{|y|/(12K)}$ holds [see the last step in (6.30)]. Assertion (b) will be shown with $M=CR^dq^{-R/4}$, where $C>0$ depends on the quantities in (5.14) only. The parameters $\delta>0$ and $\mathcal R>0$ are chosen in accordance with Lemma 5.4.

First we prove (a). Fix $z_0\in\dot Z_{\delta,\mathcal R}$ and $z\in Z\setminus Z_{\delta,\mathcal R}$. We shall decompose the path $X$ into the pieces between subsequent visits to dotted islands in $\dot Z_{\delta,\mathcal R}$ different from the starting island. To this end, we introduce the corresponding stopping time,

$$
\zeta_{\mathrm B}=\inf\{t>0\colon X_t\in Z_{\delta,\mathcal R}\setminus Z_R[X_0]\}.\tag{6.14}
$$

We repeatedly apply the strong Markov property at the time $\zeta_{\mathrm B}$ and sum over the walker's positions at these times to obtain

$$
\begin{aligned}
&\text{l.h.s. of (6.6)}\\
&\quad=\sum_{i=0}^{\infty}\sum_{x_1,\dots,x_i\in\dot Z_{\delta,\mathcal R}}\Bigg[\prod_{k=1}^{i}\mathbb E_{x_{k-1}}\exp\bigg\{\int_0^{\tau_{x_k}}\bigg\}
\end{aligned}\tag{6.15}
$$

$$\times \mathbb{1}\{\tau_{x_k} = \zeta_{\rm B} < \tau_{Z\setminus Z_{\delta,\mathcal{R}}} \wedge \tau_{\mathcal{C}_{\delta,\mathcal{R}}} \wedge \tau_{B^c}\}\Bigg]$$

$$\times \mathbb{E}_{x_i} \exp\bigg\{\int_0^{\tau_z}\bigg\}\mathbb{1}\{\tau_z = \tau_{Z\setminus Z_{\delta,\mathcal{R}}} < \zeta_{\rm B} \wedge \tau_{\mathcal{C}_{\delta,\mathcal{R}}} \wedge \tau_{B^c}\},$$

where $x_0 = z_0$. In the sum on $x_1, \dots, x_i$, we may and shall add the constraint that $x_{k-1}$ and $x_k$ lie in different islands for all $k = 1, \dots, i$, which implies that $|x_{k-1} - x_k| \geq 2R$. The summand for $i = 0$ is interpreted as just the last line for $i = 0$.

In the following we shall bound the expectations on the right-hand side of (6.15) as follows:

$$\mathbb{E}_{x_{k-1}} \exp\bigg\{\int_0^{\tau_{x_k}}\bigg\}\mathbb{1}\{\tau_{x_k} = \zeta_{\rm B} < \tau_{Z\setminus Z_{\delta,\mathcal{R}}} \wedge \tau_{\mathcal{C}_{\delta,\mathcal{R}}} \wedge \tau_{B^c}\} \leq q^{|x_{k-1} - x_k|/(4K)}, \quad 1 \leq k \leq i, \tag{6.16}$$

$$\mathbb{E}_{x_i} \exp\bigg\{\int_0^{\tau_z}\bigg\}\mathbb{1}\{\tau_z = \tau_{Z\setminus Z_{\delta,\mathcal{R}}} < \zeta_{\rm B} \wedge \tau_{\mathcal{C}_{\delta,\mathcal{R}}} \wedge \tau_{B^c}\} \leq q^{|x_i - z|/(4K)}, \quad i \geq 0. \tag{6.17}$$

In order to show (6.16)–(6.17), we need the stopping times

$$\eta_R = \inf\{t > 0 \colon X_t \notin B_R(Z_R[X_0])\}, \tag{6.18}$$

$$\sigma_R = \inf\{t > \eta_R \colon X_t \in Z_R[X_0]\}, \tag{6.19}$$

of the first exit time from the big cluster around the starting island and the next return to the same island. Let us consider the expectation on the left-hand side of (6.16) and expand according to the number of times before time $\zeta_{\rm B} \wedge \tau_{Z\setminus Z_{\delta,\mathcal{R}}} \wedge \tau_{\mathcal{C}_{\delta,\mathcal{R}}} \wedge \tau_{B^c}$ the walker leaves the big cluster around the starting island and returns to the same island. Hence, for $i \geq 1$ and $k = 1, \dots, i$, we have

$$\begin{aligned} &\text{l.h.s. of (6.16)} \\ &\quad = \sum_{m=0}^{\infty} \sum_{y_1,\dots,y_m \in \dot{Z}_R[y_0]} \Bigg[\prod_{l=1}^m \mathbb{E}_{y_{l-1}} \exp\bigg\{\int_0^{\tau_{y_l}}\bigg\} \\ &\qquad\qquad \times \mathbb{1}\{\tau_{y_l} = \sigma_R < \zeta_{\rm B} \wedge \tau_{Z\setminus Z_{\delta,\mathcal{R}}} \wedge \tau_{\mathcal{C}_{\delta,\mathcal{R}}} \wedge \tau_{B^c}\}\Bigg] \\ &\qquad \times \mathbb{E}_{y_m} \exp\bigg\{\int_0^{\tau_{x_k}}\bigg\} \\ &\qquad \times \mathbb{1}\{\tau_{x_k} = \zeta_{\rm B} < \sigma_R \wedge \tau_{Z\setminus Z_{\delta,\mathcal{R}}} \wedge \tau_{\mathcal{C}_{\delta,\mathcal{R}}} \wedge \tau_{B^c}\}, \end{aligned} \tag{6.20}$$

where $y_0 = x_{k-1}$. The summand for $m = 0$ is interpreted as just the last two lines for $m = 0$.

In each of the expectations on the right-hand side of (6.20), we use the strong Markov property at time $\eta_R$ to obtain

$$
\begin{aligned}
&\text{l.h.s. of (6.16)}\\
&\quad\le \sum_{m=0}^{\infty} \sum_{y_1,\dots,y_m\in \dot Z_R[y_0]} \Biggl[\prod_{l=1}^{m}\Biggl[\mathbb{E}_{y_{l-1}} \exp\Biggl\{\int_0^{\eta_R}\Biggr\} \mathbb{1}\{\eta_R < \tau_{\mathcal{C}_{\delta,\mathcal{R}}} \wedge \tau_{B^c}\}\\
&\qquad\qquad \times \max_{y'\in\partial B_R(Z_R[y_0])} \mathbb{E}_{y'} \exp\Biggl\{\int_0^{\tau_{y_l}}\Biggr\}\\
&\qquad\qquad \times \mathbb{1}\{\tau_{y_l} = \tau_Z < \tau_{B^c}\}\Biggr]\Biggr]\\
&\qquad \times \mathbb{E}_{y_m} \exp\Biggl\{\int_0^{\eta_R}\Biggr\} \mathbb{1}\{\eta_R < \tau_{\mathcal{C}_{\delta,\mathcal{R}}} \wedge \tau_{B^c}\}\\
&\qquad \times \max_{y'\in\partial B_R(Z_R[y_0])} \mathbb{E}_{y'} \exp\Biggl\{\int_0^{\tau_{x_k}}\Biggr\} \mathbb{1}\{\tau_{x_k} = \tau_Z < \tau_{B^c}\}.
\end{aligned}
\tag{6.21}
$$

(For a subset $A$ of $\mathbb{Z}^d$, we write $\partial A$ for the external boundary of $A\cap B$ in $B$.) On the event $\{\eta_R < \tau_{\mathcal{C}_{\delta,\mathcal{R}}} \wedge \tau_{B^c}\}$, the stopping time $\eta_R$ coincides with $\tau_{A^c}$ for $A = B_R(Z_R[y_m]) \setminus \mathcal{C}_{\delta,\mathcal{R}}$. Hence, an application of Lemma 4.2 with $\gamma = \lambda_y$ and an application of the spectral gap estimate in (5.16) of Lemma 5.4 yields that

$$
\begin{aligned}
&\mathbb{E}_{y_{l-1}} \exp\Biggl\{\int_0^{\eta_R}\Biggr\} \mathbb{1}\{\eta_R < \tau_{\mathcal{C}_{\delta,\mathcal{R}}} \wedge \tau_{B^c}\}\\
&\quad\le 1 + 2d \frac{|B_R(Z_R[y_0])|}{\lambda_y - \lambda(B_R(Z_R[y_0]) \setminus \mathcal{C}_{\delta,\mathcal{R}})}\\
&\quad\le 1 + 2d \frac{K|B_R|}{b/2 \vee 1/4(\lambda(V_\varrho) - \lambda(\dot V_\varrho))},
\end{aligned}
\tag{6.22}
$$

for any $l = 1,\dots,m+1$. Using this for the first and third expectation on the right-hand side of (6.21), and using (6.4) for the second and the fourth expectation, we obtain

$$
\begin{aligned}
&\text{r.h.s. of (6.21)}\\
&\quad\le \sum_{m=0}^{\infty} \sum_{y_1,\dots,y_m\in \dot Z_R[y_0]} \Biggl[\prod_{l=1}^{m}\Biggl[\Biggl(1 + 2d \frac{K|B_R|}{b/2 \vee 1/4(\lambda(V_\varrho) - \lambda(\dot V_\varrho))}\Biggr)\\
&\qquad\qquad \times q^{\operatorname{dist}(\partial B_R(Z_R[y_{l-1}]), y_l)}\Biggr]\Biggr]
\end{aligned}
\tag{6.23}
$$

$$\times \left(1 + 2d\frac{K|B_R|}{b/2 \vee 1/4(\lambda(V_\varrho) - \lambda(\dot{V}_\varrho))}\right)$$
$$\times\, q^{\mathrm{dist}(\partial B_R(Z_R[y_m]),x_k)}.$$

Clearly, since $y_l$ and $y_{l-1}$ belong to the same island,

$$\mathrm{dist}(\partial B_R(Z_R[y_{l-1}]), y_l) \geq R \geq \frac{R}{2} + \frac{|y_{l-1} - y_l|}{4K}, \qquad l = 1, \dots, m. \tag{6.24}$$

Furthermore,

$$\mathrm{dist}(\partial B_R(Z_R[y_m]), x_k) \geq \frac{|y_m - x_k|}{2K}. \tag{6.25}$$

Indeed, since $x_k$ and $y_m$ lie in different islands, and since $R \leq \mathrm{dist}(\partial B_R(Z_R[y_0]), x_k)$, we have

$$|y_m - x_k| \leq KR + \mathrm{dist}(\partial B_R(Z_R[y_0]), x_k) \leq 2K\, \mathrm{dist}(\partial B_R(Z_R[y_0]), x_k). \tag{6.26}$$

Combining (6.25) with $\mathrm{dist}(\partial B_R(Z_R[y_m]), x_k) \geq R$, we get

$$\mathrm{dist}(\partial B_R(Z_R[y_m]), x_k) \geq \frac{R}{2} + \frac{|y_m - x_k|}{4K}. \tag{6.27}$$

We substitute (6.24) and (6.27) in (6.23). Now assume that $R$ is so large that

$$q^{R/2}\left(1 + 2d\frac{K|B_R|}{b/2 \vee 1/4(\lambda(V_\varrho) - \lambda(\dot{V}_\varrho))}\right) < \frac{1}{2K}, \tag{6.28}$$

to obtain the bound

$$\begin{aligned} \text{l.h.s. of (6.16)} &\leq \sum_{m=0}^{\infty} \sum_{y_1,\dots,y_m \in Z_R[y_0]} \left[\prod_{l=1}^{m} \left(\frac{1}{2K} q^{|y_{l-1}-y_l|/(4K)}\right)\right] \\ &\qquad \times \frac{1}{2K} q^{|y_m - x_k|/(4K)} \\ &\leq \frac{1}{2} q^{|x_{k-1}-x_k|/(4K)} \sum_{m=0}^{\infty} \sum_{y_1,\dots,y_m \in Z_R[y_0]} (2K)^{-m} \\ &\leq q^{|x_{k-1}-x_k|/(4K)}, \end{aligned} \tag{6.29}$$

using the triangle inequality for the sequence of points $x_{k-1} = y_0, y_1, \dots, y_m, x_k$, and noting that $|Z_R[y_0]| \leq K$. Hence, we have proved (6.16). The proof of (6.17) is analogous, replacing $x_{k-1}$ by $x_i$ and $x_k$ by $z$.

Substituting (6.16) and (6.17) into (6.15), we obtain

$$
\begin{aligned}
\text{l.h.s. of (6.6)} &\le \sum_{i=0}^{\infty} \sum_{x_1,\dots,x_i\in\dot Z_{\delta,\mathcal R}\colon |x_{k-1}-x_k|\ge R\,\forall k} \Biggl[\prod_{k=1}^{i} q^{|x_{k-1}-x_k|/(4K)}\Biggr] \\
&\qquad\qquad \times q^{|x_i-z|/(4K)} \\
&\le q^{|z-z_0|/(6K)} \\
&\qquad \times \sum_{i=0}^{\infty} \sum_{x_1,\dots,x_i\in\mathbb Z^d\colon |x_{k-1}-x_k|\ge R\,\forall k} \prod_{k=1}^{i} q^{|x_{k-1}-x_k|/(12K)} \qquad (6.30)\\
&\le q^{|z-z_0|/(6K)} \sum_{i=0}^{\infty} \Biggl(\sum_{y\in\mathbb Z^d\colon |y|\ge R} q^{|y|/(12K)}\Biggr)^{i} \\
&\le q^{|z-z_0|/(8K)},
\end{aligned}
$$

where in the last step we used that $|z-z_0|\ge R$ and assumed that $R$ is sufficiently large, depending on the quantities in (5.14) only. This ends the proof of (6.6), that is, of assertion (a) of Lemma 6.2.

Now we turn to the proof of (6.7), that is, of assertion (b) of Lemma 6.2. This is analogous to the proof of (6.6). Indeed, interchange $\tau_{\mathcal C_{\delta,\mathcal R}}$ and $\tau_{Z\setminus Z_{\delta,\mathcal R}}$ in (6.15) with each other to obtain an expansion analogous to (6.15):

$$
\begin{aligned}
&\text{l.h.s. of (6.7)} \\
&\quad = \sum_{i=0}^{\infty} \sum_{x_1,\dots,x_i\in\dot Z_{\delta,\mathcal R}} \Biggl[\prod_{k=1}^{i} \mathbb E_{x_{k-1}} \exp\biggl\{\int_0^{\tau_{x_k}}\biggr\} \\
&\qquad\qquad \times \mathbb 1\{\tau_{x_k}=\zeta_{\mathrm B}<\tau_{Z\setminus Z_{\delta,\mathcal R}}\wedge\tau_{\mathcal C_{\delta,\mathcal R}}\wedge\tau_{B^{\mathrm c}}\}\Biggr] \qquad (6.31)\\
&\qquad \times \mathbb E_{x_i} \exp\biggl\{\int_0^{\tau_z}\biggr\} \\
&\qquad \times \mathbb 1\{\tau_z=\tau_{\mathcal C_{\delta,\mathcal R}}<\zeta_{\mathrm B}\wedge\tau_{Z\setminus Z_{\delta,\mathcal R}}\wedge\tau_{B^{\mathrm c}}\}.
\end{aligned}
$$

Note that the first two lines of the right-hand side of (6.31) are the same as in (6.15), hence, (6.16) can be used for this term as well. The last expectation is equal to zero if $x_i$ does not lie in $Z_R[z]$. Otherwise, we estimate

$$
\begin{aligned}
&\mathbb E_{x_i} \exp\biggl\{\int_0^{\tau_z}\biggr\}\mathbb 1\{\tau_z=\tau_{\mathcal C_{\delta,\mathcal R}}<\zeta_{\mathrm B}\wedge\tau_{Z\setminus Z_{\delta,\mathcal R}}\wedge\tau_{B^{\mathrm c}}\} \\
&\qquad \le \frac{C}{2} R^d q^{-R/4} q^{|x_i-z|/(4K)} \qquad \text{if } x_i\in Z_R[z].
\end{aligned} \qquad (6.32)
$$

In the same way as (6.6) follows from (6.15) in combination with (6.16)–(6.17), (6.7) follows from (6.31) in combination with (6.32); we omit the details.

The proof of (6.32) is similar to the proof of (6.17), but requires some changes, as we explain now. Fix $z \in \mathcal{C}_{\delta,\mathcal{R}} \cap Z_R[x_i]$. We again use an expansion as in (6.20), with $x_{k-1}$ replaced by $x_i$ and $x_k$ replaced by $z$. We also may apply the strong Markov property at the stopping time $\eta_R$ to the expectation in the first line of (6.20) and Lemma 4.2 and (6.4), in the same way as in (6.21). This yields the bound

$$
\begin{aligned}
&\text{l.h.s. of } (6.32)\\
&\quad\le \sum_{m=0}^{\infty} \sum_{y_1,\dots,y_m \in \dot Z_R[y_0]} \left[\prod_{l=1}^{m}\left[\left(1+2d\frac{K|B_R|}{b/2 \vee 1/4(\lambda(V_\varrho)-\lambda(\dot V_\varrho))}\right)\right.\right.\\
&\qquad\qquad\left.\left.\times\, q^{\operatorname{dist}(\partial B_R(Z_R[y_{l-1}]),y_l)}\right]\right]\\
&\qquad\qquad\times\, \mathbb{E}_{y_m} \exp\left\{\int_0^{\tau_z}\right\} \mathbb{1}\{\tau_z < \tau_{B_R(Z_R[y_0])^{\mathrm{c}}}\}\\
&\quad\le \sum_{m=0}^{\infty} \sum_{y_1,\dots,y_m \in \dot Z_R[y_0]} \left[\prod_{l=1}^{m}\left(\frac{1}{2K} q^{|y_{l-1}-y_l|/(4K)}\right)\right]\\
&\qquad\qquad\times\, \mathbb{E}_{y_m} \exp\left\{\int_0^{\tau_z}\right\} \mathbb{1}\{\tau_z < \tau_{B_R(Z_R[y_0])^{\mathrm{c}}}\},
\end{aligned}
\tag{6.33}
$$

where $y_0 = x_i$. Now we estimate the last expectation differently: we directly apply Lemma 4.2 with $\gamma = \lambda_y$ and $A = B_R(Z_R[y_0]) \setminus \{z\}$ and use the spectral gap in (5.16) to estimate

$$
\begin{aligned}
&\mathbb{E}_{y_m} \exp\left\{\int_0^{\tau_z}\right\} \mathbb{1}\{\tau_z < \tau_{B_R(Z_R[y_0])^{\mathrm{c}}}\}\\
&\quad\le 1 + 2d\frac{|B_R(Z_R[y_0])|}{\lambda_y - \lambda(B_R(Z_R[y_{l-1}]) \setminus \mathcal{C}_{\delta,\mathcal{R}})}\\
&\quad\le 1 + 2d\frac{K|B_R|}{b/2 \vee 1/4(\lambda(V_\varrho)-\lambda(\dot V_\varrho))} \le \frac{C}{4} R^d,
\end{aligned}
\tag{6.34}
$$

for some suitable choice of $C$ which depends on the quantities in (5.14) only. Since $y_m$ and $z$ belong to the same island, we can estimate $|y_m - z| \le RK$, and hence, we may continue the right-hand side of (6.34) with $\le \frac{C}{4} R^d q^{-R/4} q^{|y_m - z|/(4K)}$. Using this in (6.33), and arguing as above, we arrive at (6.32).

6.4. *Large-cluster expansion*: *Proof of Lemma* 6.3. In this section we carry out the details of the large-cluster expansion, that is, we prove Lemma 6.3. We pick $M > 1$ as in Lemma 6.2(b) and assume that $\delta$ and $\mathcal{R}$ are chosen in accordance with Lemma 6.2. Additionally, we require that (6.52) below holds.

Fix $z_0 \in Z \setminus Z_{\delta,\mathcal{R}}$ and $z \in \mathcal{C}_{\delta,\mathcal{R}}$. We shall divide the trajectory $X$ into the pieces between subsequent visits to different non-$(\delta,\mathcal{R})$-optimal archipelagos before time $\tau_{\mathcal{C}_{\delta,\mathcal{R}}} \wedge \tau_{B^c}$. Introduce the corresponding stopping time,

$$\zeta_{\mathrm{L}} = \inf\{t > 0 : X_t \in (Z \setminus Z_{\delta,\mathcal{R}}) \setminus Z_{\mathcal{R}}[X_0]\}. \tag{6.35}$$

We repeatedly use the strong Markov property at the time $\zeta_{\mathrm{L}}$ and sum on all the walker's positions at these times to obtain

$$\begin{aligned}
&\text{l.h.s. of (6.8)}\\
&\quad = \sum_{i=0}^{\infty} \sum_{x_1,\dots,x_i \in Z\setminus Z_{\delta,\mathcal{R}}} \Biggl[\prod_{k=1}^{i} \mathbb{E}_{x_{k-1}} \exp\biggl\{\int_0^{\tau_{x_k}}\biggr\} \\
&\qquad\qquad \times \mathbb{1}\{\tau_{x_k} = \zeta_{\mathrm{L}} < \tau_{\mathcal{C}_{\delta,\mathcal{R}}} \wedge \tau_{B^c}\}\Biggr] \\
&\qquad\qquad \times \mathbb{E}_{x_i} \exp\biggl\{\int_0^{\tau_z}\biggr\} \mathbb{1}\{\tau_z = \tau_{\mathcal{C}_{\delta,\mathcal{R}}} < \zeta_{\mathrm{L}} \wedge \tau_{B^c}\},
\end{aligned} \tag{6.36}$$

where we put $x_0 = z_0$. In the sum on $x_1, \dots, x_i$, we may and shall add the constraint that $x_{k-1}$ and $x_k$ are in different archipelagos for any $k = 1, \dots, i$ (which implies that $|x_{k-1} - x_k| \geq 2\mathcal{R}$). The summand for $i = 0$ is defined to be just the term in the last line for $i = 0$.

We are going to further estimate the expectations on the right-hand side of (6.36) as follows:

$$\begin{aligned}
&\mathbb{E}_{x_{k-1}} \exp\biggl\{\int_0^{\tau_{x_k}}\biggr\} \mathbb{1}\{\tau_{x_k} = \zeta_{\mathrm{L}} < \tau_{\mathcal{C}_{\delta,\mathcal{R}}} \wedge \tau_{B^c}\} \leq q^{|x_{k-1}-x_k|/(64K^2)}, \\
&\qquad\qquad k = 1, \dots, i,
\end{aligned} \tag{6.37}$$

$$\mathbb{E}_{x_i} \exp\biggl\{\int_0^{\tau_z}\biggr\} \mathbb{1}\{\tau_z = \tau_{\mathcal{C}_{\delta,\mathcal{R}}} < \zeta_{\mathrm{L}} \wedge \tau_{B^c}\} \leq q^{|x_i - z|/(64K^2)}, \qquad i \geq 0. \tag{6.38}$$

In order to show (6.37), we expand according to the number of times the walker leaves the large cluster $B_{\mathcal{R}}(Z_{\mathcal{R}}[X_0])$ and revisits the starting archipelago $Z_{\mathcal{R}}[X_0]$. Define

$$\eta_{\mathcal{R}} = \inf\{t > 0 : X_t \notin B_{\mathcal{R}}(Z_{\mathcal{R}}[X_0])\}, \tag{6.39}$$

$$\sigma_{\mathcal{R}} = \inf\{t > \eta_{\mathcal{R}} : X_t \in Z_{\mathcal{R}}[X_0]\}. \tag{6.40}$$

We repeatedly use the strong Markov property at time $\sigma_{\mathcal{R}}$ and sum over the sites visited by the walker at these times to obtain

$$\begin{aligned}
&\text{l.h.s. of (6.37)}\\
&\quad=\sum_{m=0}^{\infty}\sum_{y_1,\dots,y_m\in Z_{\mathcal{R}}[x_{k-1}]}\Bigg[\prod_{l=1}^{m}\mathbb{E}_{y_{l-1}}\exp\bigg\{\int_0^{\tau_{y_l}}\bigg\}\\
&\qquad\qquad\times\mathbb{1}\{\tau_{y_l}=\sigma_{\mathcal{R}}<\zeta_{\mathrm{L}}\wedge\tau_{\mathcal{C}_{\delta,\mathcal{R}}}\wedge\tau_{B^{\mathrm{c}}}\}\Bigg]\\
&\qquad\qquad\times\mathbb{E}_{y_m}\exp\bigg\{\int_0^{\tau_{x_k}}\bigg\}\\
&\qquad\qquad\times\mathbb{1}\{\tau_{x_k}=\zeta_{\mathrm{L}}<\sigma_{\mathcal{R}}\wedge\tau_{\mathcal{C}_{\delta,\mathcal{R}}}\wedge\tau_{B^{\mathrm{c}}}\},
\end{aligned}\tag{6.41}$$

where $y_0=x_{k-1}$. The summand for $m=0$ is interpreted as just the last line.

In the following we shall show that

$$\begin{aligned}
&\mathbb{E}_{y_{l-1}}\exp\bigg\{\int_0^{\tau_{y_l}}\bigg\}\mathbb{1}\{\tau_{y_l}=\sigma_{\mathcal{R}}<\zeta_{\mathrm{L}}\wedge\tau_{\mathcal{C}_{\delta,\mathcal{R}}}\wedge\tau_{B^{\mathrm{c}}}\}\le\frac{1}{2K}q^{|y_{l-1}-y_l|/(64K^2)},\\
&\qquad\qquad l=1,\dots,m,
\end{aligned}\tag{6.42}$$

$$\begin{aligned}
&\mathbb{E}_{y_m}\exp\bigg\{\int_0^{\tau_{x_k}}\bigg\}\mathbb{1}\{\tau_{x_k}=\zeta_{\mathrm{L}}<\sigma_{\mathcal{R}}\wedge\tau_{\mathcal{C}_{\delta,\mathcal{R}}}\wedge\tau_{B^{\mathrm{c}}}\}\le\tfrac{1}{2}q^{|y_m-x_k|/(64K^2)},\\
&\qquad\qquad m\ge 0.
\end{aligned}\tag{6.43}$$

Substituting (6.42) and (6.43) in (6.41), we obtain

$$\begin{aligned}
\text{l.h.s. of (6.37)}&\le\sum_{m=0}^{\infty}\sum_{y_1,\dots,y_m\in Z_{\mathcal{R}}[x_{k-1}]}\Bigg[\prod_{l=1}^{m}\frac{1}{2K}q^{|y_{l-1}-y_l|/(64K^2)}\Bigg]\\
&\qquad\times\frac{1}{2}q^{|y_m-x_k|/(64K^2)}\\
&\le q^{|x_{k-1}-x_k|/(64K^2)},
\end{aligned}\tag{6.44}$$

where we used the triangle inequality for the sequence of points $x_{k-1}=y_0,y_1,\dots,y_m,x_k$, and the fact that $|Z_{\mathcal{R}}[x_{k-1}]|\le K$. Hence, we have derived (6.37).

Now we prove (6.42) and (6.43). Apply the strong Markov property at time $\eta_{\mathcal{R}}$ to obtain, for any $l=1,\dots,m$,

$$\begin{aligned}
&\text{l.h.s. of (6.42)}\\
&\quad\le\mathbb{E}_{y_{l-1}}\exp\bigg\{\int_0^{\eta_{\mathcal{R}}}\bigg\}\mathbb{1}\{\eta_{\mathcal{R}}<\tau_{B^{\mathrm{c}}}\}\\
&\qquad\times\max_{y'\in\partial B_{\mathcal{R}}(Z_{\mathcal{R}}[y_0])}\mathbb{E}_{y'}\exp\bigg\{\int_0^{\tau_{y_l}}\bigg\}\mathbb{1}\{\tau_{y_l}=\tau_{Z\setminus Z_{\delta,\mathcal{R}}}<\tau_{\mathcal{C}_{\delta,\mathcal{R}}}\wedge\tau_{B^{\mathrm{c}}}\}.
\end{aligned}\tag{6.45}$$

In order to further estimate the first factor on the right-hand side of (6.45), use Lemma 4.2 for $\gamma = \lambda_y$ and the set $A = B_{\mathcal{R}}(Z_{\mathcal{R}}[y_0])$, and use the spectral gap in (5.17), to obtain

$$\begin{aligned}\mathbb{E}_{y_{l-1}} \exp\left\{\int_0^{\eta_{\mathcal{R}}}\right\}\mathbb{1}\{\eta_{\mathcal{R}} < \tau_{B^c}\} &\le 1 + 2d\frac{|B_{\mathcal{R}}(Z_{\mathcal{R}}[y_0])|}{\lambda_y - \lambda(B_{\mathcal{R}}(Z_{\mathcal{R}}[y_0]))} \\ &\le 1 + 2d\frac{K|B_{\mathcal{R}}|}{\delta/2}.\end{aligned} \tag{6.46}$$

In order to further estimate the second term on the right-hand side of (6.45), we use the Markov property at time $\tau_Z$ and (6.4). Furthermore, we distinguish whether the walker is in $\dot{Z}_{\delta,\mathcal{R}}$ at time $\tau_Z$ or already at $y_l$. (There are no other possibilities on the event $\{\tau_{y_l} = \tau_{Z\setminus Z_{\delta,\mathcal{R}}} < \tau_{\mathcal{C}_{\delta,\mathcal{R}}} \wedge \tau_{B^c}\}$.) In this way we obtain, for any $y' \in \partial B_{\mathcal{R}}(Z_{\mathcal{R}}[y_0])$,

$$\begin{aligned}&\mathbb{E}_{y'} \exp\left\{\int_0^{\tau_{y_l}}\right\}\mathbb{1}\{\tau_{y_l} = \tau_{Z\setminus Z_{\delta,\mathcal{R}}} < \tau_{\mathcal{C}_{\delta,\mathcal{R}}} \wedge \tau_{B^c}\} \\ &\quad\le \sum_{z'\in\dot{Z}_{\delta,\mathcal{R}}} q^{|y'-z'|}\mathbb{E}_{z'} \exp\left\{\int_0^{\tau_{y_l}}\right\}\mathbb{1}\{\tau_{y_l} = \tau_{Z\setminus Z_{\delta,\mathcal{R}}} < \tau_{\mathcal{C}_{\delta,\mathcal{R}}} \wedge \tau_{B^c}\} \\ &\qquad + q^{|y_l-y'|}.\end{aligned} \tag{6.47}$$

We can apply the big-cluster expansion, Lemma 6.2(a), to the expectation on the right-hand side and obtain

$$\begin{aligned}&\mathbb{E}_{y'} \exp\left\{\int_0^{\tau_{y_l}}\right\}\mathbb{1}\{\tau_{y_l} = \tau_{Z\setminus Z_{\delta,\mathcal{R}}} < \tau_{\mathcal{C}_{\delta,\mathcal{R}}} \wedge \tau_{B^c}\} \\ &\quad\le \sum_{z'\in\dot{Z}_{\delta,\mathcal{R}}} q^{|y'-z'|}q^{|z'-y_l|/(8K)} + q^{|y_l-y'|} \\ &\quad\le q^{|y'-y_l|/(16K)} \sum_{z'\in\mathbb{Z}^d} q^{|z'|/(16K)} + q^{|y_l-y'|} \\ &\quad\le \widetilde{C} q^{|y'-y_l|/(16K)},\end{aligned} \tag{6.48}$$

where $\widetilde{C} > 0$ depends on the quantities in (5.14) only. Use (6.46) and (6.48) in (6.45) to get, for any $l = 1, \dots, m$,

$$\text{l.h.s. of (6.42)} \le \widetilde{C}\left(1 + 2d\frac{K|B_{\mathcal{R}}|}{\delta/2}\right) \max_{y'\in\partial B_{\mathcal{R}}(Z_{\mathcal{R}}[y_0])} q^{|y'-y_l|/(16K)}. \tag{6.49}$$

In the same way, one derives

$$\text{l.h.s. of (6.43)} \le \widetilde{C}\left(1 + 2d\frac{K|B_{\mathcal{R}}|}{\delta/2}\right) \max_{y'\in\partial B_{\mathcal{R}}(Z_{\mathcal{R}}[y_0])} q^{|y'-x_k|/(16K)}. \tag{6.50}$$

Since $y_l$ and $y_0$ lie in the same archipelago, we estimate the last term in (6.49) against $q^{\mathcal{R}/(16K)}$. Furthermore, since $|y_{l-1} - y_l| \le \mathcal{R}K$, we may further estimate $q^{\mathcal{R}/(16K)} \le q^{\mathcal{R}/(32K)} q^{|y_{l-1}-y_l|/(32K^2)}$.

In (6.50) we estimate $q^{|y'-x_k|/(16K)} \le q^{\mathcal{R}/(32K)} q^{|y_m - x_k|/(64K^2)}$, where we used that $|y' - x_k| \ge \mathcal{R}$ and that $x_k$ and $y_m$ lie in different archipelagos, and we estimated

$$
\begin{aligned}
|y_m - x_k| &\le K\mathcal{R} + \operatorname{dist}(\partial B_{\mathcal{R}}(Z_{\mathcal{R}}[y_0]), x_k) \\
&\le 2K \operatorname{dist}(\partial B_{\mathcal{R}}(Z_{\mathcal{R}}[y_0]), x_k).
\end{aligned} \tag{6.51}
$$

Now we make the additional assumption that $\mathcal{R}$ is so large (depending only on $R$ and $\delta$) that

$$
\widetilde{C} M \Big(1 + 2d \frac{K|B_{\mathcal{R}}|}{\delta/2}\Big) q^{\mathcal{R}/(32K)} < \frac{1}{2K}, \tag{6.52}
$$

where $M > 1$ is as in Lemma 6.2(b). Then (6.42) and (6.43) easily follow from (6.49) and (6.50), respectively, in combination with the bounds given below (6.50). As remarked earlier, this finishes the proof of (6.37). The proof of (6.38) is analogous, the main difference being that in the expansion analogous to (6.47) (where $\tau_{\mathcal{C}_{\delta,\mathcal{R}}}$ and $\tau_{Z \setminus Z_{\delta,\mathcal{R}}}$ are interchanged), we apply Lemma 6.2(b) rather than (a). [The factor of $M$ in (6.52) is needed only here.]

Now we substitute (6.37)–(6.38) in (6.36) to obtain

$$
\begin{aligned}
&\mathbb{E}_{z_0} \exp\Big\{\int_0^{\tau_z}\Big\} \mathbb{1}\{\tau_z = \tau_{\mathcal{C}_{\delta,\mathcal{R}}} < \tau_{B^c}\} \\
&\quad \le \sum_{i=0}^{\infty} \sum_{x_1,\dots,x_i \in Z \setminus Z_{\delta,\mathcal{R}} \colon |x_{k-1}-x_k| \ge \mathcal{R}\, \forall k} \Big[\prod_{k=1}^{i} q^{|x_{k-1}-x_k|/(64K^2)}\Big] \\
&\qquad\qquad \times q^{|x_i - z|/(64K^2)} \\
&\quad \le q^{|z_0 - z|/(128K^2)} \sum_{i=0}^{\infty} \sum_{x_1,\dots,x_i \in \mathbb{Z}^d \colon |x_{k-1}-x_k| \ge \mathcal{R}\, \forall k} \prod_{k=1}^{i} q^{|x_{k-1}-x_k|/(128K^2)} \\
&\quad \le q^{|z_0 - z|/(128K^2)} \sum_{i=0}^{\infty} \Big(\sum_{y \in \mathbb{Z}^d \colon |y| \ge \mathcal{R}} q^{|y|/(128K^2)}\Big)^i \\
&\quad \le q^{|z_0 - z|/(256K^2)},
\end{aligned} \tag{6.53}
$$

if $\mathcal{R}$ is sufficiently large (recall that $|z_0 - z| \ge \mathcal{R}$). This implies the assertion in (6.8) and ends the proof of Lemma 6.3.

6.5. *Huge-cluster expansion*: *Proof of Lemma* 6.4. In this section we carry out the details of the huge-cluster expansion, that is, we prove Lemma 6.4. We pick $M > 1$ as in Lemma 6.2(b) and assume that $\delta$ and $\mathcal{R}$ are chosen in accordance with Lemmas 6.2 and 6.3.

We shall divide the path $X$ into the pieces between visits to different sites in $\mathcal{C}_{\delta,\mathcal{R}} \setminus \Gamma$. Introduce the corresponding stopping time,

$$\zeta_{\mathrm{H}} = \inf\{t > 0 : X_t \in (\mathcal{C}_{\delta,\mathcal{R}} \setminus \Gamma) \setminus \{X_0\}\}. \tag{6.54}$$

Fix $z_0 \in \mathcal{C}_{\delta,\mathcal{R}} \setminus \Gamma$ and $y \in \Gamma$. We repeatedly use the strong Markov property at time $\zeta_{\mathrm{H}}$ and sum over the sites visited at these times to obtain

$$\begin{aligned} &\text{l.h.s. of (6.9)} \\ &\quad = \sum_{i=0}^{\infty} \sum_{x_1,\dots,x_i \in \mathcal{C}_{\delta,\mathcal{R}} \setminus \Gamma} \Biggl[\prod_{k=1}^{i} \mathbb{E}_{x_{k-1}} \exp\biggl\{\int_0^{\tau_{x_k}}\biggr\} \mathbb{1}\{\tau_{x_k} = \zeta_{\mathrm{H}} < \tau_\Gamma \wedge \tau_{B^{\mathrm{c}}}\}\Biggr] \\ &\qquad\qquad \times \mathbb{E}_{x_i} \exp\biggl\{\int_0^{\tau_y}\biggr\} \mathbb{1}\{\tau_y = \tau_\Gamma < \zeta_{\mathrm{H}} \wedge \tau_{B^{\mathrm{c}}}\}, \end{aligned} \tag{6.55}$$

where $x_0 = z_0$. The summand for $i = 0$ is interpreted as just the last line with $i = 0$. In the sum on $x_1, \dots, x_i$ we may add the constraint $x_{k-1} \neq x_k$ for all $k = 1, \dots, i$. Recall that the huge clusters around the sites in $\mathcal{C}_{\delta,\mathcal{R}}$ are disjoint for $t$ large, since $\mathfrak{R} = \mathfrak{R}_t = \log^2 t \ll d_t^\delta$. Hence, we may and shall assume that $|x_{k-1} - x_k| \geq \mathfrak{R}$ for all $k = 1, \dots, i$.

We shall show the following estimates for the expectations on the right-hand side of (6.55):

$$\begin{aligned} &\mathbb{E}_{x_{k-1}} \exp\biggl\{\int_0^{\tau_{x_k}}\biggr\} \mathbb{1}\{\tau_{x_k} = \zeta_{\mathrm{H}} < \tau_\Gamma \wedge \tau_{B^{\mathrm{c}}}\} \leq q^{|x_{k-1} - x_k|/(2048K^2)}, \\ &\qquad\qquad k = 1, \dots, i, \end{aligned} \tag{6.56}$$

$$\mathbb{E}_{x_i} \exp\biggl\{\int_0^{\tau_y}\biggr\} \mathbb{1}\{\tau_y = \tau_\Gamma < \zeta_{\mathrm{H}} \wedge \tau_{B^{\mathrm{c}}}\} \leq q^{|x_i - y|/(2048K^2)}, \qquad i \geq 0. \tag{6.57}$$

We shall prove (6.56) only, the proof of (6.57) is identical.

In order to derive (6.56), we expand according to the number of times the walker leaves the huge cluster around its starting point and returns to the starting point before $\tau_{x_k}$. Denote by

$$\eta_{\mathfrak{R}} = \inf\{t > 0 : X_t \notin B_{\mathfrak{R}}(X_0)\}, \tag{6.58}$$

$$\sigma_{\mathfrak{R}} = \inf\{t > \eta_{\mathfrak{R}} : X_t = X_0\}, \tag{6.59}$$

the exit time from the huge cluster $B_{\mathfrak{R}}(X_0)$ and the next return time to $X_0$. Use the strong Markov property repeatedly at time $\sigma_{\mathfrak{R}}$ to get, for $k = 1, \dots, i$,

$$\mathbb{E}_{x_{k-1}} \exp\biggl\{\int_0^{\tau_{x_k}}\biggr\} \mathbb{1}\{\tau_{x_k} = \zeta_{\mathrm{H}} < \tau_\Gamma \wedge \tau_{B^{\mathrm{c}}}\}$$

$$(6.60)\quad \begin{aligned} &= \sum_{m=0}^{\infty}\Big[\mathbb{E}_{x_{k-1}}\exp\Big\{\int_0^{\sigma_{\mathfrak{R}}}\Big\}\mathbb{1}\{\sigma_{\mathfrak{R}}<\zeta_{\mathrm{H}}\wedge\tau_\Gamma\wedge\tau_{B^c}\}\Big]^m \\ &\qquad\times \mathbb{E}_{x_{k-1}}\exp\Big\{\int_0^{\tau_{x_k}}\Big\}\mathbb{1}\{\tau_{x_k}=\zeta_{\mathrm{H}}<\sigma_{\mathfrak{R}}\wedge\tau_\Gamma\wedge\tau_{B^c}\} \\ &= \frac{\mathbb{E}_{x_{k-1}}\exp\{\int_0^{\tau_{x_k}}\}\mathbb{1}\{\tau_{x_k}=\zeta_{\mathrm{H}}<\sigma_{\mathfrak{R}}\wedge\tau_\Gamma\wedge\tau_{B^c}\}}{1-\mathbb{E}_{x_{k-1}}\exp\{\int_0^{\sigma_{\mathfrak{R}}}\}\mathbb{1}\{\sigma_{\mathfrak{R}}<\zeta_{\mathrm{H}}\wedge\tau_\Gamma\wedge\tau_{B^c}\}}. \end{aligned}$$

The expectation in the denominator of the right-hand side of (6.60) is estimated by using the strong Markov property at time $\eta_{\mathfrak{R}}$ as follows:

$$(6.61)\quad \begin{aligned} &\mathbb{E}_{x_{k-1}}\exp\Big\{\int_0^{\sigma_{\mathfrak{R}}}\Big\}\mathbb{1}\{\sigma_{\mathfrak{R}}<\zeta_{\mathrm{H}}\wedge\tau_\Gamma\wedge\tau_{B^c}\} \\ &\quad\le \mathbb{E}_{x_{k-1}}\exp\Big\{\int_0^{\eta_{\mathfrak{R}}}\Big\}\mathbb{1}\{\eta_{\mathfrak{R}}<\tau_{B^c}\} \\ &\qquad\times\max_{z'\in\partial B_{\mathfrak{R}}(x_{k-1})}\mathbb{E}_{z'}\exp\Big\{\int_0^{\tau_{x_{k-1}}}\Big\}\mathbb{1}\{\tau_{x_{k-1}}=\tau_{\mathcal{C}_{\delta,\mathcal{R}}}<\tau_{B^c}\}. \end{aligned}$$

The numerator on the right-hand side of (6.60) is estimated in the same way:

$$(6.62)\quad \begin{aligned} &\mathbb{E}_{x_{k-1}}\exp\Big\{\int_0^{\zeta_{\mathrm{H}}}\Big\}\mathbb{1}\{\tau_{x_k}=\zeta_{\mathrm{H}}<\sigma_{\mathfrak{R}}\wedge\tau_\Gamma\wedge\tau_{B^c}\} \\ &\quad\le \mathbb{E}_{x_{k-1}}\exp\Big\{\int_0^{\eta_{\mathfrak{R}}}\Big\}\mathbb{1}\{\eta_{\mathfrak{R}}<\tau_{B^c}\} \\ &\qquad\times\max_{z'\in\partial B_{\mathfrak{R}}(x_{k-1})}\mathbb{E}_{z'}\exp\Big\{\int_0^{\tau_{x_k}}\Big\}\mathbb{1}\{\tau_{x_k}=\tau_{\mathcal{C}_{\delta,\mathcal{R}}}<\tau_{B^c}\}. \end{aligned}$$

Using Lemma 4.2 for $\gamma=\lambda_y$ and $A=B_{\mathfrak{R}}(x_{k-1})$, and using the spectral gap estimate in (5.18), we further estimate the first term on the right-hand side from above by

$$(6.63)\quad \begin{aligned} \mathbb{E}_{x_{k-1}}\exp\Big\{\int_0^{\eta_{\mathfrak{R}}}\Big\}\mathbb{1}\{\eta_{\mathfrak{R}}<\tau_{B^c}\} &\le 1+2d\frac{|B_{\mathfrak{R}}|}{\lambda_y-\lambda(B_{\mathfrak{R}}(x_{k-1}))} \\ &\le 1+2d\frac{|B_{\mathfrak{R}}|}{\mathfrak{g}^0}. \end{aligned}$$

In order to handle the expectations on the right-hand sides of (6.61) and (6.62) simultaneously, we shall prove the following.

LEMMA 6.5. *For any $x\in B$ and any $z\in\mathcal{C}_{\delta,\mathcal{R}}$,*

$$(6.64)\quad \mathbb{E}_x\exp\Big\{\int_0^{\tau_z}\Big\}\mathbb{1}\{\tau_z=\tau_{\mathcal{C}_{\delta,\mathcal{R}}}<\tau_{B^c}\}\le C' q^{|x-z|/(512K^2)},$$

*where the constant $C'>0$ depends on the quantities in (5.14) only.*

PROOF. This proof is a simple variant of the completion of the proof of Proposition 6.1 at the end of Section 6.2.

We apply the strong Markov property at time $\tau_Z$ and the estimate in (6.4). Furthermore, we distinguish the two cases that, at time $\tau_Z$, the walker is already at $z$ and that he is in the set $Z \setminus \mathcal{C}_{\delta,\mathcal{R}}$. This gives

$$\begin{aligned}
&\mathbb{E}_x \exp\biggl\{\int_0^{\tau_z}\biggr\}\mathbb{1}\{\tau_z = \tau_{\mathcal{C}_{\delta,\mathcal{R}}} < \tau_{B^c}\}\\
&\quad\le q^{|x-z|} + \sum_{z_0\in Z\setminus\mathcal{C}_{\delta,\mathcal{R}}} q^{|x-z_0|}\mathbb{E}_{z_0}\exp\biggl\{\int_0^{\tau_z}\biggr\}\mathbb{1}\{\tau_z = \tau_{\mathcal{C}_{\delta,\mathcal{R}}} < \tau_{B^c}\}.
\end{aligned} \tag{6.65}$$

In the sum over $z_0$ we distinguish the two cases (1) $z_0 \in Z \setminus Z_{\delta,\mathcal{R}}$ and (2) $z_0 \in \dot{Z}_{\delta,\mathcal{R}}$. In the second case we distinguish the two cases $\tau_{Z\setminus Z_{\delta,\mathcal{R}}} < \tau_{\mathcal{C}_{\delta,\mathcal{R}}}$ and $\tau_{\mathcal{C}_{\delta,\mathcal{R}}} < \tau_{Z\setminus Z_{\delta,\mathcal{R}}}$, and in the first case we apply the strong Markov property at time $\tau_{Z\setminus Z_{\delta,\mathcal{R}}}$. This yields

$$\begin{aligned}
&\text{l.h.s. of (6.65)}\\
&\quad\le q^{|x-z|}\\
&\qquad + \sum_{z_0\in Z\setminus Z_{\delta,\mathcal{R}}} q^{|x-z_0|}\mathbb{E}_{z_0}\exp\biggl\{\int_0^{\tau_z}\biggr\}\mathbb{1}\{\tau_z = \tau_{\mathcal{C}_{\delta,\mathcal{R}}} < \tau_{B^c}\}\\
&\qquad + \sum_{z_0\in \dot{Z}_{\delta,\mathcal{R}}} q^{|x-z_0|}\biggl[\sum_{z_1\in Z\setminus Z_{\delta,\mathcal{R}}} \mathbb{E}_{z_0}\exp\biggl\{\int_0^{\tau_{z_1}}\biggr\}\\
&\qquad\qquad \times \mathbb{1}\{\tau_{z_1} = \tau_{Z\setminus Z_{\delta,\mathcal{R}}} < \tau_{\mathcal{C}_{\delta,\mathcal{R}}} \wedge \tau_{B^c}\}\\
&\qquad\qquad \times \mathbb{E}_{z_1}\exp\biggl\{\int_0^{\tau_z}\biggr\}\mathbb{1}\{\tau_z = \tau_{\mathcal{C}_{\delta,\mathcal{R}}} < \tau_{B^c}\}\\
&\qquad\qquad + \mathbb{E}_{z_0}\exp\biggl\{\int_0^{\tau_{z_1}}\biggr\}\mathbb{1}\{\tau_{z_1}\tau_{\mathcal{C}_{\delta,\mathcal{R}}} < \tau_{Z\setminus Z_{\delta,\mathcal{R}}} \wedge \tau_{B^c}\}\biggr].
\end{aligned} \tag{6.66}$$

Now we apply Lemma 6.3, respectively Lemma 6.2, to the expectations on the right-hand side to obtain

$$\begin{aligned}
&\text{l.h.s. of (6.65)}\\
&\quad\le q^{|x-z|} + \sum_{z_0\in Z\setminus Z_{\delta,\mathcal{R}}} q^{|x-z_0|}q^{|z_0-z|/(256K^2)}\\
&\qquad + \sum_{z_0\in \dot{Z}_{\delta,\mathcal{R}}} q^{|x-z_0|}\biggl[\sum_{z_1\in Z\setminus Z_{\delta,\mathcal{R}}} q^{|z_0-z_1|/(8K)}q^{|z_1-z|/(256K^2)}\\
&\qquad\qquad + Mq^{|z_0-z|/(8K)}\biggr]
\end{aligned} \tag{6.67}$$

$$\leq C' q^{|x-z|/(512K^2)},$$

for some $C' > 0$, depending only on the quantities in (5.14), and on the constant $M > 1$ of Lemma 6.2(b). This ends the proof of the lemma. $\square$

With the help of Lemma 6.5, we may continue (6.61) [using also (6.63)] and (6.62) as follows:

$$\text{l.h.s. of (6.61)} \leq C'\left(1 + 2d\frac{|B_{\mathfrak{R}}|}{\mathfrak{g}^0}\right) \max_{z' \in \partial B_{\mathfrak{R}}(x_{k-1})} q^{|z'-x_{k-1}|/(512K^2)}, \tag{6.68}$$

$$\text{l.h.s. of (6.62)} \leq C'\left(1 + 2d\frac{|B_{\mathfrak{R}}|}{\mathfrak{g}^0}\right) \max_{z' \in \partial B_{\mathfrak{R}}(x_{k-1})} q^{|z'-x_k|/(512K^2)}. \tag{6.69}$$

On the right-hand side of (6.68) we estimate the last factor from above against $q^{\mathfrak{R}/(1024K^2)}$. The last factor on the right-hand side of (6.69) is estimated from above against $q^{\mathfrak{R}/(1024K^2)} q^{|z'-x_k|/(1024K^2)}$. Since $\mathfrak{R} \leq |z' - x_k|$, we may bound $|x_{k-1} - x_k| \leq \mathfrak{R} + |z' - x_k| \leq 2|z' - x_k|$ and, hence, $|z' - x_k|/(1024K^2) \geq |x_{k-1} - x_k|/(2048K^2)$.

Now we recall that $\mathfrak{R} = \mathfrak{R}_t = \log^2 t$ and assume $t$ so large that

$$q^{\mathfrak{R}/(1024K^2)} C'\left(1 + 2d\frac{|B_{\mathfrak{R}}|}{\mathfrak{g}^0}\right) < \frac{1}{2}. \tag{6.70}$$

Hence, we obtain that the right-hand side of (6.68) is not larger than $\frac{1}{2}$, and that the right-hand side of (6.69) is not larger than $\frac{1}{2} q^{|x_{k-1}-x_k|/(2048K^2)}$. Substituting these two bounds on the right-hand side of (6.60), we obtain that (6.56) holds.

Substituting (6.56) and (6.57) in (6.55), we obtain

$$\begin{aligned} \text{l.h.s. of (6.9)} &\leq \sum_{i=0}^{\infty} \sum_{x_1,\dots,x_i \in \mathcal{C}_{\delta,\mathcal{R}} \setminus \Gamma \colon |x_k - x_{k-1}| \geq \mathfrak{R}\, \forall k} \left[\prod_{k=1}^{i} q^{|x_k - x_{k-1}|/(2048K^2)}\right] \\ &\qquad\qquad \times q^{|x_i - y|/(2048K^2)} \\ &\leq q^{|z_0-y|/(3072K^2)} \sum_{i=0}^{\infty} \left(\sum_{|x| \geq \mathfrak{R}} q^{|x|/(6144K^2)}\right)^i \\ &\leq 2q^{|z_0-y|/(3072K^2)} \leq q^{|z_0-y|/(4096K^2)}, \end{aligned} \tag{6.71}$$

if $\mathfrak{R}$ is sufficiently large. This ends the proof of (6.9) and finishes the proof of Lemma 6.4.

6.6. *A corollary.* We add a technical result, which is a corollary of Proposition 6.1. For $y \in \Gamma$, we write $\dot{\lambda}_y^{(R)}$ for the principal eigenvalue $\lambda_{B_R(y)}(\dot{\xi})$ of $\Delta+\xi$ in the dotted set $B_R(y)\setminus\{y\}$. Furthermore, for $y\in\Gamma$, let $v_y^{(R)}:\mathbb{Z}^d\to[0,\infty)$ denote the principal $\ell^2(\mathbb{Z}^d)$-eigenfunction of $\Delta+\xi$ with zero boundary condition in $B_R(y)$, normalized by $v_y^{(R)}(y)=1$. Note that $v_y^{(R)}$ vanishes outside $B_R(y)$. The constant $c>0$ was introduced in Proposition 6.1.

Corollary 6.6. *For the choices of the parameters in Proposition 6.1, for any $y\in\Gamma$ and $x\in B_R(y)$,*

$$v_y(x)\le v_y^{(R)}(x)+q^{cR}\left(1+2d\frac{|B_R|}{\lambda_y^{(R)}-\dot\lambda_y^{(R)}}\right). \tag{6.72}$$

Proof. We go back to the representation of $v_y(x)$ in (6.2) and distinguish the two contributions from paths that stay in $B_R(y)$ by time $\tau_y$ and the remaining one. In the first term we estimate $\lambda_y\ge\lambda_y^{(R)}$, in the second term we use the strong Markov property at time $\tau_{B_R(y)^c}$. This yields

$$\begin{aligned}v_y(x)&\le\mathbb{E}_x\exp\left\{\int_0^{\tau_y}[\xi(X_s)-\lambda_y^{(R)}]\,ds\right\}\mathbb{1}\{\tau_y<\tau_{B_R(y)^c}\}\\&\quad+\mathbb{E}_x\exp\left\{\int_0^{\tau_{B_R(y)^c}}[\xi(X_s)-\lambda_y]\,ds\right\}\mathbb{1}\{\tau_{B_R(y)^c}<\tau_\Gamma\wedge\tau_{B^c}\}\\&\quad\times v_y(X_{\tau_{B_R(y)^c}}).\end{aligned}\tag{6.73}$$

The first expectation on the right-hand side is equal to $v_y^{(R)}(x)$, analogously to (6.2). Use (6.1) for the last factor on the right-hand side and note that $|X_{\tau_{B_R(y)^c}}-y|\ge R$. Then we obtain

$$\begin{aligned}v_y(x)&\le v_y^{(R)}(x)+q^{cR}\mathbb{E}_x\exp\left\{\int_0^{\tau_{B_R(y)^c}}[\xi(X_s)-\lambda_y^{(R)}]\,ds\right\}\\&\quad\times\mathbb{1}\{\tau_{B_R(y)^c}<\tau_\Gamma\wedge\tau_{B^c}\}.\end{aligned}\tag{6.74}$$

Now use Lemma 4.2 for $\gamma=\lambda_y^{(R)}$ and $A=B_R(y)\setminus\{y\}$ to finish the proof. □

**7. Localization and shape of potential.** In this section we prove that the set $\Gamma=\Gamma_{t\log^2 t}(\xi;\delta,\mathcal{R})$ defined in Section 5 satisfies all the assertions (2.8)–(2.11) of Proposition 2.2, provided that the parameters $\delta$ and $\mathcal{R}$ are chosen appropriately. As in the preceding sections, we consider the quantities in (5.14) as fixed. Recall that $q=2d/(2d+a/2)$.

Let $\varepsilon>0$ and $\eta>0$ be given and pick $r=r(\varrho,\varepsilon)\in\mathbb{N}_0$ as in (1.17) and $\mathfrak{R}_t=\log^2 t$ as in (5.9). Furthermore, let $\gamma>0$ and $R>0$ be given. Note that assertion (2.10) gets stronger if $R$ is picked larger or $\gamma>0$ smaller, respectively, and that the other three claims, (2.8), (2.9) and (2.11), do not depend on $\gamma$ nor on $R$. Hence, we are allowed to choose $R$ as large as we want and $\gamma$ as small as we want.

Let the random set $\Gamma = \Gamma_{t\log^2 t}$ be constructed (with $t$ replaced by $t\log^2 t$) as in Section 5.1, with parameters $\delta$ and $\mathcal{R}$. From Lemma 5.3 we already know that the conditions (2.9) and (2.10) are satisfied almost surely for sufficiently large $t$ if $R$ is large enough, $\delta$ small enough and $\mathcal{R}$ large enough. Hence, it only remains to show that (2.8) and (2.11) are satisfied, after possibly making the parameters $R, \delta, \mathcal{R}$ and $t$ even larger respectively smaller. This will be carried out in Sections 7.1, respectively 7.2, below. We assume that $R, \delta, \mathcal{R}$ and $t$ are chosen so large respectively, small that all the conclusions of Lemmas 5.3 and 5.4 and Proposition 6.1 hold. In addition, assume that $\delta$ is so small that $\delta < a/2$ and $2e^{-\delta/\varrho} - 1 > 1 - \eta$.

7.1. *Proof of* (2.8). Let $\alpha > 0$ be given. Because of (1.17), we may choose $\widetilde{\varepsilon} \in (0,1)$ (depending on $\alpha$ only) such that

$$(4\widetilde{\varepsilon} + \|w_\varrho\|_2)^2 \Bigg(4\widetilde{\varepsilon} + \sum_{x\in\mathbb{Z}^d\setminus B_{r'}} w_\varrho(x)\Bigg) < \varepsilon' + \alpha, \qquad \varepsilon' \in (\varepsilon, 1), \tag{7.1}$$

where $r' = r(\varrho, \varepsilon')$.

Now we assume $R$ so large (depending on $\widetilde{\varepsilon}$ only) that additionally all the following conditions are satisfied (cf. Lemma 3.3):

$$\|w_\varrho - w_\varrho^{(R)}\|_2 < \widetilde{\varepsilon} \quad \text{and} \quad \|w_\varrho - w_\varrho^{(R)}\|_1 < \widetilde{\varepsilon}, \tag{7.2}$$

$$|\lambda(V_\varrho) - \lambda_{B_R}(V_\varrho)| \vee |\lambda(\dot{V}_\varrho) - \lambda_{B_R}(\dot{V}_\varrho)| < \tfrac{1}{8}|\lambda(V_\varrho) - \lambda(\dot{V}_\varrho)|, \tag{7.3}$$

$$q^{cR}\Bigg(1 + 2d\frac{|B_R|}{1/2[\lambda(V_\varrho) - \lambda(\dot{V}_\varrho)]}\Bigg) < \frac{\widetilde{\varepsilon}}{|B_R|}, \tag{7.4}$$

$$\sum_{x\in\mathbb{Z}^d : |x|>R} q^{c|x|} < \widetilde{\varepsilon}^2, \tag{7.5}$$

where the constant $c > 0$ in (7.4) and (7.5) was introduced in Proposition 6.1.

Furthermore, we require that $\gamma$ is so small (depending on $R$ and $\widetilde{\varepsilon}$ only) that, for any $V : B_R \to \mathbb{R}$ satisfying $d_R(V, V_\varrho) < \gamma$, the following four bounds hold:

$$\begin{aligned} &|\lambda_{B_R}(V) - \lambda_{B_R}(V_\varrho)| < \tfrac{1}{8}|\lambda(V_\varrho) - \lambda(\dot{V}_\varrho)|, \\ &|\lambda_{B_R}(\dot{V}) - \lambda_{B_R}(\dot{V}_\varrho)| < \tfrac{1}{8}|\lambda(V_\varrho) - \lambda(\dot{V}_\varrho)| \end{aligned} \tag{7.6}$$

and

$$\|v - w_\varrho^{(R)}\|_{2,R} < \widetilde{\varepsilon} \quad \text{and} \quad \|v - w_\varrho^{(R)}\|_{1,R} < \widetilde{\varepsilon}, \tag{7.7}$$

where $v$ is the principal Dirichlet eigenfunction of $\Delta + V$ in $B_R$ satisfying $0 \in \arg\max(v)$ and $v(0) = 1$. (Here we have used that the principal eigenvalue and corresponding eigenfunction of $\Delta + V$ depend continuously on the potential $V$ on $B_R$.)

Lemma 5.3 implies that, with probability one, the set $\Gamma_{t\log^2 t}$ satisfies (2.10), if $t$ is sufficiently large. This means that

$$V_{y,t} = \xi(y + \cdot) - h_{t\log^2 t} \quad \text{satisfies} \quad d_R(V_{y,t}, V_\varrho) < \gamma, \qquad y \in \Gamma_{t\log^2 t}. \tag{7.8}$$

Now we apply Theorem 4.1 to see that the term on the left-hand side of (2.8) may be bounded as follows:

$$\begin{aligned} &\frac{1}{U(t)} \sum_{x \in B_{t\log^2 t} \setminus B_{r'}(\Gamma_{t\log^2 t})} u_3(t,x) \\ &\qquad \le \max_{y \in \Gamma_{t\log^2 t}} \Biggl[ \|v_y\|_2^2 \sum_{x \in B_{t\log^2 t} \setminus B_{r'}(\Gamma_{t\log^2 t})} v_y(x) \Biggr]. \end{aligned} \tag{7.9}$$

Here $v_y$ is the positive eigenfunction of $\Delta + \xi$ with zero boundary condition in $(B_{t\log^2 t} \setminus \Gamma_{t\log^2 t}) \cup \{y\}$. The corresponding eigenvalue is denoted by $\lambda_y$.

According to our choice of $\widetilde{\varepsilon}$ in (7.1), the only two things left to be proved are the following. With probability one, if $t$ is large enough, for any $y \in \Gamma_{t\log^2 t}$,

$$\sum_{x \in B_{t\log^2 t} \setminus B_{r'}(\Gamma_{t\log^2 t})} v_y(x) < 4\widetilde{\varepsilon} + \sum_{x \in \mathbb{Z}^d \setminus B_{r'}} w_\varrho(x), \tag{7.10}$$

$$\|v_y\|_2 < 4\widetilde{\varepsilon} + \|w_\varrho\|_2. \tag{7.11}$$

This is done as follows. We use the triangle inequality to estimate, for $x \in \mathbb{Z}^d$,

$$\begin{aligned} v_y(y+x) \le w_\varrho(x) &+ |w_\varrho(x) - w_\varrho^{(R)}(x)| + |w_\varrho^{(R)}(x) - v_y^{(R)}(y+x)| \\ &+ |v_y(y+x) - v_y^{(R)}(y+x)|. \end{aligned} \tag{7.12}$$

The last term is further estimated from above using Corollary 6.6 (using obvious notation):

$$0 \le v_y(y+x) - v_y^{(R)}(y+x) \le q^{cR}\biggl(1 + 2d\frac{|B_R|}{\lambda_y^{(R)} - \dot{\lambda}_y^{(R)}}\biggr), \qquad x \in B_R. \tag{7.13}$$

In order to estimate the denominator from below, use the triangle inequality to obtain

$$\begin{aligned} \lambda_y^{(R)} - \dot{\lambda}_y^{(R)} &= \lambda_{B_R(y)}(\xi - h_{t\log^2 t}) - \lambda_{B_R(y)}(\dot{\xi} - h_{t\log^2 t}) \\ &\ge -|\lambda_{B_R(y)}(\xi - h_{t\log^2 t}) - \lambda_{B_R}(V_\varrho)| \end{aligned}$$

$$
\begin{aligned}
&- |\lambda_{B_R}(V_\varrho) - \lambda(V_\varrho)| + \lambda(V_\varrho) - \lambda(\dot V_\varrho) \\
&- |\lambda(\dot V_\varrho) - \lambda_{B_R}(\dot V_\varrho)| \\
&- |\lambda_{B_R}(\dot V_\varrho) - \lambda_{B_R(y)}(\dot\xi - h_{t\log^2 t})|.
\end{aligned} \tag{7.14}
$$

(Here we used the notation $\dot\xi$ for the field $\xi$ dotted in $y$.) Use (7.3) to see that the second and the fourth of the terms on the r.h.s. are each not smaller than $-\frac{1}{8}(\lambda(V_\varrho) - \lambda(\dot V_\varrho))$. In order to see that also the first and the fifth of these terms are each not smaller than the same quantity, we recall that our requirement (7.6) applies to $v = V_{y,t}$ because of (7.8).

Hence, we obtain that $\lambda_y^{(R)} - \dot\lambda_y^{(R)} \geq \frac{1}{2}(\lambda(V_\varrho) - \lambda(\dot V_\varrho))$. Substituting this in (7.13), we obtain

$$
\begin{aligned}
v_y(y+x) - v_y^{(R)}(y+x) &\leq q^{cR}\left(1 + 2d\frac{|B_R|}{1/2[\lambda(V_\varrho) - \lambda(\dot V_\varrho)]}\right) \\
&\leq \frac{\widetilde\varepsilon}{|B_R|}, \qquad x \in B_R,
\end{aligned} \tag{7.15}
$$

according to (7.4).

For $x \notin B_R$, the last term in (7.12) is estimated by recalling that $v_y^{(R)}(y+x) = 0$ and by using Proposition 6.1 which implies that

$$
0 \leq v_y(y+x) \leq q^{c|x|}, \qquad y \in \Gamma_{t\log^2 t}, x \in \mathbb{Z}^d \setminus B_R. \tag{7.16}
$$

In order to prove (7.10), sum (7.12) over $x \in B_{t\log^2 t} \setminus B_{r'}(\Gamma_{t\log^2 t})$ to obtain

$$
\begin{aligned}
&\text{l.h.s. of (7.10)} \\
&\leq \sum_{x\in\mathbb{Z}^d\setminus B_{r'}} w_\varrho(x) + \|w_\varrho - w_\varrho^{(R)}\|_1 + \|w_\varrho^{(R)} - v_y^{(R)}(y+\cdot)\|_1 \\
&\quad + \sum_{x\in B_R} |v_y(y+x) - v_y^{(R)}(y+x)| + \sum_{x\in\mathbb{Z}^d\setminus B_R} v_y(y+x).
\end{aligned} \tag{7.17}
$$

The second and the third term on the r.h.s. do both not exceed $\widetilde\varepsilon$ by (7.2), respectively (7.7) [recall (7.8)]. The fourth term is not larger than $\widetilde\varepsilon$ by (7.15), and the last one as well by (7.16) and (7.5). This shows that (7.10) is indeed satisfied.

In order to prove (7.11), use the triangle inequality for $\|\cdot\|_2$ as in (7.12) and split the last norm into the contributions from $B_R$ and $\mathbb{Z}^d \setminus B_R$, and apply (7.15) and (7.16), respectively. This gives

$$
\begin{aligned}
\|v_y\|_2 &\leq \|w_\varrho\|_2 + \|w_\varrho - w_\varrho^{(R)}\|_2 + \|w_\varrho^{(R)} - v_y^{(R)}(y+\cdot)\|_2 \\
&\quad + \sqrt{|B_R|}\frac{\widetilde\varepsilon}{|B_R|} + \sqrt{\sum_{x\in\mathbb{Z}^d\setminus B_R} q^{2c|x|}}.
\end{aligned} \tag{7.18}
$$

Now the same arguments below (7.17) apply to show that (7.11) is satisfied. This ends the proof of (2.8).

7.2. *Proof of* (2.11). We borrow an important technical tool from [1], Lemma 4.6, which gives an estimate of the principal eigenvalue in large boxes in terms of the maximal principal eigenvalue in small subboxes:

LEMMA 7.1. *There is a constant $C>0$ such that, for any large centered box $B\subset\mathbb{Z}^d$, any potential $V:B\to\mathbb{R}\cup\{-\infty\}$ and any $n>0$,*

$$\lambda_B(V)\le\frac{C}{n^2}+\max_{x\in B}\lambda_{B_n(x)\cap B}(V). \tag{7.19}$$

Now we turn to the proof of (2.11) for $\varrho<\infty$ (the case $\varrho=\infty$ may be done in the same way and is even simpler). Assume that all the parameters are chosen as in the preceding subsection. We additionally require that the parameter $\delta$ is so small that $\varrho\log(1-e^{-(\chi(\varrho)+a)/\varrho})<-\delta/4$. And we assume that $\mathcal{R}$ is chosen so large that $4C/\mathcal{R}^2<\delta/8$, where $C$ is as in Lemma 7.1. Let us abbreviate $B=B_{t\log^2 t}$ and $h=h_{t\log^2 t}$ and so on.

By time reversion, we deduce from the definition of $u_2$ in (2.4) that

$$\sum_{x\in B}u_2(t,x)=\mathbb{E}_0\exp\biggl\{\int_0^t\xi(X_u)\,du\biggr\}\mathbb{1}\{\tau_{B^c}>t\}\mathbb{1}\{\tau_\Gamma>t\}. \tag{7.20}$$

From Schwarz' inequality and a Fourier expansion with respect to the eigenfunctions of $\Delta+\xi$ in $B\setminus\Gamma$ with zero boundary condition, also using Parseval's identity, one concludes that, for any $t>0$,

$$\sum_{x\in B}u_2(t,x)\le\sqrt{|B|}\|u_2(t,\cdot)\|_2=\sqrt{|B|}e^{t\lambda_{B\setminus\Gamma}(\xi)}\|\delta_0\|_2=\sqrt{|B|}e^{t\lambda_{B\setminus\Gamma}(\xi)}. \tag{7.21}$$

Recall from (1.12) that, almost surely, $U(t)\ge e^{th}e^{-t[\chi(\varrho)+o(1)]}$ as $t\to\infty$. Hence, for proving (2.11), it suffices to show that

$$\limsup_{t\to\infty}\lambda_{B\setminus\Gamma}(\xi-h)<-\chi(\varrho)\qquad\text{a.s.} \tag{7.22}$$

Define $\xi_\Gamma=\xi-\infty\mathbb{1}_\Gamma$, then it is clear that $\lambda_{B\setminus\Gamma}(\xi-h)=\lambda_B(\xi_\Gamma-h)$. We apply Lemma 7.1 with $B$ as above, $V=\xi_\Gamma-h$, and $n=\mathcal{R}/2$ to obtain that

$$\lambda_{B\setminus\Gamma}(\xi-h)=\lambda_B(\xi_\Gamma-h)\le\frac{4C}{\mathcal{R}^2}+\max_{x\in B}\lambda_{B_{\mathcal{R}/2}(x)}(\xi_\Gamma-h). \tag{7.23}$$

Recall that $4C/\mathcal{R}^2<\delta/8$. Hence, for proving (7.22), it suffices to show that

$$\limsup_{t\to\infty}\max_{x\in B}\lambda_{B_{\mathcal{R}/2}(x)}(\xi_\Gamma-h)\le-\chi(\varrho)-\frac{\delta}{4}\qquad\text{a.s.} \tag{7.24}$$

Let $x \in B$ and recall the definition of the set $Z$ in (5.1). We distinguish the three cases $B_{\mathcal{R}/2}(x) \cap Z = \varnothing$, $B_{\mathcal{R}/2}(x) \cap \Gamma \neq \varnothing$ and $B_{\mathcal{R}/2}(x) \cap \Gamma = \varnothing$, but $B_{\mathcal{R}/2}(x) \cap Z \neq \varnothing$.

In the first case we have $\xi_\Gamma - h \le -\chi(\varrho) - a$ in $B_{\mathcal{R}/2}(x)$, which implies that $\lambda_{B_{\mathcal{R}/2}(x)}(\xi_\Gamma - h) \le -\chi(\varrho) - a \le -\chi(\varrho) - \delta/4$.

For handling the second case, we recall the definition of the finite-space version of $\mathcal{L}$ in (3.3). According to Corollary 2.12 in [8], we have, almost surely,

$$\max_{x \in B} \mathcal{L}_{B_{K\mathcal{R}}(x)}(\xi - h) \le 1 + o(1) \qquad \text{as } t \to \infty, \tag{7.25}$$

and, consequently, uniformly in $x \in B$,

$$\mathcal{L}_{B_{\mathcal{R}/2}(x)}(\xi_\Gamma - h) \le \mathcal{L}_{B_{K\mathcal{R}}(x)}(\xi - h) - e^{-\min_Z(\xi - h)/\varrho} \le 1 + o(1) - e^{-(\chi(\varrho)+a)/\varrho}.$$

Using Lemma 3.1, we therefore obtain, almost surely,

$$\begin{aligned}\lambda_{B_{\mathcal{R}/2}(x)}(\xi_\Gamma - h) &\le -\chi(\varrho) + \varrho \log(1 + o(1) - e^{-(\chi(\varrho)+a)/\varrho}) \\ &\le -\chi(\varrho) - \delta/4 + o(1) \qquad \text{as } t \to \infty,\end{aligned}$$

uniformly in $x \in B$, according to our additional requirement on $\delta$.

In the third case, $B_{\mathcal{R}/2}(x)$ is contained in some large cluster $B'_{\mathcal{R}} = B_{\mathcal{R}}(Z_{\mathcal{R}}[z'])$, where $z'$ is its capital. If $\lambda_{B'_{\mathcal{R}}}(\xi) < h - \chi(\varrho) - \delta/4$, then we just estimate $\lambda_{B_{\mathcal{R}/2}(x)}(\xi_\Gamma - h) \le \lambda_{B'_{\mathcal{R}}}(\xi - h) \le -\chi(\varrho) - \delta/4$. Otherwise, the large cluster $B'_{\mathcal{R}}$ is $(\delta, \mathcal{R})$-optimal and the huge-cluster eigenvalue $\lambda_{B_{\mathfrak{R}}(z')}(\xi)$ lies above the spectral gap. Hence, by the definition of $\Gamma$, $z' \in \Gamma$. Now we argue as above to see that $\mathcal{L}_{B_{K\mathcal{R}}(z')}(\xi_\Gamma - h) \le 1 + o(1) - e^{-(\chi(\varrho)+a)/\varrho}$. Since $B_{\mathcal{R}/2}(x)$ is contained in $B_{K\mathcal{R}}(z')$, we may estimate $\lambda_{B_{\mathcal{R}/2}(x)}(\xi_\Gamma - h) \le \lambda_{B_{K\mathcal{R}}(z')}(\xi_\Gamma - h) \le -\chi(\varrho) - \delta/4 + o(1)$. This completes the proof of (7.24). $\square$

**8. Shape of the solution.** In this section we prove Proposition 2.3. We define a certain subset $\Gamma^*$ of the set $\Gamma$ defined in the proof of Proposition 2.2 such that (2.12) and (2.13) in Proposition 2.3 are satisfied.

Recall that the quantities in (5.14) are fixed. Let positive parameters $\varepsilon, \eta, \gamma, R$ be given. As in the proof of Proposition 2.2, we note that $\gamma$ may be made as small as we want and $R$ as large as we want. We shall pick the parameters $R, \gamma, \delta$ and $\mathcal{R}$ as in the proof of Proposition 2.2 in Section 7.1. Later we shall assume $\delta$ even smaller and $\mathcal{R}$ even larger in order that (2.12) and (2.13) can be deduced. Let $\mathfrak{R}_t = \log^2 t$, and let $\Gamma = \Gamma_{t \log^2 t}$ be defined as in Section 5.1.

Introduce

$$\Gamma^*_{t \log^2 t} = \Big\{ y \in \Gamma_{t \log^2 t} : u(t, y) \ge t^{-3\eta d} \max_{z \in \Gamma_{t \log^2 t}} u(t, z) \Big\}. \tag{8.1}$$

For any $s > 0$, choose $y^*(s) \in \Gamma$ such that

$$u(s, y^*(s)) = \max_{z \in \Gamma_{t \log^2 t}} u(s, z). \tag{8.2}$$

Certainly, $y^*(t)$ lies in $\Gamma^*_{t \log^2 t}$.

Recall from (7.11) that, for $t$ large,

$$\max_{z \in \Gamma_{t \log^2 t}} \|v_z\|_2^2 \le C, \tag{8.3}$$

where the constant $C > 1$ depends on the quantities in (5.14) only.

Let us now show that both (2.12) and (2.13) are satisfied.

8.1. *Proof of* (2.12). Fix $\varepsilon > 0$ and $\varepsilon' \in (0, \varepsilon)$. We use the split $u = u_1 + u_2 + u_3$ with $u_1$, $u_2$ and $u_3$ as in (2.3)–(2.5) with the set $\Gamma = \Gamma_{t \log^2 t}$ as defined in Section 5.1.

Apply (4.2) for $\Gamma = \Gamma_{t \log^2 t}$ and $w = u_3$ to estimate, for any $y \in \Gamma_{t \log^2 t} \setminus \Gamma^*_{t \log^2 t}$ and any $x \in B_{r(\varrho, \varepsilon')}(y)$,

$$u_3(t, x) \le \|v_y\|_2^2 v_y(x) u(t, y) + \sum_{z \in \Gamma_{t \log^2 t} \setminus \{y\}} \|v_z\|_2^2 v_z(x) u(t, y^*(t)). \tag{8.4}$$

Apply (8.3) to estimate the norm on the right-hand side of (8.4). We also estimate $v_y(x) \le C^{1/2} \le C$ on the right-hand side of (8.4). Use that $y$ is not in $\Gamma^*_{t \log^2 t}$ to estimate $u(t, y) \le t^{-3\eta d} u(t, y^*(t))$. By Proposition 6.1, $v_z(x)$ decays exponentially in the distance between $x$ and $z$. For $z \in \Gamma_{t \log^2 t} \setminus \{y\}$ and $x \in B_{r(\varrho, \varepsilon')}(y)$, since $\Gamma_{t \log^2 t}$ is spread out, $|x - z|$ is not smaller than a positive power of $t$. Hence, we may roughly estimate $v_z(x) \le t^{-4\eta d}$ for large $t$ in the second term on the right-hand side of (8.4). Furthermore, recall that $|\Gamma_{t \log^2 t}| \le t^{\eta d}$ as $t \to \infty$. Hence, for large $t$, we can continue (8.4) with

$$\begin{aligned} u_3(t, x) &\le C^2 t^{-3\eta d} u(t, y^*(t)) + \sum_{z \in \Gamma_{t \log^2 t} \setminus \{y\}} C t^{-4\eta d} u(t, y^*(t)) \\ &\le (C^2 t^{-3\eta d} + |\Gamma_{t \log^2 t}| C t^{-4\eta d}) u(t, y^*(t)) \\ &\le 2C^2 t^{-3\eta d} u(t, y^*(t)). \end{aligned} \tag{8.5}$$

This implies, for large $t$, the estimate

$$\begin{aligned} \sum_{x \in B_{r(\varrho, \varepsilon')}(\Gamma_{t \log^2 t} \setminus \Gamma^*_{t \log^2 t})} u_3(t, x) &\le |\Gamma_{t \log^2 t}| |B_{r(\varrho, \varepsilon')}| 2C^2 t^{-3\eta d} u(t, y^*(t)) \\ &\le t^{-\eta d} \sum_{x \in B_{r(\varrho, \varepsilon')}(\Gamma^*_{t \log^2 t})} u(t, x). \end{aligned} \tag{8.6}$$

Hence, we have derived (2.12).

8.2. *Proof of* (2.13). We shall prove (2.13) with the metric $d_R$ replaced by $\|\cdot\|_{\infty,R}$, the supremum norm on $B_R$. We introduce large auxiliary parameters $\mathtt{R}$ and $\mathtt{T}>0$. Fix $y\in\Gamma^*_{t\log^2 t}$. We write $\tau_{\mathtt{R}}$ for the exit time from $B_{\mathtt{R}}(y)$. Use the strong Markov property at time $\mathtt{T}\wedge\tau_{\mathtt{R}}$ to obtain that $u=u^{\mathrm{I}}+u^{\mathrm{II}}$, where, for any $x\in\mathbb{Z}^d$, and any $t>\mathtt{T}$,

$$u^{\mathrm{I}}(t,x)=\mathbb{E}_x\exp\biggl\{\int_0^{\mathtt{T}}\xi(X_s)\,ds\biggr\}u(t-\mathtt{T},X_{\mathtt{T}})\mathbb{1}\{\tau_{\mathtt{R}}>\mathtt{T}\},\tag{8.7}$$

$$u^{\mathrm{II}}(t,x)=\mathbb{E}_x\exp\biggl\{\int_0^{\tau_{\mathtt{R}}}\xi(X_s)\,ds\biggr\}u(t-\tau_{\mathtt{R}},X_{\tau_{\mathtt{R}}})\mathbb{1}\{\tau_{\mathtt{R}}\le\mathtt{T}\}.\tag{8.8}$$

In the following, we shall show that $u^{\mathrm{I}}$ gives the main contribution to $u$ in $B_R(y)$ and that $u^{\mathrm{II}}$ is negligible with respect to $u^{\mathrm{I}}$ in $B_R(y)$, provided that $\mathtt{T}$ and $\mathtt{R}$ are chosen large enough, $\delta$ small enough and afterward $\mathcal{R}$ large enough.

We write $(\cdot,\cdot)_{y,2,\mathtt{R}}$ and $\|\cdot\|_{y,2,\mathtt{R}}$ for the $\ell^2$-inner product and the corresponding norm in $B_{\mathtt{R}}(y)$, respectively. Denote the Dirichlet eigenvalue and eigenfunction pairs of $\Delta+\xi$ in $B_{\mathtt{R}}(y)$ by $(\lambda_k(\mathtt{R}),\mathrm{e}_k^{\mathtt{R}})$ for $k\in\mathbb{N}_0$, such that $\lambda_0(\mathtt{R})>\lambda_1(\mathtt{R})\ge\cdots$ and such that $(\mathrm{e}_k^{\mathtt{R}})_{k\in\mathbb{N}_0}$ is an orthonormal basis of $\ell^2(B_{\mathtt{R}}(y))$. A Fourier expansion yields that, for $x\in B_{\mathtt{R}}(y)$ and any $t>\mathtt{T}$,

$$\begin{aligned}u^{\mathrm{I}}(t,x)&=e^{\mathtt{T}\lambda_0(\mathtt{R})}(u(t-\mathtt{T},\cdot),\mathrm{e}_0^{\mathtt{R}}(\cdot))_{y,2,\mathtt{R}}\mathrm{e}_0^{\mathtt{R}}(x)\\&\quad+\sum_{k\in\mathbb{N}}e^{\mathtt{T}\lambda_k(\mathtt{R})}(u(t-\mathtt{T},\cdot),\mathrm{e}_k^{\mathtt{R}}(\cdot))_{y,2,\mathtt{R}}\mathrm{e}_k^{\mathtt{R}}(x).\end{aligned}\tag{8.9}$$

Applications of the Cauchy–Schwarz inequality and Parseval's identity yield that the second term is bounded in absolute value by

$$\begin{aligned}&e^{\mathtt{T}\lambda_1(\mathtt{R})}\sum_{k\in\mathbb{N}_0}|(u(t-\mathtt{T},\cdot),\mathrm{e}_k^{\mathtt{R}}(\cdot))_{y,2,\mathtt{R}}\mathrm{e}_k^{\mathtt{R}}(x)|\\&\quad\le e^{\mathtt{T}\lambda_1(\mathtt{R})}\sqrt{\sum_{k\in\mathbb{N}_0}(u(t-\mathtt{T},\cdot),\mathrm{e}_k^{\mathtt{R}}(\cdot))^2_{y,2,\mathtt{R}}}\sqrt{\sum_{k\in\mathbb{N}_0}(\mathrm{e}_k^{\mathtt{R}}(\cdot),\delta_x)^2_{y,2,\mathtt{R}}}\\&\quad=e^{\mathtt{T}\lambda_1(\mathtt{R})}\|u(t-\mathtt{T},\cdot)\|_{y,2,\mathtt{R}}.\end{aligned}\tag{8.10}$$

Hence, we may summarize (8.9) by writing

$$u^{\mathrm{I}}(t,y+x)=F(t,\mathtt{R},\mathtt{T},y)w_\varrho(x)[1+f(t,\mathtt{R},\mathtt{T},y,x)],\qquad x\in B_{\mathtt{R}},\tag{8.11}$$

where

$$\begin{aligned}F(t,\mathtt{R},\mathtt{T},y)&=e^{\mathtt{T}\lambda_0(\mathtt{R})}(u(t-\mathtt{T},\cdot),\mathrm{e}_0^{\mathtt{R}}(\cdot))_{y,2,\mathtt{R}}\frac{1}{\|w_\varrho\|_2},\\|f(t,\mathtt{R},\mathtt{T},y,x)|&\le\frac{|w_\varrho(x)-\mathrm{e}_0^{\mathtt{R}}(y+x)\|w_\varrho\|_2|}{w_\varrho(x)}\end{aligned}\tag{8.12}$$

$$+\frac{e^{\mathtt{T}(\lambda_1(\mathtt{R})-\lambda_0(\mathtt{R}))}}{w_\varrho(x)/\|w_\varrho\|_2}\frac{\|u(t-\mathtt{T},\cdot)\|_{y,2,\mathtt{R}}}{(u(t-\mathtt{T},\cdot),\mathrm{e}_0^{\mathtt{R}}(\cdot))_{y,2,\mathtt{R}}}. \tag{8.13}$$

If $R<\mathtt{R}$, this implies that (recall that $w_\varrho$ is maximal at 0 with value 1)

$$\begin{aligned}\left\|\frac{u(t,y+\cdot)}{u(t,y)}-w_\varrho(\cdot)\right\|_{\infty,R} &\le \max_{x\in B_R}\left|\frac{f(t,\mathtt{R},\mathtt{T},y,x)-f(t,\mathtt{R},\mathtt{T},y,0)}{1+f(t,\mathtt{R},\mathtt{T},y,0)}\right| \\ &\quad+\max_{x\in B_R}\frac{u^{\mathrm{II}}(t,y+x)}{u(t,y)}.\end{aligned} \tag{8.14}$$

Hence, in order to finish the proof of (2.13), it suffices to show that

$$\limsup_{t\to\infty}\max_{y\in\Gamma^*_{t\log^2 t}}\max_{x\in B_R}|f(t,\mathtt{R},\mathtt{T},y,x)|\le\frac{\gamma}{8}, \tag{8.15}$$

$$\limsup_{t\to\infty}\max_{y\in\Gamma^*_{t\log^2 t}}\max_{x\in B_R}\frac{u^{\mathrm{II}}(t,y+x)}{u(t,y)}\le\frac{\gamma}{2}, \tag{8.16}$$

for $\mathtt{T}$ and $\mathtt{R}$ sufficiently large, if $\delta$ is small enough and $\mathcal{R}$ large enough.

Let us first determine how to choose the parameters $\mathtt{T}$ and $\mathtt{R}$. Recall the optimal potential shape $V_\varrho$ in Assumption (M), and let $\lambda_0>\lambda_1$ denote the two leading eigenvalues of $\Delta+V_\varrho$ in $\mathbb{Z}^d$. Furthermore, recall that $w_\varrho$ is the positive eigenfunction to the eigenvalue $\lambda_0$ satisfying $w_\varrho(0)=1$. We first pick $\mathtt{T}$ so large that

$$e^{-\mathtt{T}(\lambda_0-\lambda_1)/2}12C^3<\frac{\gamma}{32}\inf_{x\in B_R}w_\varrho(x), \tag{8.17}$$

with $C>1$ as in (8.3). It follows from [8], Lemma 2.5(a), that $\mathbb{P}_x(\tau_{\mathtt{R}}\le\mathtt{T})\le 2^{d+1}e^{1/2\mathtt{R}-1/2\mathtt{R}\log(\mathtt{R}/(2d\mathtt{T}))}$, for any $x\in B_{\mathtt{R}/2}(y)$. Hence, it is possible to choose $\mathtt{R}>2R$ so large that $\mathtt{R}>r(\varrho,\varepsilon)$ and

$$4C^3e^{2d\mathtt{T}}|B_{\mathtt{R}/2}|\max_{x\in B_{\mathtt{R}/2}(y)}\mathbb{P}_x(\tau_{\mathtt{R}}\le\mathtt{T})<\frac{\gamma}{64}, \tag{8.18}$$

$$C^3e^{2d\mathtt{T}}\sum_{x\in\mathbb{Z}^d\setminus B_{\mathtt{R}/2}}w_\varrho(x)<\frac{1}{64}. \tag{8.19}$$

Recall the notation from the beginning of Section 6.6. Observe that $v_y^{\mathtt{R}}(\cdot)=\mathrm{e}_0^{\mathtt{R}}(\cdot)/\mathrm{e}_0^{\mathtt{R}}(0)$ and compare (7.2) and (7.7). Replacing $R$ in Section 7.1 by $\mathtt{R}$, we see that it is possible to choose the parameters $\delta$ and $\mathcal{R}$ so small respectively large that, almost surely, if $t$ is large enough, for any $\widetilde{y}\in\Gamma_{t\log^2 t}$, the potential $\xi(\widetilde{y}+\cdot)-h$ is sufficiently close to $V_\varrho$ in $B_{\mathtt{R}}$ such that

$$\lambda_0(\mathtt{R})-\lambda_1(\mathtt{R})>\tfrac{1}{2}(\lambda_0-\lambda_1), \tag{8.20}$$

$$\max_{x\in B_{\mathtt{R}}}|w_\varrho(x)-\mathrm{e}_0^{\mathtt{R}}(y+x)\|w_\varrho\|_2|<\frac{\gamma}{16}\inf_{x\in B_{\mathtt{R}}}w_\varrho(x), \tag{8.21}$$

$$\max_{x\in B_{\mathtt{R}}}|v_{\widetilde{y}}(\widetilde{y}+x)-w_\varrho(x)|\le\tfrac12\inf_{x\in B_{\mathtt{R}}}w_\varrho(x). \tag{8.22}$$

[For (8.22), see (7.12) and the arguments below (7.12).] In particular, we may estimate $v_{\widetilde{y}}(\widetilde{y}+x)\le 2$ since $w_\varrho$ is maximal at zero with value one.

As a preparation for the proof of (8.15)–(8.16), we now present a number of lemmas.

Lemma 8.1. *Almost surely, for $t$ large, for any $y\in\Gamma^*_{t\log^2 t}$, $x\in B_{\mathtt{R}}(y)$ and $s\in[0,\mathtt{T}]$,*

$$u_3(t-s,x)\le 2C^2e^{-\lambda_y s}u_3(t,y)v_y(x), \tag{8.23}$$

*where $C$ is as in (8.3).*

Proof. Fix $y\in\Gamma^*_{t\log^2 t}$. Use (4.2) for $t-s$ instead of $t$ and recall (8.3) to obtain, for $x\in B_{\mathtt{R}}(y)$ and $t$ large,

$$\begin{aligned} u_3(t-s,x)&\le\|v_y\|_2^2u_3(t-s,y)v_y(x)\\ &\qquad+\sum_{\widetilde{y}\in\Gamma_{t\log^2 t}\setminus\{y\}}\|v_{\widetilde{y}}\|_2^2u_3(t-s,\widetilde{y})v_{\widetilde{y}}(x)\\ &\le Cu_3(t-s,y)v_y(x)\\ &\qquad+|\Gamma_{t\log^2 t}|C\max_{\widetilde{y}\in\Gamma_{t\log^2 t}}u_3(t-s,\widetilde{y})e^{-\sqrt{t}}, \end{aligned} \tag{8.24}$$

where we estimated, according to Proposition 6.1, $v_{\widetilde{y}}(x)\le q^{c|x-\widetilde{y}|}$ and used that $|x-\widetilde{y}|\ge|y-\widetilde{y}|-\mathtt{R}$ is not smaller than a power of $t$ that is close to one. Use Lemma 4.3 and again (8.3) to bound

$$u_3(t-s,\widetilde{y})\le Ce^{-\lambda_{\widetilde{y}}s}u_3(t,\widetilde{y}),\qquad \widetilde{y}\in\Gamma_{t\log^2 t}. \tag{8.25}$$

Because $y\in\Gamma^*_{t\log^2 t}$, we have

$$u_3(t,\widetilde{y})\le u(t,\widetilde{y})\le u(t,y^*(t))\le t^{3\eta d}u(t,y),\qquad \widetilde{y}\in\Gamma_{t\log^2 t}. \tag{8.26}$$

Substituting (8.26) in (8.25) and (8.25) in the second term of (8.24), and recalling (5.19) and that $|\Gamma_{t\log^2 t}|\le t^{\eta d}$, we obtain, for $t$ large,

$$u_3(t-s,x)\le Cu_3(t-s,y)v_y(x)+t^{5\eta d}e^{-\sqrt{t}}e^{(-\lambda_y+\chi+\delta/2)s}u(t,y). \tag{8.27}$$

Now use (8.25) for $\widetilde{y}=y$ to estimate the first term on the right-hand side of (8.27). Bound the second, for $t$ large, as follows: $u(t,y)=u_1(t,y)+u_3(t,y)\le$

$Cu_3(t,y)$ [recall (2.7)]. This gives, uniformly in $y \in \Gamma_{t\log^2 t}$ and $x \in B_{\mathtt{R}}(y)$,

$$u_3(t-s,x) \le e^{-\lambda_y s} u_3(t,y) v_y(x)\left[C^2 + o(1)\frac{1}{\inf_{B_{\mathtt{R}}(y)} v_y}\right]. \tag{8.28}$$

Now use (8.22) to obtain (8.23). □

Lemma 8.2. *As $t \to \infty$,*

$$\min_{y \in \Gamma^*_{t\log^2 t}} u(t,y) \ge \exp\{t(h_{t\log^2 t} - \chi(\varrho) - o(1))\}. \tag{8.29}$$

Proof. Recall that $\sum_{z\in\mathbb{Z}^d} u(t,z) = e^{t(h_{t\log^2 t} - \chi(\varrho) + o(1))}$ as $t \to \infty$ almost surely, and put $r = r(\varrho, \varepsilon)$. From (2.12), we therefore have that

$$e^{t(h_{t\log^2 t} - \chi(\varrho) - o(1))} = \sum_{\widetilde{y} \in \Gamma^*_{t\log^2 t}} \sum_{x \in B_r(\widetilde{y})} u_3(t,x). \tag{8.30}$$

Use Lemma 8.1 for $s = 0$ and recall (8.3) to see that $u_3(t,x) \le 2C^3 u_3(t,\widetilde{y})$ for $x \in B_r(\widetilde{y})$ and $\widetilde{y} \in \Gamma^*_{t\log^2 t}$. In addition, use (8.26) in (8.30) to estimate, for any $y \in \Gamma^*_{t\log^2 t}$ and $t$ large,

$$e^{t(h_{t\log^2 t} - \chi(\varrho) - o(1))} \le 2C^3 t^{3\eta d} |\Gamma^*_{t\log^2 t}|\,|B_r| u(t,y) \le t^{5\eta d} u(t,y). \tag{8.31}$$

This implies the assertion, noting that the term $t^{5\eta d}$ may be absorbed in the term $e^{o(t)}$ on the left-hand side. □

Recall from Section 7.2 that, if $\mathcal{R}$ is large enough (depending on $\delta$ only), for $t$ large, and any $s \in [0, \mathtt{T}]$,

$$\begin{aligned}\sum_{x \in B_{t\log^2 t}} u_2(t-s,x) &\le \sqrt{|B_{t\log^2 t}|}\, e^{-(t-s)(h_{t\log^2 t} - \chi(\varrho) + o(1) - \delta/8)} \\ &\le e^{-(t-s)(h_{t\log^2 t} - \chi(\varrho) - \delta/16)}.\end{aligned} \tag{8.32}$$

Lemma 8.3. *Almost surely, for $t$ large, the following holds. For any $y \in \Gamma^*_{t\log^2 t}$, any $x \in B_{\mathtt{R}}$,*

$$u^{\mathrm{II}}(t, y+x) \le u(t,y)[o(1) + 2C^3 e^{2d\mathtt{T}} \mathbb{P}_x(\tau_{\mathtt{R}} \le \mathtt{T})], \tag{8.33}$$

*where $o(1)$ is uniform in $y$ and $x$.*

PROOF. We again split $u = u_1 + u_2 + u_3$ on the right-hand side of (8.8), which yields a sum of three terms. The first one is estimated as follows. On $\{\tau_{\mathtt{R}} \le \mathtt{T}\}$ we estimate $u_1(t-\tau_{\mathtt{R}}, X_{\tau_{\mathtt{R}}}) \le \sum_{z\in\mathbb{Z}^d} u_1(t-\tau_{\mathtt{R}}, z) \le o(1)$ [see (2.7)]. Furthermore, we bound $\xi(X_s) \le h_{t\log^2 t}$ in the exponent. Thus, the first summand is not larger than $o(1)e^{\mathtt{T}h_{t\log^2 t}}$, which is $o(1)u(t,y)$, in accordance with Lemma 8.2.

In the second summand we bound, with the help of (8.32), $u_2(t-\tau_{\mathtt{R}}, X_{\tau_{\mathtt{R}}}) \le e^{(t-\tau_{\mathtt{R}})(h_{t\log^2 t}-\chi(\varrho)-\delta/16)}$ and again bound $\xi(X_s) \le h_{t\log^2 t}$ in the exponent. Thus, the second summand is, for $t$ large, not larger than

$$\begin{aligned}
&\mathbb{E}_x e^{\int_0^{\tau_{\mathtt{R}}} h_{t\log^2 t}\,ds} e^{(t-\tau_{\mathtt{R}})(h_{t\log^2 t}-\chi(\varrho)-\delta/16)} \mathbb{1}\{\tau_{\mathtt{R}} \le \mathtt{T}\}\\
&\quad\le e^{t(h_{t\log^2 t}-\chi(\varrho)-\delta/16)} e^{\mathtt{T}(\chi(\varrho)+\delta/16)}\\
&\quad\le u(t,y)e^{-t\delta/32},
\end{aligned} \tag{8.34}$$

where we again used Lemma 8.2.

In the third and last summand, we use Lemma 8.1 for bounding the $u_3$-term. Hence, the third summand is not larger than

$$\begin{aligned}
&\mathbb{E}_x e^{\int_0^{\tau_{\mathtt{R}}} \xi(X_s)\,ds} 2C^3 u(t,y) e^{-\lambda_y \tau_{\mathtt{R}}} v_y(X_{\tau_{\mathtt{R}}}) \mathbb{1}\{\tau_{\mathtt{R}} \le \mathtt{T}\}\\
&\quad\le 2C^3 u(t,y) \mathbb{E}_x e^{\int_0^{\tau_{\mathtt{R}}} [\xi(X_s)-\lambda_0(\mathtt{R})]\,ds} \mathbb{1}\{\tau_{\mathtt{R}} \le \mathtt{T}\},
\end{aligned} \tag{8.35}$$

where we used (8.3) and that $\lambda_y \ge \lambda_0(\mathtt{R})$. Now estimate $\xi(X_s) - \lambda_0(\mathtt{R}) \le 2d$ in the exponent and summarize to obtain the upper bound $2C^3 e^{2d\mathtt{T}} \mathbb{P}_x(\tau_{\mathtt{R}} \le \mathtt{T})$. Collecting the upper bounds for the other two terms, we arrive at the assertion. □

Let us now finish the proof of (8.15)–(8.16). The estimate in (8.16) is now immediate from Lemma 8.3 in combination with (8.18) (recall that $\mathtt{R}/2 > R$).

We turn to the proof of (8.15). We shall show the following bounds for $t$ large:

$$\|u(t-\mathtt{T},\cdot)\|_{y,2,\mathtt{R}} \le 3C^3 e^{-\lambda_0(\mathtt{R})\mathtt{T}} u(t,y), \tag{8.36}$$

$$(u(t-\mathtt{T},\cdot), \mathrm{e}_0^{\mathtt{R}})_{y,2,\mathtt{R}} \ge \frac{1}{4\|w_\varrho\|_2} e^{-\lambda_0(\mathtt{R})\mathtt{T}} u(t,y). \tag{8.37}$$

Combining (8.36) and (8.37), we see that the last quotient on the right-hand side of (8.13) is bounded by $12C^3\|w_\varrho\|_2$ for $t$ sufficiently large. Using (8.20)–(8.21) and (8.17), we see that (8.15) is indeed satisfied. This ends the proof of (8.15), subject to (8.36) and (8.37).

Let us now derive (8.36). Apply the triangle inequality to bound

$$\begin{aligned}
\|u(t-\mathtt{T},\cdot)\|_{y,2,\mathtt{R}} &\le \|u_1(t-\mathtt{T},\cdot)\|_{y,2,\mathtt{R}} + \|u_2(t-\mathtt{T},\cdot)\|_{y,2,\mathtt{R}}\\
&\quad + \|u_3(t-\mathtt{T},\cdot)\|_{y,2,\mathtt{R}}.
\end{aligned} \tag{8.38}$$

Recall (2.7). Since, in particular, $\liminf_{t\to\infty}\min_{y\in\Gamma^*_{t\log^2 t}} u(t,y)=\infty$ by Lemma 8.2, it is clear that, as $t\to\infty$,

$$\|u_1(t-\mathtt{T},\cdot)\|_{y,2,\mathtt{R}} \le \sum_{x\in\mathbb{Z}^d} u_1(t-\mathtt{T},x) \le o(1) \le o(1)e^{-\lambda_0(\mathtt{R})\mathtt{T}}u(t,y), \tag{8.39}$$

where $o(1)$ is uniform in $y\in\Gamma^*_{t\log^2 t}$.

With the help of (8.32), we estimate the second summand on the right-hand side of (8.38) as follows. For $t$ large and any $s\in[0,\mathtt{T}]$,

$$\begin{aligned} \|u_2(t-s,\cdot)\|_{y,2,\mathtt{R}} &\le \sum_{x\in B_{t\log^2 t}} u_2(t-s,x) \\ &\le e^{(t-s)(h_{t\log^2 t}-\chi-\delta/16)} \\ &\le o(1)e^{-\lambda_0(\mathtt{R})\mathtt{T}}u(t,y), \end{aligned} \tag{8.40}$$

where we used that $\lambda_0(\mathtt{R})\le h_{t\log^2 t}\le o(t)$, and $o(1)$ is uniform in $y\in\Gamma^*_{t\log^2 t}$, according to Lemma 8.2.

In order to bound the third term on the right-hand of (8.38), use Lemma 8.1 to estimate

$$\|u_3(t-\mathtt{T},\cdot)\|_{y,2,\mathtt{R}} \le 2C^2e^{-\lambda_y\mathtt{T}}u_3(t,y)\|v_y\|_2 \le 2C^3e^{-\lambda_0(\mathtt{R})\mathtt{T}}u(t,y), \tag{8.41}$$

where we recall (8.3). Now substitute (8.39), (8.40) and (8.41) in (8.38) to see that (8.36) holds.

Now we prove (8.37). Recall $u=u^{\mathrm{I}}+u^{\mathrm{II}}$ [see (8.7)–(8.8)] and the Fourier expansion in (8.9) for $u^{\mathrm{I}}$. Multiply (8.9) with $\mathrm{e}_0^{\mathtt{R}}$ and sum over $x\in B_{\mathtt{R}}(y)$, and use that $(\mathrm{e}_k^{\mathtt{R}})_{k\in\mathbb{N}_0}$ is an orthonormal basis, to obtain

$$\begin{aligned} (u(t,\cdot),\mathrm{e}_0^{\mathtt{R}})_{y,2,\mathtt{R}} &= (u^{\mathrm{I}}(t,\cdot),\mathrm{e}_0^{\mathtt{R}})_{y,2,\mathtt{R}} + (u^{\mathrm{II}}(t,\cdot),\mathrm{e}_0^{\mathtt{R}})_{y,2,\mathtt{R}} \\ &= e^{\lambda_0(\mathtt{R})\mathtt{T}}(u(t-\mathtt{T},\cdot),\mathrm{e}_0^{\mathtt{R}})_{y,2,\mathtt{R}} + (u^{\mathrm{II}}(t,\cdot),\mathrm{e}_0^{\mathtt{R}})_{y,2,\mathtt{R}}. \end{aligned} \tag{8.42}$$

Hence,

$$\begin{aligned} (u(t-\mathtt{T},\cdot),\mathrm{e}_0^{\mathtt{R}})_{y,2,\mathtt{R}} &= e^{-\lambda_0(\mathtt{R})\mathtt{T}}[(u(t,\cdot),\mathrm{e}_0^{\mathtt{R}})_{y,2,\mathtt{R}} - (u^{\mathrm{II}}(t,\cdot),\mathrm{e}_0^{\mathtt{R}})_{y,2,\mathtt{R}}] \\ &\ge e^{-\lambda_0(\mathtt{R})\mathtt{T}}[u(t,y)\mathrm{e}_0^{\mathtt{R}}(y) - (u^{\mathrm{II}}(t,\cdot),\mathrm{e}_0^{\mathtt{R}})_{y,2,\mathtt{R}}]. \end{aligned} \tag{8.43}$$

Recall from (8.21) that $\mathrm{e}_0^{\mathtt{R}}(y)\ge 1/(2\|w_\varrho\|_2)$. We use now Lemma 8.3 in order to show that the last term in (8.43) is not larger than $u(t,y)/(4\|w_\varrho\|_2)$. Indeed, use (8.33) to obtain

$$\begin{aligned} &(u^{\mathrm{II}}(t,\cdot),\mathrm{e}_0^{\mathtt{R}})_{y,2,\mathtt{R}} \\ &\quad\le u(t,y)2C^3e^{2d\mathtt{T}}\sum_{x\in B_{\mathtt{R}}(y)}\mathbb{P}_x(\tau_{\mathtt{R}}\le\mathtt{T})\mathrm{e}_0^{\mathtt{R}}(x) + o(1)u(t,y). \end{aligned} \tag{8.44}$$

Now split the sums in the sum over $x \in B_{\mathtt{R}/2}(y)$ and over $x \in B_{\mathtt{R}}(y) \setminus B_{\mathtt{R}/2}(y)$. In the first sum estimate $\mathrm{e}_0^{\mathtt{R}}(x) \le 2/\|w_\varrho\|_2$ [recall (8.21)] and use (8.18), to see that this sum is not larger than $u(t,y)/(16\|w_\varrho\|_2)$. In the second sum estimate $\mathbb{P}_x(\tau_{\mathtt{R}} \le \mathtt{T}) \le 1$ and $\mathrm{e}_0^{\mathtt{R}}(x) \le 2w_\varrho(x)/\|w_\varrho\|_2$ [see (8.21)] to obtain the upper bound

$$
\begin{aligned}
& u(t,y) 2C^3 e^{2d\mathtt{T}} \sum_{x \in B_{\mathtt{R}}(y) \setminus B_{\mathtt{R}/2}(y)} \mathbb{P}_x(\tau_{\mathtt{R}} \le \mathtt{T}) \mathrm{e}_0^{\mathtt{R}}(x) \\
& \quad \le u(t,y) \frac{4}{\|w_\varrho\|_2} C^3 e^{2d\mathtt{T}} \sum_{x \in \mathbb{Z}^d \setminus B_{\mathtt{R}/2}} w_\varrho(x) \qquad (8.45) \\
& \quad < \frac{1}{16\|w_\varrho\|_2} u(t,y),
\end{aligned}
$$

according to (8.19). Hence, we have obtained that

$$
(u^{\mathrm{II}}(t,\cdot), \mathrm{e}_0^{\mathtt{R}})_{y,2,\mathtt{R}} \le \frac{1}{4\|w_\varrho\|_2} u(t,y). \tag{8.46}
$$

Using (8.46) in (8.43) and recalling that $\mathrm{e}_0^{\mathtt{R}}(y) \ge 1/(2\|w_\varrho\|_2)$, we see that (8.37) is satisfied.

J. Gärtner
Technische Universität Berlin
Institut für Mathematik
MA7-5, Strasse des 17. Juni 136
D-10623 Berlin
Germany

W. König
Universität Leipzig
Mathematisches Institut
Augustusplatz 10/11
D-04109 Leipzig
Germany
E-mail: koenig@math.uni-leipzig.de

S. Molchanov
Department of Mathematics
University of North Carolina
Charlotte, North Carolina 28223
USA